\documentclass[a4paper]{article}

\usepackage[all]{xy}
\usepackage{latexsym,amssymb,xy}
\usepackage{amsopn}
\usepackage{amsmath}
\usepackage{tensor}
\usepackage{mathrsfs}
\usepackage{makeidx}

\usepackage[shortcuts]{extdash} 

\usepackage[pdftex]{graphicx}
\usepackage[pdftex]{color}

\definecolor{BlueViolet}{rgb}{0.26,0.12,0.66}

\def\dem{\noindent{\bf Proof.} }

\def\f{\mbox{\sf f}} 

\def\ham{\mbox{$\mathfrak{X}_{\textrm{ham}}^{\bullet}(M)$}}
\def\ins{\raisebox{\depth}{\rotatebox{180}{$\neg$}}}

\def\Linf{\mbox{$L_{\infty}$}\=/algebra }
\def\Linfm{\mbox{$L_{\infty}$}\=/morphism }
\def\Linfms{\mbox{$L_{\infty}$}\=/morphisms }
\def\Linfs{\mbox{$L_{\infty}$}\=/algebras }
\def\Linfb{\mbox{\boldmath ${L_{\infty}}$}\=/algebra }

\def\Linfsb{\mbox{\boldmath ${L_{\infty}}$}\=/algebras }

\def\NA{\mathbb N}

\def\nth{\mbox{$n\mbox{th }$}}

\def\qed{\hfill\ensuremath{\blacksquare}}

\DeclareMathOperator{\dec}{dec}

\newcommand{\beq}{\begin{equation}}
\newcommand{\eeq}{\end{equation}}
\newcommand{\beqa}{\begin{eqnarray}}
\newcommand{\eeqa}{\end{eqnarray}}

\newtheorem{tm}{Theorem}[section]
\newtheorem{rk}[tm]{Remark}
\newtheorem{pp}[tm]{Proposition}
\newtheorem{lm}[tm]{Lemma}
\newtheorem{cl}[tm]{Corollary}
\newtheorem{df}[tm]{Definition}
\newtheorem{dfpp}[tm]{Definition/Proposition}

\usepackage{tikz}
\usetikzlibrary{matrix,arrows,decorations.pathmorphing}

\author{\Large{N\'estor Le\'on Delgado}}
\title{\textcolor{BlueViolet}{ \huge{Multisymplectic Structures and Higher Momentum Maps}}}
\date{\today}

\makeindex

\usepackage[pdftex]{hyperref}

\begin{document}

\maketitle

\tableofcontents


\color{BlueViolet}
\section*{Acknowledgments}
\color{black}

This document is a reorganization of the results on the Master Thesis: {\it Multisymplectic Structures and Higher Momentum Maps} written by N\'estor Le\'on Delgado under the supervision of Dr. Christian Blohmann at the University of Bonn in $2014$. I express my gratitude to Christian Blohmann for his time and his ideas which make a huge impact in the work here presented.

The reason behind this reorganization and upload into the arXiv is so that it is available for future references. I would like to thank Marco Zambon for encouraging me to put this document together and to Camille Laurent-Gengoux for having repeatedly push this from happening. The fact that him and his collaborators used some of the results here presented in their paper \cite{AL}, was key in me deciding to make these results available for the general public.

Back when I was writing the Master Thesis report, I benefited greatly from the Higher Differential Geometry Seminar at the Max Planck Institute for Mathematics.  Special thanks to those speakers at the seminar who shared ideas with me: Marco Zambon, Benjamin Cooper and Mohammad Azimi. I also thank the organizers of that seminar: Christian Blohmann and David Carchedi for making the meetings possible.

\color{BlueViolet}
\section*{Introduction}
\color{black}

There are three main results in this document. 

\begin{itemize}
	\item The first one, Theorem \ref{main}, constructs an \mbox{$L_{\infty}$}\=/algebra structure for every Gerstenhaber algebra which involves non\=/vanishing higher brackets.
	\item The second one, Theorem \ref{mainTransfer}, uses a cochain null\=/homotopy from the Chevalley\=/Eilenberg complex $C(\mathfrak{g})$ associated to a graded Lie algebra $\mathfrak{g}$ to any cochain complex $D$ to induce (transfer) the \Linf structure from $\mathfrak{g}$ to $Q(D)$ a resolution of $D$. 
	\item The third one, Corollary \ref{central}, establishes the equivalence between cochain\=/homotopies to $Q$ and \Linf lifts satisfying certain properties.
\end{itemize}

All the three results are motivated from the study of multi\=/symplectic structures. In particular:

\begin{itemize}
	\item The Gerstenhaber algebra we will have in mind is the Gertenhaber algebra of multi\=/vector fields in a smooth manifold. If the manifold is pre\=/multi\=/symplectic, symplect and hamiltonian multi\=/vector fields are sub\=/\Linfs of the \mbox{$L_{\infty}$}\=/algebra given by Theorem \ref{main}.
	\item The construction of Theorem \ref{mainTransfer} is applied to $L$ multi\=/vector fields in a pre\=/multi\=/symplectic manifold and $D$ the complex of differential forms on the manifold. The trasferred structre defines an \Linf on hamiltonian forms that generalizes the work of Rogers in \cite{R}.
	\item Applying Corollary \ref{central} to that  example has a geometric interpretation: the maps can be understood as comomentum maps for a hamiltonian action.
\end{itemize}


\color{BlueViolet}
\section{Gerstenhaber Algebras and \mbox{\boldmath ${L_{\infty}}$}\=/algebras}
\color{black}

In this section we present a brief summary of the conventions we will be using when talking about graded vector spaces. We put some emphasis on the decalage isomorphism since we present an interpretation in terms of category theory of the isomorphism due to C. Blohmann (private communication).

Later in the Chapter we include the definitions of Gerstenhaber algebras and \Linf , necessary to understand the main theorem in this first part of the document. 

All the concepts here defined are common in the literature, we refer to the papers \cite{LM} of Lada and Martin or \cite{LS} of Lada and Stasheff, as references for \Linfs and graded vector spaces in general.

\subsection{Basics on Graded Vector Spaces}
We begin by fixing some notations regarding graded vector spaces over a field $\mathbb{K}$ of characteristic zero, most of the times $\mathbb{K} = \mathbb{R}$. We will work with $\mathbb{Z}$-graded vector spaces and we will omit to mention the group, just saying that the objects are graded vector spaces or objects in {\sf grVec}.

The {\bf total vector space\index{Total vector space} associated} to a graded vector space $V$ will be denoted by $V^{\oplus} := \bigoplus_{i \in \mathbb{N}}{V_i}$. The opposite procedure of starting with the total vector space $U$ and splitting into degrees will be done via a {\bf grading}\index{Grading on a vector space}: a pair $(\varphi, V)$ where $V$ is a graded vector space and $\varphi \colon U \longrightarrow V^{\oplus}$ is an isomorphism. The elements $u \in \, U$ such that $\varphi(u) \in \, V_i$ are called {\bf homogeneous elements of degree} $\deg(u) := i$. We will use graded vector spaces with an additional grading, that is why it is important to fix this notation.

A sign rule\index{Sign rule} is a group homomorphism $|\cdot| \colon \mathbb{Z} \longrightarrow \mathbb{Z}_2$. Usually we will use the parity morphism. Many times we will work with expressions of the kind $(-1)^{|v|}$. We want to point out that in this particular case $(-1)^{|v|}=(-1)^{\deg{v}}$. A {\it sign convention} is a family of isomorphisms $V \otimes W \cong W \otimes V$  that takes care of the permutation of the elements with the given sign rule: $v \otimes w \cong (-1)^{|v||w|} w \otimes v$.

A graded vector space is said to be {\bf concentrated} in degrees $H \subset \mathbb{Z}$ if $V_i = \{0\}$ for all $i \notin \, H$. If $V$ is a graded vector space we define the {\bf truncation} of $V$ in degrees $H \subset \mathbb{Z}$ to be the graded vector space $\mathsf{tr}_{H}(V)$ which is given degree\=/wise by $\mathsf{tr}_{H}(V)_i = V_i$ if $i \in H$ and  $\mathsf{tr}_{H}(V)_i = 0$ otherwise. If $H = \{0, \ldots, i\}$ we write $\mathsf{tr}_i(V)$ instead of $\mathsf{tr}_{\{0, \ldots, i\}}(V)$. 

We will denote the reversing degree functor by a set-minus sign: $\smallsetminus \colon {\sf grVec} \longrightarrow {\sf grVec}$ defined on objects by $(\smallsetminus V)_i := V_{-i}$ for all graded vector spaces $V$ and all integers $i$. The choice of this minus sign and not the usual one is so that there is no ambiguity when talking about degree reversed morphisms $\smallsetminus f$.

Degree shifts will be denoted as follows: $V[i]$ is called $V$ shifted by $i$ and it is defined degree\=/wise by $V[i]_j := V_{i+j}$ for all $j \in \, \mathbb{Z}$. We can view the degree shift $[i]$ as a functor from graded vector spaces to itself. The degree shift functor {\bf is not compatible} with the tensor representation of $\mathcal{S}_n$. That means that the following diagram {\bf is not commutative}:

\begin{center}
\begin{tikzpicture}[description/.style={fill=white,inner sep=2pt}]\label{dec:tres}
\matrix (m) [matrix of math nodes, row sep=3em,
column sep=2.5em, text height=1.5ex, text depth=0.25ex]
{ (V)^{\otimes n} & (V[1]^{\otimes n})[-n] \\
  (V)^{\otimes n} & (V[1]^{\otimes n})[-n] \\};
\path[->,font=\scriptsize]
(m-1-1) edge node[auto] {$$} (m-1-2)
(m-2-1) edge node[auto] {$$} (m-2-2)
(m-1-1) edge node[auto] {$\sigma$} (m-2-1)
(m-1-2) edge node[auto] {$\sigma[-n]$} (m-2-2);
\end{tikzpicture}
\end{center}

If $v$ is an element of a $G$\=/graded vector space $V$, $v_{[g]}$ will denote the image of $v$ after $[g]$, $v_{[g]} \in \, V[g]$.

As we mentioned previously there exist an isomorphism between $\mathbb{R}[i] \otimes V$ and $V[i]$ for every $i \in \, \mathbb{Z}$ and $V$ a graded vector space. This conventions is called {\bf shifted to the left} as in \cite{FM}. 

One of the first consequences that the choice of a shifting convention brings is that $V[i][j]$ and $V[j][i]$ are non\=/trivially isomorphic:
$$V[i][j] \cong \mathbb{K}[j] \otimes V[i] \cong \mathbb{K}[j] \otimes \mathbb{K}[i] \otimes V \stackrel{\tau_{\mathbb{K}[i],\mathbb{K}[j]} \otimes id_V}{\cong} \mathbb{K}[i] \otimes \mathbb{K}[j] \otimes V \cong V[j][i].$$
That is $(v_{[i]})_{[j]} \stackrel{\cong}{\mapsto} (-1)^{i j} (v_{[j]})_{[i]}$ for every $v \in \, V$.

\subsubsection{Decalage}

It is convenient to give a name to an isomorphism that takes care of all the signs associated with a permutation of the factors:

\begin{df}[Decalage Isomorphism\index{Decalage}]
For every integer $n$ and any graded vector space $V$ we have the following isomorphism, called decalage:
\begin{eqnarray*}
	\dec \colon V[-1] \otimes \cdots \otimes V[-1] & \longrightarrow & (V \otimes \cdots \otimes V)[-1]\stackrel{n}{\cdots}[-1] \\
	{v_1}_{[-1]} \otimes \cdots \otimes {v_n}_{[-1]} & \mapsto & (-1)^{\sum_{i=1}^n {(n-i)|v_i|}} (v_1 \otimes \cdots \otimes v_n)_{[-1]\stackrel{n}{\cdots}[-1]}.
\end{eqnarray*}
\end{df}

Let us fix a graded vector space $V$. The free graded algebra generated by $V$ will be denoted by $T(V)_g := \bigoplus_{n=0}^{\infty}{(V^{\otimes n})_g}$ with the multiplication given by the tensor product. This is called the {\bf graded tensor algebra\index{Tensor algebra}} and it has two gradings, the original in $V$ and the polynomial one.

\begin{df}[Koszul sign\index{Koszul sign}]
Let $V$ be a graded vector space and let $\sigma \in \, \mathcal{S}_n$ be a permutation. Let $v = (v_1 \otimes \cdots \otimes v_n) \in \, V^{\otimes n}$. We define the Koszul sign of $\sigma$ and $v$ to be $\epsilon({\sigma}) = \epsilon(\sigma; v ) \in \, \{-1, 1\}$ where 
$$ \sigma \cdot (v_1 \otimes \cdots \otimes v_n) = \epsilon(\sigma) \, (v_{\sigma(1)} \otimes \cdots \otimes v_{\sigma(n)} ).$$
\end{df}

The Koszul sign depends on the sign rule $\tau$ and we should better write $\epsilon_{\tau}(\sigma)$. As we have fixed $\tau$ we write just $\epsilon(\sigma)$. If we take as the isomorphism $-\tau$ we get $\epsilon_{-\tau}(\sigma) = (-1)^{\sigma}\epsilon(\sigma)$ where $(-1)^{\sigma}$ denotes the usual sign of a permutation.

We can view the free graded commutative\index{Tensor algebra!Free graded commutative algebra} algebra as symmetrization $S(V)$ of the tensor algebra. The anti\=/symmetrization of the tensor algebra is denoted by $\wedge V$\index{Tensor algebra!Exterior algebra}. As in the non\=/graded world, those algebras are degree\=/wise isomorphic to the co\=/invariants of the canonical action of $\mathcal{S}_n$ on $V^{\otimes n}$ and of the action given by $(-1)^{\sigma} \sigma$ respectively. The product on $\wedge V$ is usually denoted by $\wedge$.

We will understand the decalage isomorphism in the case where all the $V_i$ are the same graded vector space as a natural transformation between two functors. First we have to define the target category of those two functors. We fix $n \geqslant 1$ an integer. The category of {\sf grRep}($\mathcal{S}_n$) is given by degree\=/wise linear representations of $\mathcal{S}_n$ as objects and $\mathcal{S}_n$\=/equivariant morphisms in {\sf grVec} as morphisms.

The following proposition gives an answer to the non\=/commutativity of \ref{dec:tres}. This point of view is original of C. Blohmann (private communication):

\begin{pp}\label{nati}\index{Decalage!natural isomorphism}
For every $n \geq 1$ integer, the decalage isomorphism is a natural isomorphism between the functors $F, G \colon {\sf grVec} \longrightarrow {\sf grRep}(\mathcal{S}_n)$ where for every graded vector space $V$ we have $F(V) := (\mathbb{K}[-1] \otimes V)^{\otimes n}$ with the $\mathcal{S}_n$ action induced by $\tau$ and $G(V) := (\mathbb{K}[-1])^{\otimes n} \otimes V^{\otimes n}$ with the induced  action by the co\=/product of the canonical actions of $\mathcal{S}_n$ on $(\mathbb{K}[-1])^{\otimes n}$ and $V^{\otimes n}$ respectively.
\end{pp}

The proof of this fact is equivalent to the usual statement that the decalage isomorphism makes the following diagram commute:

\begin{center}
\begin{tikzpicture}[description/.style={fill=white,inner sep=2pt}]
\matrix (m) [matrix of math nodes, row sep=3em,
column sep=2.5em, text height=1.5ex, text depth=0.25ex]
{ (V[-1])^{\otimes n} & (V^{\otimes n})[-1]\stackrel{n}{\cdots}[-1] \\
  (V[-1])^{\otimes n} & (V^{\otimes n})[-1]\stackrel{n}{\cdots}[-1] \\};
\path[->,font=\scriptsize]
(m-1-1) edge node[auto] {$\dec$} (m-1-2)
(m-1-1) edge node[auto] {$\sigma$} (m-2-1)
(m-2-1) edge node[auto] {$\dec$} (m-2-2)
(m-1-2) edge node[auto] {$(-1)^{\sigma}\sigma[-n]$} (m-2-2);
\end{tikzpicture}
\end{center}

If we quotient the left and the right side of the previous diagram by the co\=/invariants of the action in each of the cases we get the following result which is the relevant one when talking about Lie and \Linfs .

\begin{cl}
The decalage isomorphism induces an isomorphism also called decalage for every graded vector space $V$ and every $n\geq 1$ integer:
$$\dec \colon S^n(V[-1]) \stackrel{\sim}{\longrightarrow} (\wedge^{n} V)[-1]\stackrel{n}{\cdots}[-1].$$
\end{cl}

\subsubsection{Morphisms}

Morphisms of graded vector spaces can be enriched in graded vector spaces, this will be denoted by $\underline{\mathrm{Hom}}(V, W)$ for every pair of $G$\=/graded vector spaces $V$ and $W$:
\begin{eqnarray*}
	\underline{\mathrm{Hom}}(V, W)_j &\cong \prod_{i \in G}{\mathrm{Hom}_{\mathcal{V}ec}(V_{i-j}, W_i)}. 
\end{eqnarray*} 

The elements of $\underline{\mathrm{Hom}}(V, W)_i$ are called {\bf morphisms of degree $g$} from $V$ to $W$. Observe that morphisms of degree zero $\underline{\mathrm{Hom}}(V, W)_0 \cong \mathrm{Hom}(V, W)$ as vector spaces. 

The graded endomorphism space of a graded vector space $V$ is defined to be $\underline{\mathrm{End}}(V) := \underline{\mathrm{Hom}}(V,V)$. It has the structure of a graded associative algebra (non\=/commutative in general), where the product is given by composition. If $A$ is an associative algebra and $a \in A_{\deg{a}}$ we can construct $\underline{a} \in \underline{\mathrm{End}}(A)_{\deg{a}}$ given by left multiplication by $a$.

It is possible to define a graded commutator on $\underline{\mathrm{End}}(A)$ by setting  for every $f$ and $g$ in $\underline{\mathrm{End}}(A)$, $[f,g] := f \circ g -(-1)^{|f||g|}g \circ f$. This construction can be done in general for a graded associative algebra. With this notation, a {\bf graded derivation} of degree $i \in \mathbb{N}$ of $A$ is an element $D$ in $\underline{\mathrm{End}}(A)_i$ such that $[D,\underline{a}] = \underline{D(a)}$ for every $a \in \, A$.

Let $V$ and $W$ be two $G$\=/graded vector spaces. A morphism $f \in \, \underline{\mathrm{Hom}}(V^{\otimes n}, W)$ is symmetric if and only if it descends to a linear map $f \in \, \underline{\mathrm{Hom}}(S(V)^{n}, W)$ and if it is anti\=/symmetric whenever it descends to a corresponding $f \in \, \underline{\mathrm{Hom}}(\wedge^{n} V, W)$. Therefore, we denote then graded vector spaces of symmetric and anti\=/symmetric maps between $V^{\otimes n}$ and $W$ by $\underline{\mathrm{Hom}}(S^{n}(V), W)$ and $\underline{\mathrm{Hom}}(\wedge^{n} V, W)$, respectively.

Yet another point to take into account: for every $V$ and $W$ graded vector spaces and every pair of integers $(i,j)$: 
$$\underline{\mathrm{Hom}}(V, W)_{j} \cong \mathrm{Hom}(V, W[j]) \cong \underline{\mathrm{Hom}}(V[-i], W[-i])_j.$$

The first isomorphism is trivial, but not the second one which picks the Koszul sign of $\tau_{\mathbb{R}[i],\mathbb{R}[j]}$. Combining this fact and the pullback of the decalage isomorphism we get for every $n \geqslant 1$ and $i$ integers the following isomorphism also called decalage:
$$\underline{\mathrm{Hom}}(V^{\otimes n}, W)_{i} \stackrel{\scriptsize{(-1)^{ni}}}{\cong} \underline{\mathrm{Hom}}((V^{\otimes n})[-1]\stackrel{n}{\cdots}[-1], W[-n])_{i} \stackrel{\scriptsize{\dec^{*}}}{\cong} \underline{\mathrm{Hom}}((V[-1])^{\otimes n}, W[-n])_i .$$

\begin{pp}\label{stas}\index{Decalage!isomorphism of maps}
The pullback of the decalage isomorphism induces for every pair of graded vector spaces $V$ and $W$, every pair of integers $(i, j)$ and every integer $n \geq 1$ a natural isomorphism
$$\dec \colon \underline{\mathrm{Hom}}(\wedge^n V, W)_{i} \, \stackrel{\sim}{\longrightarrow} \, \underline{\mathrm{Hom}}(S^n(V[-1]), W[j])_{i-j-n}$$
$$\dec (f) ({v_1}_{[-1]} \otimes \cdots \otimes {v_n}_{[-1]}) := (-1)^{ni} (-1)^{\sum_{k=1}^n {(n-k)|v_k|}} f(v_1 \wedge \cdots \wedge v_n)_{[-j]}.$$
\end{pp}


\subsection{The Richardson\=/Nijenhuis bracket and graded operations}

The Richardson\=/Nijenhuis bracket is a bracket on automorphisms of the tensor algebra. Our approach is to define it making use of some notation coming from the theory of operads. This will simplify many of the computations in this document. 

Very often we have to work with a very specific kind of permutation: {\bf the unshuffles}. 

\begin{df}[Unshuffles\index{Unshuffle}]
Let $\sigma \in \mathcal{S}_n$ be a permutation and $\{1 \leqslant p_i \leqslant n \}_{i=1}^{q}$ be a collection of integers such that $\sum_{i=1}^{q} p_i = n$, we say that $\sigma$ is a $(p_1, \ldots, p_q)$\=/unshuffle if $\sigma(i) < \sigma (i+1)$ for all $i \notin \{ k_j \}_{j=1}^{q}$ where $k_j = \sum_{i=1}^{j} p_i$.
\end{df}

The set of $(p_1, \ldots, p_q)$\=/unshuffles is denoted by $\mathcal{S}{h}(p_1, \ldots, p_q)$. When $q=2$ the notation most commonly found on the literature is $\mathcal{S}{h}_{p_2}^{p_1}$.

\begin{rk}\label{MUS}
	As far as we will be working with unshuffles it is important to understand their behavior. Given $\sigma \in \mathcal{S}h_j^i$, $\tau_1 \in \mathcal{S}h_{i-k_1}^{k_1}$ and $\tau_2 \in \mathcal{S}h_{j-k_2}^{k_2}$, the composition $(\tau_1, \tau_2) \circ \sigma$ belongs to $\mathcal{S}h(k_1,i-k_1,k_2,j-k_2)$. The reason is simple: the composition of two increasing functions is increasing. Conversely, given $\alpha \in \mathcal{S}h(p_1, p_2, p_3, p_4)$ it is possible to find $\sigma \in \mathcal{S}h_{p_3 + p_4}^{p_1 + p_2}$, $\tau_1 \in \mathcal{S}h_{p_2}^{p_1}$ and $\tau_2 \in \mathcal{S}h_{p_4}^{p_3}$ such that $\alpha = (\tau_1, \tau_2) \circ \sigma$ just by taking the respectively ordered terms in the complementary sets $\{\alpha(1), \ldots, \alpha(p_1 + p_2)\}$ and $\{\alpha(p_1 + p_2), \ldots, \alpha(n)\}$.
\end{rk}

\begin{df}[$n$\=/ary operations\index{nary @ $n$\=/ary operation}]
Let $(\mathcal{C}, \otimes)$ be a monoidal category. Let $n$ and $m$ be two positive integers and $V_i$, $W_j$ objects in $\mathcal{C}$ for every $i \in \, \{1, \ldots, n\}$ and every $j \in \, \{1, \ldots, m\}$. A morphism $f \in \, \mathrm{Hom}_{\mathcal{C}}\left( V_1 \otimes \cdots \otimes V_n, W_1 \otimes \cdots \otimes W_m \right)$ is called an $n$\=/ary operation.
\end{df}

There are several ways of composing an $n$\=/ary and a $p$\=/ary operation. 

\begin{df}
Let $(\mathcal{C}, \otimes)$ be a monoidal category. Let $n, m$ and $p$ be positive integers, $1 \leqslant q \leqslant n$ and $i, j$ be positive integers such that $n = i-1 + q + j-1$. Let $V_a$, $V^{\prime}_{a^{\prime}}$, $V^{\prime \prime}_{a^{\prime \prime}}$, $W_b$ and $X_c$ be objects in $\mathcal{C}$ for $a \in \, \{1, \ldots, i-1\}$, $a^{\prime} \in \, \{1, \ldots, q\}$, $a^{\prime \prime} \in \, \{1, \ldots, j-1\}$, $b \in \, \{1, \ldots, m\}$ and $c \in \, \{1, \ldots, p\}$. Let $f$ be an $n$\=/ary operation and $g$ be a $p$\=/ary operation, where $g \colon X_1 \otimes \cdots \otimes X_p \longrightarrow V^{\prime}_1 \otimes \cdots \otimes V^{\prime}_{q}$
$$f \colon V_1 \otimes \cdots \otimes V_{i-1} \otimes V^{\prime}_1 \otimes \cdots \otimes V^{\prime}_{q} \otimes V^{\prime \prime}_1 \otimes \cdots \otimes V^{\prime \prime}_{j-1} \longrightarrow W_1 \otimes \cdots \otimes W_m$$
We define the morphism $f \circ_i g$ by
$$f \circ_i g \colon V_1 \otimes \cdots \otimes V_{i-1} \otimes X_1 \otimes \cdots \otimes X_p \otimes V^{\prime \prime}_1 \otimes \cdots \otimes V^{\prime \prime}_{j-1} \longrightarrow W_1 \otimes \cdots \otimes W_m$$
$$f \circ_i g \left( v_1 \otimes \cdots \otimes v_{i-1} \otimes x_1 \otimes \cdots \otimes x_p \otimes v^{\prime \prime}_1 \otimes \cdots \otimes v^{\prime \prime}_{j-1} \right) := \phantom{aaaaaaaaaiaaaaaa}$$
$$\phantom{aaaaaaaaaaaiaaaa} = f \left( v_1 \otimes \cdots \otimes v_{i-1} \otimes g \left( x_1 \otimes \cdots \otimes x_p \right) \otimes v^{\prime \prime}_1 \otimes \cdots \otimes v^{\prime \prime}_{j-1} \right)$$ for every homogeneous element $v_1 \otimes \cdots \otimes v_{i-1} \otimes x_1 \otimes \cdots \otimes x_p \otimes v^{\prime \prime}_1 \otimes \cdots \otimes v^{\prime \prime}_{j-1}$.
\end{df}

In particular when $i=1=j$ we have $f\circ_1 g = f \circ g$. When all the objects involved in the previous definition are the same, we get different compositions of the same two maps letting $i$ vary from $1$ to $n-q+1$, if we fix $r$ to be $r:= n+p-q$, $f \circ_{i} g \colon V^{\otimes r} \longrightarrow V^{\otimes m}$ given for every $v = v_1 \otimes \cdots \otimes v_{r} \in \, V^{\otimes r}$ 
$$f \circ_{i} g (v) := f(v_1 \otimes \cdots \otimes v_{i-1} \otimes \, g(v_{i} \otimes \cdots \otimes v_{i+p-1}) \, \otimes v_{i+p} \otimes \cdots \otimes v_{r}).$$

\begin{df}[Insertion\index{Insertion}]
Let $(\mathcal{C}, \otimes, \tau)$ be a symmetric monoidal category. We define $f \ins g$ to be the insertion operator of $f \colon V^{\otimes n} \longrightarrow V^{\otimes m}$ and $g \colon V^{\otimes p} \longrightarrow V^{\otimes q}$ with $q \leq n$ by setting for every $v \in \, V^{\otimes n+p-q}$
$$f \ins g (v) := \sum_{\sigma \in \, \mathcal{S}{h}_{p}^{n-q}} \left( f \circ_1 g \right) \, (\sigma \cdot v).$$
\end{df}

Observe that the insertion operator depends on the sign rule chosen, $\tau$. Recall that we have a preferred choice of an isomorphism $\tau$. We write, as we did for the Koszul sign, $f \ins g$ for the insertion operator associated to $\tau$ and we will call it {\bf symmetric insertion} of $g$ into $f$ or just {\bf insertion}. If we use $-\tau$ we will call $f \stackrel{\scriptscriptstyle -\tau}{\ins} g$ the {\bf anti\=/symmetric insertion} of $g$ into $f$.

The following is a well-known result:

\begin{pp}\index{Insertion!symmetric}
Let $(\mathcal{C}, \otimes, \tau)$ be a symmetric monoidal category. Given any two morphisms $f \colon V^{\otimes n} \longrightarrow V$ and $g \colon V^{\otimes p} \longrightarrow V$. Then
\begin{enumerate}
	\item If $f$ and $g$ are both symmetric then so it is $f \ins g$.
	\item If $f$ and $g$ are both anti\=/symmetric then so it is $f \stackrel{\scriptscriptstyle -\tau}{\ins} g$.
\end{enumerate}
\end{pp}

Note that the insertion operation is not associative in general,
$$(f \ins g) \ins h \ncong f \ins (g \ins h).$$ 

For every pair of graded vector spaces $V$ and $W$, we can consider maps in $\underline{\mathrm{Hom}}(V,W)$. We can define the insertion operator also on morphisms of degree $n$ not only for $n = 0$.

\begin{df}[Richardson\=/Nijenhuis bracket\index{Richardson\=/Nijenhuis bracket}]
Let $V$ be a graded vector space. Let also $m \leqslant n$ and $m \leqslant p$ be three positive integers. Then for every $f \in \underline{\mathrm{Hom}}(V^{\otimes n}, V^{\otimes m})$ and every $g \in \underline{\mathrm{Hom}}(V^{\otimes p}, V^{\otimes m})$ we define the Richardson\=/Nijenhuis bracket $[f, g]_{\mathrm{RN}} \in \, \underline{\mathrm{Hom}}(V^{\otimes n+p-m}, V^{\otimes m})_{|f|+|g|}$ to be:
$$[f,g]_{\mathrm{RN}} := f \ins g -(-1)^{|f||g|} g \ins f.$$
\end{df}

As before, the Richardson\=/Nijenhuis bracket depends on the isomorphism $\tau$. We always assume that we are working with our preferred $\tau$ and if we want to use $(-\tau)$ we will denote the bracket $[-,-]_{\mathrm{RN}}^{-\tau}$ and call it the {\bf anti\=/symmetric Richardson\=/Nijenhuis bracket}.

When $m=p=n=1$ the Richardson\=/Nijenhuis bracket is defined on $\underline{\mathrm{End}}(V)$ and it is the usual {\bf commutator}. $(\underline{\mathrm{End}}(V), [-,-]_{\mathrm{RN}})$ is a graded Lie algebra.

With this notation, we are ready to state one of the main tools in proving the equivalence of the various definitions of \Linfs .

\begin{pp}\label{insins}
Let $V$ be a graded vector space, $i$ and $j$ two positive integers. Let also $f \in \underline{\mathrm{Hom}}(V^{\otimes j}, V)$ and $g \in \underline{\mathrm{Hom}}(V^{\otimes i}, V)$. Let $n+1=i+j$, then
$$ (-1)^{|f|(i-1)} \dec \left( f \stackrel{\scriptscriptstyle -\tau}{\ins} g \right) =  \dec(f) \ins \dec(g).$$
\end{pp}

We include the proof to get the reader comfortable with the use of the newly introduced notation.

\dem
Let $w = {v_1}_{[-1]} \otimes \cdots \otimes {v_n}_{[-1]} \in \, (V[-1])^{\otimes n}$
\begin{eqnarray*}
&& (\dec f \circ_{1} \dec g) (w) = (-1)^{|f|j +|g|i} (-1)^{\sum_{k=1}^{i}(i-k)|v_k|} (-1)^{\sum_{k=i+1}^{n} (j-k)|v_k|} \cdot \\
&\cdot& (-1)^{(j-1)\left(\sum_{k=1}^{i} (i-k)|v_k| +|g| \right)} f \left(g({v_1} \wedge \cdots \wedge {v_i}) \wedge {v_{i+1}} \wedge \cdots \wedge v_{n} \right)_{[-1]} \\
&=& (-1)^{|f|j +|g|n} (-1)^{\sum_{k=1}^{n}(n-k)|v_k|} f \left(g({v_1} \wedge \cdots \wedge {v_i}) \wedge {v_{i+1}} \wedge \cdots \wedge v_{n} \right)_{[-1]} \\
&=& (-1)^{|f|j + |g|n +(|f| + |g|)n} \dec(f \circ_{1} g) ({v_1}_{[-1]} \otimes \stackrel{n}{\cdots} \otimes {v_n}_{[-1]}).
\end{eqnarray*}
Since $|f|j + |g|n +(|f| + |g|)n = |f|(i-1)$ we have that 
$$(\dec f) \circ_{1} (\dec g) = (-1)^{|f|(i-1)} \dec(f \circ_1 g). \textrm{ Now}$$
\begin{eqnarray*}
\dec(f \stackrel{\scriptscriptstyle{-\tau}}{\ins} g) w &=& (-1)^{n(n-(|f|+|g|))} (f \stackrel{\scriptscriptstyle{-\tau}}{\ins} g) (\dec(w)_{[n]})_{[-1]} \\
&=& (-1)^{n(n-(|f|+|g|))} \left( \sum_{\mathcal{S}h_{j-1}^{i}} (-1)^{\sigma} (f \circ_1 g) ((\sigma \cdot (\dec(w)))_{[n]} \right)_{[-1]} \\
&=& (-1)^{n(n-(|f|+|g|))} \left( \sum_{\mathcal{S}h_{j-1}^{i}} (f \circ_1 g) (\dec(\sigma \cdot w)_{[n]}) \right)_{[-1]}\\
&=& \sum_{\mathcal{S}h_{j-1}^{i}} \dec(f \circ_1 g) (\sigma \cdot w) = \sum_{\mathcal{S}h_{j-1}^{i}} (-1)^{|f|(i-1)} (\dec f \circ_1 \dec g) (\sigma \cdot w) \\
&=& (-1)^{|f|(i-1)} (\dec f \ins \dec g) w.
\end{eqnarray*}
\qed\\

In particular, the last Proposition tells us that $f \ins g = 0$ if and only if the corresponding insertion $\dec(f) \stackrel{\scriptscriptstyle -\tau}{\ins} dec(g) = 0$. We remark that if $|f|$ is odd $[f,f]_{\mathrm{RN}} = 2 f \ins f$ and if it is even $[f,f]_{\mathrm{RN}} = 0$.

\begin{rk}\label{LA}
Graded Lie algebras (Lie algebras inner to graded vector spaces) can be rephrased in terms of the anti\=/symmetric insertion operator. A graded Lie algebra is a pair $(A, \lambda)$ where $A$ is a graded vector space and $\lambda \in \mathrm{Hom}(A\wedge A,A)$ such that $\lambda \stackrel{\scriptscriptstyle -\tau}{\ins} \lambda = 0$.
\end{rk}


\subsection{Gerstenhaber Algebras}

The definition of a Gerstenhaber algebra is the next and last step on this first chapter. The usual definition of a Gerstenhaber algebra involves shifting by one and minus one the underlying graded vector space often, we will use the arrow notation of Stasheff in \cite{LS}.

\begin{df}[Gerstenhaber Algebra\index{Gerstenhaber!algebra}]\label{GAD}
A triple $(A, [-,-], \cdot)$ is called a Gerstenhaber algebra if $(A[1], [-,-])$ is a graded Lie algebra and $(A, \cdot)$ is an associative symmetric graded algebra satisfying the Leibniz rule, that is, for every $a, b \textrm{ and } c$ in $A$:
$$\uparrow [\downarrow a , \downarrow (bc)] = (\uparrow [\downarrow a , \downarrow b])c + (-1)^{|\downarrow a||b|}b(\uparrow [\downarrow a , \downarrow c]).$$
\end{df}

The Leibniz rule can be stated by saying that $\uparrow [\downarrow a, \downarrow-]$ is a derivation of degree $|\downarrow a|$.

\begin{cl}\label{GA2}
A triple $(A, \lambda, \mu)$ is a Gerstenhaber algebra if and only if $(A, \nu, \mu)$ where $\nu := \dec(\lambda)$ satisfy the following conditions:
$(A, \mu)$ is an associative symmetric graded algebra, $\nu$ is a symmetric morphism of degree $-1$ such that:
\begin{enumerate}
	\item $[\nu, \nu]_{\mathrm{RN}} = 0$.
	\item For every $a \in \, A$, $\nu(a,-) \in \, \underline{\mathrm{Der}}(A)_{\deg(a)-1}$.
\end{enumerate}
\end{cl}

The first condition is called {\it the $3$\=/Jacobi identity} and the second one is the {\it Leibniz rule}.
 
\dem 
By Remark \ref{LA} $\lambda \in \, \mathrm{Hom}(A[1] \wedge A[1], A[1])$ is such that $\lambda \stackrel{\scriptscriptstyle -\tau}{\ins} \lambda = 0$. Using the isomorphism \ref{nati} we get $\nu \in \, \underline{\mathrm{Hom}}(S_2(A),A)_{-1}$ and by proposition \ref{insins} $\nu \ins \nu = \lambda \stackrel{\scriptscriptstyle -\tau}{\ins} \lambda = 0$. Even more, since $\nu$ is of odd degree $[\nu,\nu]_{\mathrm{RN}} = 2 \nu \ins \nu = 0$. The fact that we can translate the Leibniz rule in terms of graded derivation was already mentioned after the definition of a Gerstenhaber algebra (\ref{GAD}). 
\qed\\

A justification for the term $3$\=/Jacobi will be given in the following section. From now on every time we work with a Gerstenhaber algebra it would be presented by $(A, \nu, \mu)$ in the hypothesis of Corollary \ref{GA2} instead of $(A, \lambda, \mu)$ as a consequence of the Corollary.

We can set for every $n \in \, \mathbb{Z}$, $n \geq 2$ the morphism $(\cdot)^n \colon A^{\otimes n} \longrightarrow A$ defined recursively by $(\cdot)^2 = (\cdot)$ and $(\cdot)^{n+1} = (\cdot) \circ_1 (\cdot)^n$ when $n \geq 2$. Since $(\cdot)$ is associative it does not matter the order in which the compositions are made. Thus we can write $(\cdot)^{n}(a_1 \otimes \cdots \otimes a_n) = a_1 \cdot \ldots \cdot a_n = a_1 \cdots a_n$ without any need of specifying parenthesis. The same goes when we denote it by $\mu$.

We want to be able to calculate composition of the different maps that define a Gerstenhaber algebra structure. We have the following lemma.

\begin{lm}[Gerstenhaber Calculus\index{Gerstenhaber!calculus}]\label{GHC}
Let $(A, \nu, \mu)$ be a Gerstenhaber algebra, $n \geqslant 2$, for every $a= a_1 \otimes \cdots \otimes a_n \in \, A^{\otimes n}$ we denote by $a_{i,j}:=a_i \otimes \cdots \otimes a_j$ where $1 \leqslant i \leqslant j \leqslant n$. Then:
\begin{enumerate}
	\item The following generalized Leibniz rule holds
	$$(\nu \circ_2 \mu^{n-1}) \, (a) = (\mu^{n-1} \ins \nu(a_1, -)) a_{2,n}.$$
	\item Let $f_a \colon A \longrightarrow A$ be the linear map defined by $\f_a(c) := \nu(a_1 \otimes a_2) \cdot \nu(a_n \otimes c)$. Then for every $n \geq 4$ the following equality holds:
	\begin{eqnarray*}
	(-1)^{|a_n||a_{3,n-1}|} (\nu \circ_1 (\mu^{n-2} \circ_1 \nu)) \, (a) &=& (\mu^{n-2} \circ_1 (\nu \circ_1 \nu)) \, (a_{1,2} \otimes a_n \otimes a_{3,n-1}) \\
	&-& (-1)^{|a_{1,2}|} (\mu^{n-2} \ins f_a) \, (a_{3,n-1}).
	\end{eqnarray*}
\end{enumerate}
\end{lm}

\dem
\begin{enumerate}
	\item The proof will be given by induction. The induction hypothesis is given by the {\it Leibniz rule}
		\begin{eqnarray*}
			\nu(a_1 \otimes (a_2 \cdot a_3)) &=& \nu(a_1 \otimes a_2) \cdot a_3 +(-1)^{(|a_1 -1|)|a_2|} a_2 \cdot \nu(a_1 \otimes a_3) \\
			&=& \nu(a_1 \otimes a_2) \cdot a_3 +(-1)^{|a_2||a_3|} \nu(a_1 \otimes a_3) \cdot a_2.
		\end{eqnarray*}
		For the induction step we take $a= a_1 \otimes \cdots \otimes a_{n+1} \in \, A^{\otimes (n+1)}$ and $b_1 = a_1$, $b_2 = a_2 \cdot a_3$ and $b_i = a_{i+1}$ for every $i \in \, \{3, \ldots, n\}$.
		\end{enumerate}
		\begin{eqnarray*}
			&& \nu(a_1 \otimes (a_2 \cdots a_{n+1})) = \nu(b_1 \otimes (b_2 \cdots b_{n})) = \\
			&=& \sum_{\sigma \in \mathcal{S}{h}_{n-3}^1} \epsilon(\sigma) (-1)^{|b_{\sigma(1)+2}||b_2|} \nu(b_1 \otimes b_{\sigma(1)+2}) \cdot b_2 \cdot b_{\sigma(2)+2} \cdots b_{\sigma(n-2)+2} \\
			&+& \nu(a_1 \, \otimes \, (a_2 \, \cdot \, a_3)) \, \cdot \, a_4 \, \cdot \, \ldots \, \cdot \, a_{n+1} \\
			&=& \sum_{\sigma \in \mathcal{S}{h}_{n-3}^1} \epsilon(\sigma) \nu(b_1 \otimes b_{\sigma(1)+2}) \cdot b_2 \cdot b_{\sigma(2)+2} \cdots b_{\sigma(n-2)+2} \\
			&=& \sum_{\sigma \in \mathcal{S}{h}_1^1} \epsilon(\sigma) \nu (a_1 \otimes a_{\sigma(1) +1}) \cdot a_{\sigma(2) + 1} \cdot a_4 \cdots a_{n+1} \\
			&+& \sum_{\sigma \in \, \mathcal{S}{h}_{n-1}^1} \epsilon(\sigma) \nu(a_1 \otimes a_{\sigma(1)+1}) \cdot a_{\sigma(2)+1} \cdots a_{\sigma(n)+1}.
		\end{eqnarray*}
		\begin{enumerate}
\item[2.] Let $g_a \colon A \longrightarrow A$ be $g_a (c) := (-1)^{|c|(|a_{1,2}| -1)} \nu(a_n \otimes c) \cdot \nu(a_1 \otimes a_2)$ for every $c \in \, A$. Then:
\end{enumerate}
\begin{eqnarray*}
&& (-1)^{|a_n||a_{3,n-1}|} (\nu \circ_1 (\mu^{n-2} \circ_1 \nu)) \, (a) = \\ 
&=& (-1)^{|a_n|(|a_{1,2}| -1)} (\nu \circ_2 (\mu^{n-2} \circ_1 \nu)) \, (a_n \otimes a_{1,2} \otimes a_{3,n-1}) \\
&=& (-1)^{|a_n|(|a_{1,2}| -1)} \left( (\mu^{n-2} \circ_1 (\nu \circ_2 \nu)) \, (a_n \otimes a_{1,2} \otimes a_{3,n-1}) + (\mu^{n-2} \ins g_a) \, (a_{3,n-1}) \right) \\
&=& (\mu^{n-2} \circ_1 (\nu \circ_1 \nu)) \, (a_{1,2} \otimes a_n \otimes a_{3,n-1}) - (-1)^{|a_{1,2}|} (\mu^{n-2} \ins f_a) \, (a_{3,n-1}).
\end{eqnarray*}

	Since 
	\begin{eqnarray*}
	(-1)^{|a_{1,2}| +1} f_a (c) &=& (-1)^{|a_{1,2}| +1} \nu(a_1 \otimes a_2) \cdot \nu (a_n \otimes c) \\ 
	&=&(-1)^{(|a_{1,2}| +1)(|a_n| + |c|)} \nu(a_n \otimes c) \cdot \nu (a_1 \otimes a_2) \\
	&=& (-1)^{(|a^{1,2}| +1)|a_n|} g_a (c).
	\end{eqnarray*}
\qed


\subsection{\mbox{\boldmath ${L_{\infty}}$}\=/algebras}

We adopt the notion of an \mbox{$L_{\infty}$}\=/algebra on its graded\=/symmetric version. For the relations between the different definitions of \mbox{$L_{\infty}$}\=/algebras we refer to \cite[Appendix A]{FRZ} for instance. The relation between the symmetric and anti-symmetric definitions are due to the decalege isomorphism mentioned earlier.

\begin{df}[\mbox{\boldmath ${L_{\infty}}$}\=/algebra\index{Linfinity algebra @ $L_{\infty}$\=/algebra}]
An \mbox{$L_{\infty}$}\=/algebra is a graded vector space $L$ concentrated in positive degrees together with a collection of symmetric morphisms of degree $-1$ $\{ l_n \colon L^{\otimes n} \longrightarrow L \}_{n\in\mathbb{N}}$ satisfying the generalized Jacobi identity, that is, for every $n \in \, \NA$, $n \geqslant 1$:
$$ \sum_{i+j=n+1} l_j \ins l_i = 0.$$
\end{df}

We will denote for every $n \geqslant 1$, $J(n) := \sum_{i+j=n+1} l_j \ins l_i$. The relation $J(n) = 0$ will be called $n$\=/Jacobi. Often the maps $l_i$ are called brackets.

It is possible to define a co\=/differential on the co\=/free co\=/commutative connected co\=/algebra generated by $L$ out of $l$ (we will mention this in Proposition \ref{cochain}). The Jacobi identity reads in that case as $[l,l]_{\mathrm{RN}}=0$.

\mbox{$L_{\infty}$}\=/algebras with $l_i = 0$ unless $i = 2$ are called Lie$[1]$\=/algebras\index{Lie algebra!Lie$[1]$\=/algebra} as in \cite{MZ}. Ordinary Lie algebras are viewed as Lie$[1]$\=/algebras concentrated in degree $\{1\}$.
 
Hence \mbox{$L_{\infty}$}\=/algebras generalize graded Lie algebras and thus also conventional Lie algebras. Again we find the usual Lie algebra if $L_i = 0$ unless $i=1$ and where $l_i = 0$ unless $i = 2$ by applying Corollary \ref{GA2}.

\mbox{$L_{\infty}$}\=/morphism are morphisms of the corresponding coalgebras. Splitting into a family of maps we have the following compatibility conditions:

\begin{df}[\mbox{\boldmath ${L_{\infty}}$}-algebra morphism\index{Linfinity morphism @ $L_{\infty}$\=/morphism}]
An \mbox{$L_{\infty}$}\=/morphism between two \mbox{$L_{\infty}$}\=/algebras $(L, \{ l_n\})$ and $(V, \{ v_n\})$, $f \colon (L, \{ l_n\}) \longrightarrow (V, \{ v_n\})$ is a collection of symmetric morphisms (of degree $0$) $f:=\{f_k \colon L^{\otimes k} \longrightarrow V\}$ such that for every $n \in \NA$, $n \geqslant 1$ and every $x = x_1 \otimes \cdots \otimes x_n$ in $L^{\otimes n}$
$$\sum_{i+j =n +1} (f_j \ins l_i) (x) = \sum_{p_1 + \ldots p_q = n}^{p_i \geqslant 1} \sum_{\sigma \in \mathcal{S}{h}(p_1,\ldots, p_q)} \frac{1}{q!} (v_q \circ (f_{p_1} \otimes \cdots \otimes f_{p_q})) \, (\sigma \cdot x).$$
\end{df}

We recover the usual definition of morphism of Lie algebras if both $(L, \{ l_n\})$ and $(V, \{ v_n\})$ are Lie algebras after applying Corollary \ref{GA2}.

\Linfs and \mbox{$L_{\infty}$}\=/morphisms form a category.

\begin{df}[Strict \mbox{\boldmath ${L_{\infty}}$}\=/morphism\index{Linfinity morphism @ $L_{\infty}$\=/morphism!strict}]
An \mbox{$L_{\infty}$}\=/morphism between $(L, \{ l_n\})$ and $(V, \{ v_n\})$ is said to be strict if the only non\=/zero map is $f_1$. In that case a morphism $f_1 \in \textrm{Hom}(L, V)$ is a strict \Linfm when for every positive $n \in \, \mathbb{Z}$
$$ f_1 \circ l_n = v_n \circ (f_1 \otimes \stackrel{\scriptscriptstyle{n}}{\cdots} \otimes f_1).$$
\end{df}

Observe that given $(L, \{l_n\})$ an \Linf, $(\smallsetminus(L), \smallsetminus(l_1))$ is a cochain complex.

\begin{df}[\mbox{\boldmath ${L_{\infty}}$}\=/quasi\=/isomorphism\index{Linfinity morphism @ $L_{\infty}$\=/morphism!quasi\=/isomorphism}]
	An \Linfm between $(L, \{l_n\})$ and $(V, \{v_n\})$ is said to be a quasi\=/isomorphism if it induces an isomorphism at the level of homology ($(L, l_1)$ to $(V, v_1)$).
\end{df}

The other cochain complex is the one in \cite[Theorem 2.3]{LM}. The translation of that result to our conventions and notations is given in the following Proposition.

\begin{pp}\label{cochain}\index{Chevalley\=/Eilenberg!generalized co\=/differential}
	Let $(L, \{l_n\})$ be an \Linf concentrated in a finite number of degrees. Each $l_n$ can be extended on $\smallsetminus (S L)$ by 
	$$ \hat{l}_i := \otimes^{n-i+1} \ins l_i \colon (S^{n} L)_k \longrightarrow (S^{n-i+1} L)_{k-1}.$$
	They define linear maps $\hat{l}_i \colon \smallsetminus (S L)_k \longrightarrow \smallsetminus (S L)_{k+1}$. This construction gives cochain complex structure on  $(\smallsetminus (S L), \hat{l}:=\sum_{i=1}^{\infty}\hat{l}_i)$.
\end{pp}

That construction is a generalization of the Chevalley\=/Eilenberg cochain complex for a Lie algebra. Recall that given a (non\=/graded) Lie algebra $(\mathfrak{g}, [-,-])$, for every $x = x_1 \wedge \cdots \wedge x_n \in \, \wedge^{n} \mathfrak{g}$ there is a co\=/differential, called Chevalley\=/Eilenberg co\=/differential\index{Chevalley\=/Eilenberg!co\=/differential} which is given by:
\begin{eqnarray}
d_{\mathrm{CE}}^{*}(x) &:=& \sum_{\sigma \in \, \mathcal{S}h_{n-2}^2} (-1)^{\sigma(1)+\sigma(2)+1} [x_{\sigma(1)}, x_{\sigma(2)}] \wedge x_{\sigma(3)} \wedge \cdots \wedge x_{\sigma(n)} \\
&=& \sum_{\sigma \in \, \mathcal{S}h_{n-2}^2} \epsilon(\sigma) l_2(x_{\sigma(1)} \otimes x_{\sigma(2)}) \otimes x_{\sigma(3)} \otimes \cdots \otimes x_{\sigma(n)} \\
&=& (\otimes^{n-1} \ins l_2)(x) = \hat{l} (x).
\end{eqnarray}

For this reason, we will denote $d_{\textrm{CE}} := \widehat{l}$ when there is a single \Linf structure on a graded vector space $L$.


\section{Higher \mbox{\boldmath ${L_{\infty}}$}\=/algebra induced by a Gerstenhaber algebra}

Any Gerstenhaber algebra is a Lie algebra is we forget about the associative operation. But it is possible to use that extra piece of information to create another \Linf structure on the underlying graded vector space with non-trivial higher brackets.

Later in this section we look at the example of multi-vector fields and review the Cartan calculus on multi-vector fields using a different convention for the insertion of a multi-vector field on a differential form than the usual. This is done for a purpose that will be clear in the next section.

\subsection{Theorem about Gerstenhaber algebras}

If $(A, \nu, \mu)$ is a Gerstenhaber algebra, $(A, \nu)$ restricted to positive degrees is an \mbox{$L_{\infty}$}\=/algebra with vanishing higher brackets. But this structure is only obtained from $\nu$ forgetting the extra structure given by $\mu$. The main result of this section is the following proposition that asserts that it is possible to associate to a Gerstenhaber algebra an \mbox{$L_{\infty}$}\=/algebra with non\=/vanishing higher brackets. 

\begin{tm}\label{main}\index{Gerstenhaber!$L_{\infty}$\=/algebra}
Let $(A, \nu, \mu)$ be a Gerstenhaber algebra, then $(A, \{\nu_n\})$ restricted to positive degrees is an \mbox{$L_{\infty}$}\=/algebra where $\nu_1 := 0$, $\nu_2 := \nu$ and for every $n \geqslant 3$ $\nu_n := \mu^{n-1} \ins \nu$.
\end{tm}

\dem

\noindent{\it Proof of the degree, linearity and symmetry.}

Linearity is satisfied since both $\nu$ and $\mu$ are linear. The condition about the degree of the maps $\{\nu_n\}$ follows from the fact that $\nu$ is of degree $-1$ and $\mu$ is a morphism.

In order to see that each $\nu_n$ is symmetric we fix $n$ a positive integer and $\tau \in \mathcal{S}_n$. Since $\tau$ decomposes in products of permutations of the type $(i,i+1)$ and since the Koszul signs behaves as follows: $\epsilon(\sigma_1 \circ \sigma_2; x) = \epsilon(\sigma_1; \sigma_2 \cdot x) \epsilon(\sigma_2; x)$, it is enough to show that $\nu_n (i, i+1) \cdot x = \epsilon((i,i+1); x) \nu_n x$ for every $i \in \{1, \ldots, n-1 \}$. Given $\sigma  \in \mathcal{S}h_{n-2}^2$  let $\nu_n^{\sigma} x := \epsilon(\sigma) (\nu, \textrm{ id}^{\otimes (n-2)}) \sigma \cdot x$ we want to find $\sigma^{\prime} \in \mathcal{S}h_{n-2}^2$ such that $\nu_n^{\sigma} (i,i+1) \cdot x = \epsilon((i,i+1); x) \nu_n^{\sigma^{\prime}} x$. We have different options. If $\{ \sigma(i), \sigma(i+1) \}$ is fully in $\{ 1, 2 \}$ or in $\{ 3, \ldots , n \}$ we take $\sigma^{\prime} = \sigma$ since the fact that both $\nu$ and $\mu$ are symmetric gives the result that $\nu_n^{\sigma} (i,i+1) x = \epsilon((i,i+1); x)\nu_n^{\sigma} x$. If not we take $\sigma^{\prime} = \sigma \circ (i,i+1)$ which is a $(2,n-2)$\=/unshuffle.

Now it is clear that $\nu_n^{\sigma} (i,i+1) \cdot x = \epsilon((i,i+1); x) \nu_n^{\sigma^{\prime}} x$. In this way we have shown that $\nu_n (i,i+1) \cdot x = \epsilon(i,i+1; x) \nu_n x$ and since composition of permutation gives multiplication of Koszul signs we have that $\nu_n$ is symmetric.

\noindent{\it Proof of the Generalized Jacobi identity:}

During this proof, whenever $l = l_1 \otimes \cdots \otimes l_n \in L^{\otimes n}$ and $\sigma \in \, \mathcal{S}_n$ we will denote $l_{\sigma} := l_{\sigma(1)} \otimes \cdots \otimes l_{\sigma(n)}$ so that $\sigma \cdot l = \epsilon(\sigma; l) l_{\sigma}$. Given $a$ and $b$ two integers such that $1 \leqslant a \leqslant b \leqslant n$, we call $l_{\sigma(a,b)} := l_{\sigma(a)} \wedge \cdots \wedge l_{\sigma(b)}$.
In order to prove the generalized Jacobi identity we need to understand what $\nu_j \ins \nu_i$ means. Let us fix some $i, j \geq 1$ and let $n := i + j -1$. We further fix $x = x_1 \otimes \cdots \otimes x_n \in A^{\otimes n}$. 
$$(\nu_j \ins \nu_i)(x) = \sum_{\sigma_1 \in \, \mathcal{S}{h}_{j-1}^{i}} \epsilon(\sigma_1; x) (\nu_j \circ_1 \nu_i) \, (x_{\sigma_1}).$$
For every $\sigma_1 \in \, \mathcal{S}{h}_{j-1}^{i}$ and according to the previous equation we have to calculate 
$$\nu_i \, (x_{\sigma_1(1,i)}) = (\mu^{i-1} \ins \nu) \, (x_{\sigma_1(1,i)}) = \sum_{\sigma_2 \in \, \mathcal{S}{h}_{i-2}^{2}} \epsilon(\sigma_2; x_{\sigma_1(1,i)}) (\mu^{i-1} \circ_1 \nu) \, ((x_{\sigma_1(1,i)})_{\sigma(2)}).$$
For every $\sigma_1 \in \, \mathcal{S}{h}_{j-1}^{i}$ and every $\sigma_2 \in \, \mathcal{S}{h}_{i-2}^{2}$ we fix $$y^{\sigma_1, \sigma_2} := ((\mu^{i-1} \circ_1 \nu)((x_{\sigma_1(1,i)})_{\sigma(2)})) \wedge x_{\sigma_1(i+1,n)}.$$ We have to calculate
$$\nu_j \, (y^{\sigma_1, \sigma_2}) = (\mu^{j-1} \ins \nu) \, (y^{\sigma_1, \sigma_2}) = \sum_{\sigma_3 \in \, \mathcal{S}{h}_{j-2}^2} \epsilon(\sigma_3; y^{\sigma_1, \sigma_2}) (\mu^{j-1} \circ_1 \nu) \, (y_{\sigma_3}^{\sigma_1, \sigma_2}).$$
To sum up, we have that

$$(\nu_j \ins \nu_i) \, (x) := \sum_{\begin{smallmatrix} 
\sigma_1 \in \, \mathcal{S}{h}_{j-1}^{i} \\ 
\sigma_2 \in \, \mathcal{S}{h}_{i-2}^{2} \\ 
\sigma_3 \in \, \mathcal{S}{h}_{j-2}^2 
\end{smallmatrix}} \left( ((\mu^{j-1} \circ_1 \nu) \circ \sigma_3) \circ_1 ((\mu^{i-1} \circ_1 \nu) \circ \sigma_2) \right) \, (\sigma_1 \cdot x) .$$

We want to split the last sum in the cases $1=\sigma_3(1)$ and $1 \neq \sigma_3(1)$ (in that case $1 = \sigma_3(3)$). We denote the corresponding sums by $(\nu_j \ins \nu_i)_1 (x)$ and $(\nu_j \ins \nu_i)_3 (x)$ respectively.

{\bf Case $1 = \sigma_3(3)$.} 

Letting $\sigma_1$, $\sigma_2$  and $\sigma_3$ such that $1 = \sigma_3(3)$ varying all over the possibilities in the corresponding unshuffles gives the same that taking $\sigma \in \, \mathcal{S}{h}(2,i-2,2,j-3)$ by Remark \ref{MUS} by letting $\sigma := (\sigma_2, \sigma_3^{\prime}) \circ \sigma_1$ where $\sigma_3^{\prime}(k) := \sigma_3(k)$ for $k \in \, \{1, 2\}$ and $\sigma_3^{\prime}(k) := \sigma_3(k+1)$ whenever $k \in \, \{3, \ldots, j-1\}$. The resulting term on the sum without sign is
$$\nu(x_{\sigma(i+1)} \otimes x_{\sigma(i+2)}) \cdot \nu(x_{\sigma(1)} \otimes x_{\sigma(2)}) \cdot x_{\sigma(3)} \cdots x_{\sigma(i)} \cdot x_{\sigma(i+3)} \cdots x_{\sigma(n)}.$$
If we permute $\nu(x_{\sigma(i+1)} \otimes x_{\sigma(i+2)})$ with $\nu(x_{\sigma(1)} \otimes x_{\sigma(2)}) \cdot x_{\sigma(3)} \cdots x_{\sigma(i)}$ we get a total sign:
$$(-1)^{(\sum_{k=1}^{i} |x_{\sigma(k)}| -1)(|x_{\sigma(i+1)}|+|x_{\sigma(i+2)}|+1+|x_{\sigma(i+1)}|+|x_{\sigma(i+2)}|)}.$$
So we have that
$$(\nu_j \ins \nu_i)_3 (x) = \sum_{\begin{smallmatrix} 
\sigma_1 \in \, \mathcal{S}{h}_{j-1}^{i} \\ 
\sigma_2 \in \, \mathcal{S}{h}_{i-2}^{2} \\ 
\sigma_3 \in \, \mathcal{S}{h}_{j-2}^2 \\
\sigma_3(1) = 3
\end{smallmatrix}} \left( ((\mu^{j-1} \circ_1 \nu) \circ \sigma_3) \circ_1 ((\mu^{i-1} \circ_1 \nu) \circ \sigma_2) \right) \, (\sigma_1 \cdot x) =$$
$$\sum_{\sigma \in \mathcal{S}{h}(2,i-2,2,j-3)} (-1)^{\sum_{k=1}^{i} |x_{\sigma(k)}| -1} \left( (\mu^{i-1} \circ_1 \nu) \cdot (\mu^{j-2} \circ_1 \nu) \right) \, (\sigma \cdot x).$$


We are going to show that the sum over $i+j=n+1$ of $(\nu_j \ins \nu_i)_3 (x)$ is zero. Given any $\sigma \in \mathcal{S}{h}(2,i-2,2,j-3)$, let $i^{\prime} := j -1$ and $j^{\prime} := i+1$. There exists a unique $\tau \in \mathcal{S}{h}(2,i^{\prime}-2,2,j^{\prime}-3)$ such that $\sigma(k) = \tau(k+j-1)$ for every $k \in \, \{1, \ldots, i\}$ and $\sigma(k+i)=\tau(k)$ for every $k \in \, \{1, \ldots, j-1\}$. If $j \neq 1$ then $(i^{\prime}, j^{\prime})$ define a member on the sum $J(n) = \sum_{i+j=n+1} \nu_j \ins \nu_i$; if $j = 1$ $(\nu_1 \ins \nu_n)_3 = 0$ since $\nu_1 = 0$.

On $(i,j)$ and $\sigma$ we get on $(\nu_1 \ins \nu_n)_3 (x)$ the term
$$\epsilon(\sigma) (-1)^{\sum_{k=1}^{i} |x_{\sigma(k)}| -1} \left( (\mu^{i-1} \circ_1 \nu) \cdot (\mu^{j-2} \circ_1 \nu) \right) \, (x_{\sigma}).$$

On $(i^{\prime}, j^{\prime})$ and $\tau$ we get on $(\nu_1 \ins \nu_n)_3 (x)$ the term
$$\epsilon(\tau) (-1)^{\sum_{k=1}^{j-1} |x_{\sigma(k+i)}| -1} \left( (\mu^{j-2} \circ_1 \nu) \cdot (\mu^{i-1} \circ_1 \nu) \right) \, (x_{\tau}).$$

Let $\alpha \in \, \mathcal{S}_n$ be such that $\tau = \alpha \circ \sigma$ ($\alpha(k+i)=i$ for $k \in \, \{1, \ldots , j-1\}$ and $\alpha(k)=k+j-1$ for $k \in \, \{1, \ldots, i\}$). $\alpha$ can also be understood as the non\=/trivial element in $\mathcal{S}_2$ acting on $(x_{\sigma(1,i)}, x_{\sigma(i+1,n)})$. Thus, the last equation can be written in the following way:
$$\epsilon(\sigma)\epsilon(\alpha) (-1)^{|x_{\sigma(i+1,n)}| -1}(-1)^{(|x_{\sigma(1,i)}| -1)(|x_{\sigma(i+1,n)}| -1)} \left( (\mu^{i-1} \circ_1 \nu) \cdot (\mu^{j-2} \circ_1 \nu) \right) \, (x_{\sigma}).$$

We are going to prove that $$\epsilon(\alpha) (-1)^{|x_{\sigma(i+1,n)}| -1} (-1)^{(|x_{\sigma(1,i)}| -1)(|x_{\sigma(i+1,n)}| -1)} = -(-1)^{|x_{\sigma(1,i)}| -1}$$ and hence the terms in $\sigma$ and $\tau$ cancel each other. This comes from the fact that $\epsilon(\alpha) = (-1)^{|x_{\sigma(1,i)}||x_{\sigma(i+1,n)}|} $ and hence 
$$\epsilon(\alpha)(-1)^{|x_{\sigma(i+1,n)}| -1} (-1)^{(|x_{\sigma(1,i)}| -1)(|x_{\sigma(i+1,n)}| -1)}=(-1)^{|x_{\sigma(1,i)}|}.$$

We have shown that $\sum_{i+j=n+1} (\nu_j \ins \nu_i)_3 \, x = 0$ for all $x$ and all $n$. In order to show that $J(n) = 0$ we have to study the other case.

{\bf Case $1 = \sigma_3(1)$.}

As we said before, the only two chances for $\sigma_3^{-1}(1)$ are $1$ and $3$. This is true because $\sigma_3$ is a $(2,j-2)$ unshuffle.

Letting $\sigma_1$, $\sigma_2$  and $\sigma_3$ such that $1 = \sigma_3(1)$ varying all over the possibilities in the corresponding unshuffles gives the same that taking $\sigma \in \, \mathcal{S}{h}(2,i-2,1,j-2)$ by Remark \ref{MUS} by letting $\sigma := (\sigma_2, \sigma_3^{\prime}) \circ \sigma_1$ where $\sigma_3^{\prime}(k) := \sigma_3(k+1)$ for $k$ in $\{1, \ldots, j-2\}$. The resulting term on the sum without sign is

$$\nu((\nu(x_{\sigma(1)} \otimes x_{\sigma(2)}) \cdot x_{\sigma(3)} \cdots x_{\sigma(i)}) \otimes x_{\sigma(i+1)}) \cdot x_{\sigma(i+2)} \cdots x_{\sigma(n)}.$$

Applying Lemma \ref{GHC} to the last expression we get two different summands, which are:
$$(-1)^{|x_{\sigma(i+1)}||x_{\sigma(3,i)}|}(-1)^{|x_{\sigma(1,2)}| + 1} (\mu^{i-1} \ins f_{x_{\sigma(1,i+1)}}) \, x_{\sigma(3,i)}$$
$$\textrm{ and } (-1)^{|x_{\sigma(i+1)}||x_{\sigma(3,i)}|} (\mu^{i-1} \circ_1 (\nu \circ_1 \nu)) \, (x_{\sigma(1,2)} \otimes x_{\sigma(i+1)} \otimes x_{\sigma(3,i)}).$$

{\bf First summand.}

Let $\sigma_4 \in \, \mathcal{S}{h}_{i-3}^{1}$ and let $\sigma^{\prime} := (id^{2}, \sigma_4, id^{j-1}) \circ \sigma$. Now $\sigma^{\prime}$ runs over all the unshuffles $\mathcal{S}{h}(2,1,i-3,1,j-2)$. The resulting term without the sign is
$$\nu(x_{\sigma^{\prime}(1)} \otimes x_{\sigma^{\prime}(2)}) \cdot \nu(x_{\sigma^{\prime}(i+1)} \otimes x_{\sigma^{\prime}(3)}) \cdot x_{\sigma^{\prime}(4)} \cdots x_{\sigma^{\prime}(i)} \wedge x_{\sigma^{\prime}(i+2)} \cdots x_{\sigma^{\prime}(n)}.$$

If we permute $\nu(x_{\sigma^{\prime}(i+1)} \otimes x_{\sigma^{\prime}(3)})$ with $x_{\sigma^{\prime}(4)} \cdots x_{\sigma^{\prime}(i)}$ and we take $\sigma^{\prime \prime}$ before $\alpha$, which is again the non\=/trivial element in $\mathcal{S}_2$ acting on $((x_{\sigma^{\prime}(i+1)} \otimes x_{\sigma^{\prime}(3)}), (x_{\sigma^{\prime}(4,i)}))$ we get a sum with total sign:
$$ \sum_{\sigma^{\prime \prime} \in \mathcal{S}{h}(2,i-3,1,1,j-2)} \epsilon(\sigma^{\prime \prime}) (-1)^{\sum_{k=1}^{i-1} |x_{\sigma^{\prime \prime}(k)}| -1} \left( (\mu^{i-2} \circ_1 \nu) \cdot (\mu^{j-1} \circ_1 \nu) \right) \, (x_{\sigma^{\prime \prime}}).$$

A justification similar to the one used in last case is also possible here, nevertheless an easier justification can be found now. Using Remark \ref{MUS} we can write $\sigma^{\prime \prime} = (\tau_2, \tau_3) \circ \tau_1$ where $\tau_1 \in \, \mathcal{S}{h} (i-1, 2, j-2)$, $\tau_2 \in \, \mathcal{S}{h}_{i-3}^2$ and $\tau_3 \in \, \mathcal{S}{h}_1^1 = \mathcal{S}_2 $.  Now the previous line translates into
\begin{eqnarray*}
\sum_{\begin{smallmatrix} \tau_1 \in \, \mathcal{S}{h} (i-1, 2, j-2) \\ \tau_2 \in \, \mathcal{S}{h}_{i-3}^2 \end{smallmatrix}} &&  \epsilon(\tau^{\prime})(-1)^{\sum_{k=1}^{i-1} |x_{\tau^{\prime}(k)}| +1} (\mu^{i-2} \circ_1 \nu) \, ((x_{\tau_1(1,i-1)})_{\tau_2}) \cdot \\
&\cdot& \left( \sum_{\tau_3 \in \mathcal{S}_2} \nu \, (\tau_3 \cdot x_{\tau_1(i,i+1)}) \right) \mu^{j-2} \, (x_{\tau_1(i+2,n)}).
\end{eqnarray*}

Where $\tau^{\prime} := (\tau_2,id)\circ \tau_1 $, for every $\tau_1 \in \, \mathcal{S}{h} (i-1, 2, j-2)$. Since $\nu$ is anti\=/symmetric 
$$\sum_{\tau_3 \in \mathcal{S}_2} \nu \, (\tau_3 \cdot x_{\tau_1(i,i+1)}) = \nu(x_{\tau_1(i)} \otimes x_{\tau_1(i+1)}) - \nu(x_{\tau_1(i)} \otimes x_{\tau_1(i+1)}) = 0.$$

In this way we have shown that the first summand is zero.

{\bf Second summand.}

We are going to work now with the second summand. We have
$$ \epsilon(\sigma) (-1)^{|x_{\sigma(i+1)}||x_{\sigma(3,i)}|} \nu(\nu( x_{\sigma(1)} \otimes x_{\sigma(2)}) \otimes x_{\sigma(i+1)}) \cdot x_{\sigma(3,i)} \cdot x_{\sigma3(i+2,n)}.$$

We take again $\alpha$ the non\=/trivial element in $\mathcal{S}_2$ and we let it act on $x_{\sigma(3,i)} \otimes x_{\sigma(i+1)}$, it is clear that $\epsilon(\sigma^{\prime}) = \epsilon(\sigma) (-1)^{|x_{\sigma(i+1)}||x_{\sigma(3,i)}|}$ where $\sigma^{\prime}$ is $\sigma$ followed by the action of $\alpha$ on $x_{\sigma(3,i)} \otimes x_{\sigma(i+1)}$. $\sigma^{\prime}$ varies all over $\mathcal{S}{h}(2,1,i-2,j-2)$. We have then
$$\sum_{\sigma^{\prime} \in \mathcal{S}{h}(2,i-2,1,j-2)} \left( \mu^{n-2} \circ_1 (\nu \circ_1 \nu) \right) \, (\sigma \cdot \, x).$$

Again by Remark \ref{MUS} we can express $\sigma^{\prime}$ as the composition $(\tau_1, \tau_2) \circ \tau_3$ where $\tau_1 \in \, \mathcal{S}{h}_1^2$. In particular the previous line is equal to
$$\sum_{\begin{smallmatrix} \tau_3 \in \mathcal{S}{h}_{n-3}^{3} \\ \tau_2 \in \mathcal{S}{h}_{j-2}^{i-2} \end{smallmatrix}} \epsilon((id, \tau_2) \circ \tau_1) \left(\nu \ins \nu \, (x_{\tau_3(1)} \otimes x_{\tau_3(2)} \otimes x_{\tau_3(3)}) \right) \cdot \mu^{n-3} \, (x_{(id, \tau_2) \circ \tau_1}).$$

And hence that vanishes since $3$\=/Jacobi holds.

In this way we have shown that $\sum_{i+j=n+1} (\nu_j \ins \nu_i)_3 \, x = 0$ for all $x$ and all $n$. This completes the proof of $J(n) = 0$ for all $n$ and hence the proof of Theorem \ref{main}.
\qed


\subsection{Up to the coalgebra level}

We can view \Linfs as coderivations on the associated calgebras as we did in Proposition \ref{cochain}. In this section we compare the associated coderivation to the \Linf with higher brackets coming from a Gertenhaber algebra from Theorem \ref{main} and the one coming from the Lie algebra structure given my forgetting the multiplicative structure $\mu$. 

More precisely, given $(L, \nu, \mu)$ a Gertenhaber algebra, we can consider the coderivation associated to $\nu$. $(\smallsetminus (S L), \hat{\nu}=\hat{\nu_2})$ is a cochain complex where 
$$(\hat{\nu})_n := (\otimes^{n-1} \ins \nu) \colon (S^n L)_{-k} \subset \smallsetminus (S L)_{k} \longrightarrow (\otimes^{n-1} L)_{-(k+1)} \subset \smallsetminus (S L)_{k+1}.$$

On the other hand, applying Proposition \ref{cochain} to $(L, \{\nu_n\})$ gives a very similar codifferential:

\begin{eqnarray*}
\sum_{i=1}^{n}\mu^{n-i+1} \circ (\otimes^{n-i+1} \ins \nu_i) &=& \sum_{i=1}^{n}(\mu^{n-i+1} \ins (\mu^{i-1} \ins \nu)) = \sum_{i=1}^{n} \binom{n-2}{i-2} (\mu^{n-1} \ins \nu) \\
&=& \nu_{n} \left(\sum_{i=1}^{n} \binom{n-2}{i-2} \right) = 2^{n-2} \nu_n = 2^{n-2} \mu^{n-2} \circ (\hat{\nu})_n.
\end{eqnarray*}

In the case of the \Linf coming from a Gerstenhaber algebra we will not use the notation $d_{\textrm{CE}}$ for the codifferential, but we will specify $\hat{\nu_2}$ or $\hat{\nu}$ depending if we are only using the Lie bracket or all the other ones.


\subsection{Example: Multivector fields}

The working example of a Gerstenhaber algebra in this document will be the Gertenhaber algebra of multi-vector fields on a smooth manifold $M$ with the Schouten bracket. We review its definition in this section.

Multi\=/vector fields of degree $n$ are sections of the \nth exterior power of the tangent bundle of $M$.
$$\mathfrak{X}^{n}(M) = \Gamma(M, \wedge^n TM) \cong (\wedge_{\mathcal{C}^{\infty}(M)}^n \mathfrak{X})(M) \cong \wedge_{\mathcal{C}^{\infty}(M)}^n \mathfrak{X}(M).$$

Usually $\mathfrak{X}^{1}(M)$ will be denoted just by $\mathfrak{X}(M)$. The graded vector space given degree\=/wise by $\mathfrak{X}^{n}(M)$ is denoted by $\mathfrak{X}^{\bullet}(M)$.

\begin{df}[Algebra of multi\=/vector fields\index{Multi\=/vector field}]
	$\mathfrak{X} (M)$ is the vector space of vector fields on $M$. ($\mathfrak{X}^{\bullet}(M), \wedge)$ is a graded commutative algebra (with respect to point\=/wise exterior multiplication) which is called the algebra of multi\=/vector fields on $M$.
\end{df}

The product on the algebra will be often denoted just by $x y$ instead of $x \wedge y$. The Lie bracket between vector fields ($\mathfrak{X}(M)[1], [-,-]$) extend to a unique Lie algebra structure on $(\mathfrak{X}(M)[1])^{\bullet}$ that after the decalage isomorphism makes the triple $(\mathfrak{X}^{\bullet} (M), \nu, \wedge)$ a Gerstenhaber algebra. This bracket is called the Schouten bracket:

\begin{df}[Schouten bracket\index{Schouten bracket}]
	Given a smooth manifold $M$, the Schouten bracket on multi\=/vector fields is the unique bracket $[-,-] \in \mathrm{End}(\mathfrak{X}[1])^{\bullet}(M)$ such that the triple $(\mathfrak{X}^{\bullet}(M), \nu=\dec([-,-]), \wedge)$ is a Gerstenhaber algebra . It is given on homogeneous elements $x = (x_1 \cdots x_n) \in \, (\mathfrak{X}[1])^{n}(M)$ and $y = (y_1 \cdots y_m)$ in $(\mathfrak{X}[1])^{m}(M)$ by the explicit formula:
	\begin{eqnarray*}
		[x,y] &:=& \sum_{\begin{smallmatrix} \sigma \in \mathcal{S}h_{n-2}^2 \\ \tau \in \mathcal{S}h_{m-2}^2 \end{smallmatrix}}
		\epsilon(\sigma) \epsilon(\tau) [x_{\sigma(1)},y_{\tau(1)}] x_{\sigma(2)} \cdots x_{\sigma(n)} y_{\tau(2)} \cdots y_{\sigma(m)}\\
		&=& \sum_{\begin{smallmatrix} 1 \leqslant i \leqslant n \\ 1 \leqslant j \leqslant n \end{smallmatrix}}
		(-1)^{i+j} [x_i,y_j] x_1 \cdots x_{i-1} x_{i+1} \cdots x_n y_1 \cdots y_{j-1} y_{j+1} \cdots y_m.
	\end{eqnarray*}
\end{df}

Multi-vector fields are equipped with a Gerstenhaber algebra structure given by $\left( \mathfrak{X}^{\bullet}(M), \wedge, [-,-] \right)$. Theorem \ref{main} can be applied to this Gerstenhaber algebra.

\subsubsection{Multi\=/vector calculus}

There are natural extensions to multi\=/vector fields of both the concept of contraction on differential forms and the Lie derivative along a vector field. We need to observe first that the algebra of multi\=/vector fields is similar to the algebra of differential forms.  Differential forms of degree $n$ are sections of the \nth exterior power of the co\=/tangent bundle of $M$ and they are denoted by $\Omega^{n} (M)$.

\begin{df}[Algebra of differential forms\index{Differential form}]
($\Omega^{\bullet}(M), \wedge)$ is a graded commutative algebra with the grading given by the degree of the form.
\end{df}

The differential $d \colon \Omega^{\bullet}(M) \longrightarrow \Omega^{\bullet +1}(M)$ is going to be understood as an element of $\underline{\textrm{End}}(\Omega^{\bullet}(M))_{1}$. Thus, $(\Omega^{\bullet}(M), d)$ is in fact a cochain complex since $d \cdot d = 0$. Here we are using $(\cdot)$ to denote the product on $\underline{\textrm{End}}(\Omega^{\bullet}(M))$. 

We recall that the complex $(\Omega^{\bullet}(M), d)$ is called the de Rham complex. The co\=/cycles are called closed forms and the co\=/boundaries, exact forms. 

We need to introduce the concept of the {\it total Koszul sign} on the graded tensor algebra of a graded algebra.

\begin{df}[Total Koszul sign\index{Koszul sign!total Koszul sign}]\label{TKS}
Let $A$ be a commutative algebra. For every positive integer $n$ and every $a := a_1 \otimes \cdots \otimes a_n \in \, T_n(A)$ we define the total Koszul sign of $a$ as the integer $\epsilon(a) \in \, \{-1, 1\}$ such that:
$$a_1 \otimes \cdots \otimes a_n = \epsilon(a) a_n \otimes \cdots \otimes a_1.$$
\end{df}

The total Koszul sign is given by the explicit formula $\epsilon(a) = (-1)^{\sum_{i=1}^{n}\sum_{j=1}^{i-1}|a_i||a_j|}$. In the particular case where $|a_i|=1$ for every $i \in \, \{1, \ldots, n\}$ the total Koszul sign is simply $\epsilon(a) = (-1)^{\frac{n(n-1)}{2}}$. 

The reason why we have introduced the total Koszul sign is that we are not going to use the definition of the contraction operation mostly found on the literature.

\begin{df}[Contraction operator\index{Contraction}]\label{CO}
Every $x=x_1 \cdots x_n \in \mathfrak{X}^{n}(M)$ defines an endomorphism $\iota_x \in \underline{\textrm{End}}(\Omega^{\bullet}(M))_{-n}$ which is given for every $\alpha \in \, \Omega^{\bullet}(M)$ by:
$$\iota_{x_1 \cdots x_n} \alpha := \alpha(x_n, \ldots, x_1, -).$$ 
\end{df}

The contraction operator in multi\=/vector fields is a generalization to the usual contraction operator of vector fields which is defined to be $\iota_{x_1} \alpha = \alpha(x_1, -)$. It is extended to multi\=/vector fields by setting $\iota_{x_1 \cdots x_n} = \iota_{x_1} \cdots \iota_{x_n}$. Usually the contraction operator is defined to be $\hat{\iota}_{x_1 \cdots x_n} \alpha := \alpha(x_1, \ldots, x_n, -)$ as in \cite{FPR}, that convention follows the spirit $\hat{\iota}_{x_1 \cdots x_n} = \hat{\iota}_{x_n} \cdots \hat{\iota}_{x_1}$.

Our convention makes $\iota$ into a morphism of graded commutative algebras from $\smallsetminus (\mathfrak{X}^{\bullet}(M))$ to $\underline{\textrm{End}}(\Omega^{\bullet}(M))$. Clearly, the relation between the two definitions is the total Koszul sign, for $x \in \, \mathfrak{X}^{n}(M)$ we have:
$$ \iota_{x} = (-1)^{\frac{n(n-1)}{2}} \hat{\iota}_x.$$

We have decided to take the convention adopted in Definition \ref{CO} because in this way the contraction operator is a cochain map between multi\=/vector fields and differential forms as we will see in the following section.

We continue the generalization of the operations involving vector fields to multi\=/vector fields. Now we extend the notion of the Lie derivative to multi\=/vector fields just by requiring the Cartan's formula to hold. 

\begin{df}[Lie derivative\index{Lie derivative}]
We define the Lie derivative along a multi\=/vector field $x \in \, \mathfrak{X}^{ n}(M)$ to be:
$$\mathcal{L}_x := [d,\iota_x]_{\textrm{RN}} = [d,\iota_x] \in \, \underline{\textrm{End}}(\Omega^{\bullet}(M))_{1-n}.$$
\end{df}

We warn the reader again about our convention in the definition of the contraction operation which differs from the usual one by the total Koszul sign. This difference remains in the definition of the Lie derivative. The usual notion of Lie derivative $\widehat{\mathcal{L}}_x = [d, \hat{\iota}_x]$ is related to ours by the total Koszul sign 
$$\widehat{\mathcal{L}}_x = [d, \hat{\iota}_x] = (-1)^{\frac{n(n-1)}{2}} [d, \iota_x] = (-1)^{\frac{n(n-1)}{2}} \mathcal{L}_x.$$


\subsubsection{Cartan Calculus}

We have defined several endomorphisms of $\Omega^{\bullet}(M)$: the contraction operator, the Lie derivative and the usual de Rahm differential. The study of the commutativity of those operations is usually called Cartan calculus. In the following Proposition, we compute explicitly the relations between those operators with our new convention for the insertion of multi-vector fields. The result for the usual convention is common in the literature (see \cite[Proposition A.3]{FPR} for instance).

\begin{pp}[Cartan calculus\index{Cartan calculus}]\label{CCs}
Let $x$ and $y$ be two multi\=/vector fields. Then
\begin{enumerate}
	\item $[\iota_x, \iota_y] = 0$.
	\item $[d, \mathcal{L}_x] = 0$.
	\item $[\mathcal{L}_x, \iota_y] = \iota_{\nu(x,y)}$.
	\item $[\mathcal{L}_x, \mathcal{L}_y] = (-1)^{|x|+1} \mathcal{L}_{\nu(x,y)}$.
	\item $\mathcal{L}_{x \otimes y} = \mathcal{L}_x \iota_y + (-1)^{|x|} \iota_x \mathcal{L}_y$.
\end{enumerate}
\end{pp}

\dem
\begin{enumerate}
	\item From the definition of the contraction operation 
	$$\iota_x \iota_y = \iota_{x y} = (-1)^{|x||y|} \iota_{y x} = (-1)^{|x||y|} \iota_y \iota_x.$$ Now it is clear that $[\iota_x, \iota_y]= \iota_x \iota_y -(-1)^{|x||y|} \iota_y \iota_x = \iota_x \iota_y -\iota_x \iota_y = 0$.
	\item $[d, \mathcal{L}_x] = [d, [d, \iota_x]] = [d, d\iota_x -(-1)^{|x|} \iota_x d] =-(-1)^{|x|} d \iota_x d -(-1)^{|x|-1} d\iota_x d = 0.$
	\item This result when $|x|=1=|y|$ reads as follows: 
	$$[\widehat{\mathcal{L}}_x, \hat{\iota}_y] = [\mathcal{L}_x, \iota_y] = \iota_{\nu(x,y)} = \hat{\iota}_{[x,y]} = \hat{\iota}_{[x,y]}.$$
	That result is well known for vector fields (see \cite[section 18.3]{C} for instance), and we only need to show it by induction on $n=|x|$ and $m=|y|$. First we assume that the result holds for $n$ and $m$ and we want to show it for $m+1$. We will later deal with the step from $n$ to $n=1$. 

We decompose $y = y_1 y_2$ where $y_2 \in \, \mathfrak{X}(M)$ and $y_1 \in \, \mathfrak{X}^{ m}(M)$. Applying the Leibniz rule we get $$\nu(x,y_1 y_2) = \nu(x,y_1)y_2 + (-1)^{(n-1)m}y_1\nu(x,y_2).$$ Now we see that:
		\begin{eqnarray*}
			\iota_{\nu(x,y)} &=& \iota_{\nu(x,y_1)}\iota_{y_2} + (-1)^{(n-1)m}\iota_{y_1}\iota_{\nu(x,y_2)} \\
			&=& [\mathcal{L}_x, \iota_{y_1}]\iota_{y_2} +(-1)^{(n-1)m}\iota_{y_1}[\mathcal{L}_x, \iota_{y_2}] \\
			&=& \mathcal{L}_x \iota_{y} -(-1+1)(-1)^{(n-1)m}\iota_{y_2}\mathcal{L}_x\iota_{y_1} -(-1)^{(n-1)(m+1)}\iota_{y}\mathcal{L}_x \\
			&=& [\mathcal{L}_x, \iota_y].
		\end{eqnarray*}
		To show the result for $n+1$ we use the fact that the bracket on $\underline{\textrm{End}}(\Omega^{\bullet}(M))$ is a Lie bracket. Let $x=x_1 x_2$ where $x_2 \in \, \mathfrak{X}(M)$ and $x_1 \in \, \mathfrak{X}^{ n}(M)$.
		\begin{eqnarray*}
			[\mathcal{L}_{x}, \iota_y] &=& [[d,\iota_x], \iota_y] = (-1)^{mn +1} [\iota_y, [d, \iota_x]] = \\
			&=& (-1)^{m}\left((-1)^{m} [d, [\iota_x, \iota_y]] + (-1)^{n+1} [\iota_x, [\iota_y, d]]\right) \\
			&=& (-1)^{n} [\iota_x, \mathcal{L}_y] = (-1)^{m(n+1)} [\mathcal{L}_y, \iota_x] \\
			&=& (-1)^{m(n+1)} \iota_{\nu(y,x)} = \iota_{\nu(x,y)}. 
		\end{eqnarray*}
	\item Using the previous results and the Jacobi identity:
		\begin{eqnarray*}
			[\mathcal{L}_x, \mathcal{L}_y] &=& [\mathcal{L}_x, [d, \iota_y]] = \\
			&=& (-1)^{(|x|-1)|y| +1}\left( (-1)^{|x|-1}[d,[\iota_y, \mathcal{L}_x]] + (-1)^{|y|}[\iota_y,[\mathcal{L}_x, d]] \right) \\
			&=& (-1)^{|x||y| + |y| + |x|}[d, [\iota_y, \mathcal{L}_x]] = (-1)^{|x| + 1} [d, [\mathcal{L}_x, \iota_y]] \\
			&=& (-1)^{|x|+1} [d, \iota_{\nu(x,y)}] = (-1)^{|x|+1} \mathcal{L}_{\nu(x,y)}.
		\end{eqnarray*}
	\item A straight forward calculation shows that:
				\begin{eqnarray*}
				\mathcal{L}_{xy} &=& d\iota_{xy} -(-1)^{n+m} \iota_{xy}d \\
				&=& d\iota_x \iota_y -(-1)^{n+m} \iota_x \iota_y d + (-1)^{n} \iota_x d \iota_y - (-1)^{n} \iota_x d \iota_y \\
				&=& [d, \iota_x] \iota_y + (-1)^{n} \iota_x [d, \iota_y] = \mathcal{L}_x \iota_y + (-1)^{n} \iota_x \mathcal{L}_y.
				\end{eqnarray*}
\end{enumerate}
\qed

It is possible to extend the definition of the contraction operation and the Lie derivative to $T (\mathfrak{X}^{ \bullet}(M))$ by regarding $x_1 \cdots x_n \in \, \mathfrak{X}^{ j_1}(M) \otimes \cdots \otimes \mathfrak{X}^{ j_n}(M)$ as an element in $\mathfrak{X}^{ \sum_{i=1}^n j_i}(M)$. 

We know that $\mathcal{L}_{x_1 \otimes x_2} = \mathcal{L}_{x_1} \iota_{x_2} + (-1)^{|x_1|} \iota_{x_1} \mathcal{L}_{x_2}$. We would like to know what happens about $\mathcal{L}_{x_1 \cdots x_n}$ for a general $n$. The next result answer this question and is a generalization to multi\=/vector fields of \cite[Lemma 3.2.1]{FRS}.

\begin{pp}\label{FR}
Let $x = x_1 \cdots x_n$ be a tensor product of $n$ multi\=/vector fields. Then
$$\mathcal{L}_x - \iota_{\nu_n(x_1, \ldots, x_n)} = (-1)^{|x|}\sum_{\sigma \in \mathcal{S}h_1^{n-1}} \epsilon(\sigma)(-1)^{|x_{\sigma(n)}|} \iota_{x_{\sigma(1,n-1)}} \mathcal{L}_{x_{\sigma(n)}}.$$ 
\end{pp}

We recall that the notation $x_{\sigma(1,n-1)}$ means $x_{\sigma(1)} \cdots x_{\sigma(n-1)}$.

\dem
The proof will be given by induction on $n$. For $n=1$ we have that $\nu_1 = 0$ and the statement of the proposition is just $\mathcal{L}_x = \mathcal{L}_x$.

For the induction step we use Proposition \ref{CCs}, the numbers above the equalities in the following lines refer to the different statements in Proposition \ref{CCs}. We assume the statement is true for $n$ and we are going to show it for $(n+1)$:
\begin{eqnarray*}
\mathcal{L}_x &=& \mathcal{L}_{x_{1,n}\otimes x_{n+1}} \stackrel{\scriptscriptstyle{(5)}}{=} \mathcal{L}_{x_{1,n}} \iota_{x_{n+1}} + (-1)^{|x|-|x_{n+1}|}\iota_{x_{1,n}}\mathcal{L}_{x_{n+1}} \\
&=& \iota_{\nu_{n}(x_1,\ldots, x_n)} \iota_{x_{n+1}} + (-1)^{|x|-|x_{n+1}|}\sum_{\sigma \in \mathcal{S}h_1^{n-1}} \epsilon(\sigma)(-1)^{|x_{\sigma(n)}|} \iota_{x_{\sigma(1,n-1)}} \mathcal{L}_{x_{\sigma(n)}} \iota_{x_{n+1}} \\
&+& (-1)^{|x| - |x_{n+1}|}\iota_{x_{1,n}}\mathcal{L}_{x_{n+1}} \\
&\stackrel{\scriptscriptstyle{(3)}}{=}& \iota_{\nu_{n}(x_1,\ldots, x_n)} \iota_{x_{n+1}} + (-1)^{|x| - |x_{n+1}|}\iota_{x_{1,n}}\mathcal{L}_{x_{n+1}} + (-1)^{|x|-|x_{n+1}|} \cdot \\
&\cdot&\sum_{\sigma \in \mathcal{S}h_1^{n-1}} \epsilon(\sigma)(-1)^{|x_{\sigma(n)}|} \iota_{x_{\sigma(1,n-1)}} \left( \iota_{\nu(x_{\sigma(n)}, x_{n+1})} + (-1)^{(|x_{\sigma(n)}|-1)|x_{n+1}|} \iota_{x_{n+1}} \mathcal{L}_{x_{\sigma(n)}} \right)\\
&=& \iota_{\nu_{n+1}(x_1,\ldots, x_{n+1})} \iota_{x_{n+1}} + (-1)^{|x|}\sum_{\sigma \in \mathcal{S}h_1^{n}} \epsilon(\sigma)(-1)^{|x_{\sigma(n+1)}|} \iota_{x_{\sigma(1,n)}} \mathcal{L}_{x_{\sigma(n+1)}}. 
\end{eqnarray*}
\qed


\color{BlueViolet}
\section{Multi\=/symplectic geometry}
\color{black}


We review the basic definitions on multi-symplectic geometry and concentrate on what happens at the level of multi-symplectic and hamiltonian multi-vector fields 

In the second subsection we conclude from one side that multi-symplectic multi-vector fields and hamiltonian multi-vector fields are \mbox{$L_{\infty}$}\=/subalgebras of the \Linf of multi-vector fields with non-trivial higher brackets from Theorem \ref{main}. On the other side, we interpret the contraction of multi-symplectic multi-vector fields against the pre-multisymplectic form as a cochain map, justifying the choice of convention for the insertion of multi-vector fields.

\subsection{Basic definitions on multi\=/symplectic geometry}

Here we include the basic definitions on multi-symplectic geometry, present for example in the work of Rogers \cite{R}: what a multi\=/symplectic manifold is, what a multi\=/symplectic or hamiltonian multi\=/vector field is and what a hamiltonian form is. All those concepts are generalizations of the corresponding ones in symplectic geometry.

\begin{df}[Pre\=/multi\=/symplectic manifold\index{Multi\=/symplectic!Pre\=/multi\=/symplectic manifold}]
A differential form $\omega$ of degree $(n+1)$ on a manifold $M$ is said to be pre\=/$n$\=/symplectic if it is closed. The pair $(M,\omega)$ is called a pre\=/$n$\=/symplectic manifold.
\end{df}

We recover a pre\=/symplectic manifold if $n=1$. The step from pre\=/multi\=/symplectic to multi\=/symplectic structures is like the one from pre\=/symplectic to multi\=/symplectic manifolds.

\begin{df}[Multi\=/symplectic manifold\index{Multi\=/symplectic!manifold}]
A form $\omega \in \Omega^{n+1}(M)$ is called non-degenerate if $\widetilde{\omega} : \mathfrak{X}(M) \longrightarrow \Omega^{n}(M)$ given by $\widetilde{\omega}(x):=\iota_x \omega$ is injective. A form $\omega$ in $\Omega^{n+1}(M)$ is called $n$\=/symplectic if it is both closed and non\=/degenerate. The pair $(M, \omega)$ is then called an $n$\=/symplectic manifold.
\end{df}

Given $(M,\omega)$ a $n$\=/symplectic manifold, $\widetilde{\omega}$ might be extended to 
$$\widetilde{\omega} \colon \mathfrak{X}^{ m}(M) \longrightarrow \Omega^{n+1-m}(M)$$
by the same formula ($\widetilde{\omega}(x):=\iota_x \omega$) for every $m \in \, \{1, \ldots, n+1\}$.

It is possible to extend the definition of symplectic and hamiltonian vector fields to the multi\=/symplectic case. This step is crucial in our approach to the problem of studying the generalization of momentum maps to the multi\=/symplectic world. In \cite{FRZ} the study of purely hamiltonian vector fields is used to extend the concept of momentum maps. In the following chapters we will use hamiltonian and multi\=/symplectic multi\=/vector fields to get to a variety of definitions of momentum maps.

\begin{df}[Multi\=/symplectic multi\=/vector field\index{Multi\=/symplectic!multi\=/vector field}\index{Multi\=/vector field!multi\=/symplectic}]
Let $(M,\omega)$ be a pre- $n$\=/symplectic manifold. A multi\=/vector field $x\in\, \mathfrak{X}^{ \bullet}(M)$ is called $n$\=/symplectic if $\mathcal{L}_{x}\omega=0$.
\end{df}

Since $\omega$ is closed $\mathcal{L}_x \omega = d\iota_x \omega$, and hence $n$\=/symplectic multi\=/vector fields are precisely those for which $\widetilde{\omega}(x)$ is closed. Non\=/trivial $n$\=/symplectic multi\=/vector fields only exist in degrees $\{1,\ldots, n+1\}$. $1$\=/symplectic vector fields are the symplectic vector fields.

The set of $n$\=/symplectic $m$\=/vector fields is a vector sub\=/space of $\mathfrak{X}^{m}(M)$ and it is denoted by $\mathfrak{X}_{\textrm{sym}}^{m}(M)$. The  graded vector space of multi\=/symplectic multi\=/vector fields is denoted by $\mathfrak{X}_{\textrm{sym}}^{ \bullet}(M)$.

\begin{df}[Hamiltonian multi\=/vector fields\index{Multi\=/vector field!hamiltonian}\index{Hamiltonian!multi\=/vector field}]
A multi\=/vector field $x$ in $\mathfrak{X}^{ m}(M)$, on a given a pre\=/$n$\=/symplectic manifold $(M, \omega)$, is called hamiltonian when there exists $\alpha \in \, \Omega^{n-m}(M)$ such that $\iota_{x}\omega=d\alpha$.
\end{df}

Hamiltonian multi\=/vector fields are precisely those such that $\widetilde{\omega}(x)$ is exact. Hence hamiltonian multi\=/vector fields are in particular multi\=/symplectic. Non\=/trivial hamiltonian multi\=/vector fields only exist in degrees $\{1,\ldots, n\}$ on a pre\=/$n$\=/symplectic manifold. Hamiltonian vector fields in a $1$\=/symplectic manifold correspond to hamiltonian vector fields in symplectic geometry.

The set of hamiltonian $m$\=/vector fields is a vector sub\=/space of $\mathfrak{X}^{m}(M)$ and it is denoted by $\mathfrak{X}_{\textrm{ham}}^{m}(M)$. The  graded vector space of hamiltonian multi\=/vector fields is denoted by $\mathfrak{X}_{\textrm{ham}}^{ \bullet}(M)$.

The \mbox{$4\mbox{th}$} statement in Proposition \ref{CCs} tells us that $\mathfrak{X}_{\textrm{sym}}^{ \bullet}(M)$ is closed under $\nu$, and hence it is a Lie sub\=/algebra of $(\mathfrak{X}^{ \bullet}(M), \nu)$. 

\begin{df}[Hamiltonian forms and pairs]\index{Hamiltonian!form}\index{Hamiltonian!pair}
Let $(M, \omega)$ be a pre\=/$n$\=/symplectic manifold. A differential form $\alpha$ on $M$ is called hamiltonian if there exists a hamiltonian multi\=/vector field such that $\iota_x \omega = d \alpha$. The pair $(x,\alpha)$ is called a hamiltonian pair\index{Hamiltonian!pair}.
\end{df}

Hamiltonian forms form a sub\=/cochain complex of the de Rham complex on $M$ denoted by $\Omega_{\mathrm{ham}}(M)$. This is true because closed forms (in particular exact ones) are trivially hamiltonian. Observe that non\=/closed hamiltonian forms only exists in degrees $\{0, \ldots, n-1\}$: if $\alpha \in \Omega_{\textrm{ham}}^{i \geqslant n}(M)$ then $d \, \alpha = \iota_x \omega \in \, \Omega^{i+1 \geqslant n+1}(M)$ which means that $x \in \mathfrak{X}_{\textrm{ham}}^{n-i \leqslant 0}(M)=0$ so that $\iota_{0} \omega = 0 = d\, \alpha$.


\subsection{Contraction against the pre-multisymplectic form is a cochain map}

The vector space of hamiltonian pairs of $m$\=/vector fields and $(n-m)$\=/forms is denoted by $\widetilde{\mathfrak{X}}^{ m}(M)$. The  graded vector space of hamiltonian pairs is then $\widetilde{\mathfrak{X}}^{ \bullet}(M)$. Again, the non\=/trivial part of this graded vector space is in degrees $\{1, \ldots, n \}$.

The pairing $\widetilde{\omega} \colon \mathfrak{X}^{ m}(M) \longrightarrow \Omega^{n+1-m}(M)$ is linear, but it does not define a morphism of \Linfs or of co\=/chain complexes by degree reasons. We fix that by viewing $\Omega^{n+1-m} (M) = \Omega[n+1] (M)_{-m} = \smallsetminus(\Omega[n+1] (M))_m$.

We are going to use the notation $\smallsetminus[n+1]\Omega(M)$ to refer to the graded vector space given by $\smallsetminus \left( \Omega[n+1](M) \right)$ without the brackets, to make the reading lighter. For a general graded vector space $V$ we denote by $\smallsetminus [n]V$ the graded vector space given by $\smallsetminus(V[n])$.

Proposition \ref{FR} applied to a multi\=/symplectic form gives interesting conditions.

\begin{pp}\label{maincl}
Let $(M, \omega)$ be a pre\=/$m$\=/symplectic manifold. Let $x = x_1 \cdots x_n$ be a exterior power of $n$ multi\=/vector fields. Then
$$d \iota_x \omega - \iota_{\nu_n(x_1,\ldots, x_n)} \omega = (-1)^{|x|}\sum_{\sigma \in \mathcal{S}h_1^{n-1}} \epsilon(\sigma)(-1)^{|x_{\sigma(1)}|} \iota_{x_{\sigma(2,n)}} d \iota_{x_{\sigma(1)}} \omega.$$
Moreover, if $x_i$ is multi\=/symplectic for each $i \in \, \{1, \ldots, n\}$ then
$$d \iota_x \omega - \iota_{\nu_n(x_1,\ldots, x_n)} \omega = 0.$$
\end{pp}

\dem We only need to observe that if $\omega$ is closed $\mathcal{L}_{y} \omega = d \iota_y \omega$ for every $y$ multi\=/vector field and $d \iota_{y} \omega = 0$ for every $y$ multi\=/symplectic multi\=/vector field.
\qed\\

This result has several important consequences. First, the bracket $\nu$ of two multi\=/symplectic multi\=/vector fields is not only multi\=/symplectic, but hamiltonian. It follows that \ham is closed under $\nu$ and then it defines a Lie algebra. This extends to an arbitrary number of multi\=/symplectic multi\=/vector fields, its bracket is hamiltonian. 

Proposition \ref{cochain} might be applied both to the Lie sub\=/algebras of multi\=/symplectic and hamiltonian multi\=/vector fields giving the following cochain complexes.

\begin{pp}
Let $(M, \omega)$ be a pre\=/multi\=/symplectic manifold. Then the following sequences define cochain complexes.
$$\cdots \longrightarrow (S \ham)_{-(m)} \stackrel{\scriptscriptstyle{(\hat{\nu_2})}}{\longrightarrow} (S \ham)_{-(m+1)} \rightarrow \cdots ,$$
$$\colon \cdots \longrightarrow (S \mathfrak{X}_{\textrm{sym}}^{ \bullet} (M))_{-(m)} \stackrel{\scriptscriptstyle{(\hat{\nu_2})}}{\longrightarrow} (S \mathfrak{X}_{\textrm{sym}}^{ \bullet} (M))_{-(m+1)} \rightarrow \cdots ,$$
$$ \cdots \longrightarrow S^{-m} (\ham) \stackrel{\scriptscriptstyle{(\hat{\nu_2})_{-m}}}{\longrightarrow} S^{-(m+1)} (\ham) \rightarrow \cdots \textrm{ and }$$
$$\cdots \longrightarrow S^{-m} (\mathfrak{X}_{\textrm{sym}}^{ \bullet} (M)) \stackrel{\scriptscriptstyle{(\hat{\nu_2})_{-m}}}{\longrightarrow} S^{-(m+1)} (\mathfrak{X}_{\textrm{sym}}^{ \bullet} (M)) \rightarrow \cdots .$$
\end{pp}

It is clear that $(\smallsetminus (S \ham), \hat{\nu_2})$ is a cochain sub\=/complex of $(\smallsetminus (S \mathfrak{X}_ {\mathrm{sym}}^{\bullet}(M)), \hat{\nu_2})$.

We recall that $(\hat{\nu_2})_n := \otimes^{n-1} \ins \nu$ and so $\nu_n = \mu^{n-1} \circ (\hat{\nu_2})_n$. It is clear that $\widetilde{\omega} = \widetilde{\omega} \circ \mu^{\bullet}$ because $\iota_{x \otimes y} = \iota_x \iota_y$. This implies that $\widetilde{\omega}(\hat{\nu_2}(x)) = \widetilde{\omega}(\nu_{n}(x))$ for all $x \in \, S^{n}(\mathfrak{X}^{ \bullet}(M))$.

One of the main results of this section, which will be very important in the future development of the theory, is the following:

\begin{pp}\label{ccmap}
Let $(M, \omega)$ be a pre\=/$m$\=/symplectic manifold. Then 
$$\widetilde{\omega} \colon (\smallsetminus[m+1](S \ham), \hat{\nu_2}) \longrightarrow (\Omega^{\bullet}(M), d) \textrm{ and }$$
$$\widetilde{\omega} \colon (\smallsetminus[m+1](S \mathfrak{X}_{\textrm{sym}}^{\bullet}(M)), \hat{\nu_2}) \longrightarrow (\Omega^{\bullet}(M), d)$$
are cochain maps. In fact, since there only exists positive $(\mathbb{N} \cup \{0\})$ degree forms and there not exist multi\=/vector fields of degree less than $0$, we can restict the cochain maps to:
$$\widetilde{\omega} \colon \mathrm{tr}_{m}(\smallsetminus[m+1](S \ham), \hat{\nu_2}) \longrightarrow \mathrm{tr}_{m}(\Omega^{\bullet}(M), d) \textrm{ and }$$
$$\widetilde{\omega} \colon \mathrm{tr}_{m}(\smallsetminus[m+1](S \mathfrak{X}_{\textrm{sym}}^{\bullet}(M)), \hat{\nu_2}) \longrightarrow \mathrm{tr}_{m}(\Omega^{\bullet}(M), d).$$
\end{pp}

\dem
The map $\widetilde{\omega}$ is linear. We fix $x \in \, (S^{n} (\mathfrak{X}_{sym}^{ \bullet}(M)))_{m+1-k} \subset \smallsetminus[m+1](S \mathfrak{X}^{\bullet}(M))_k$. Now $\widetilde{\omega}(x) \in \, \Omega^{m+1-(m+1-k)}(M) = \Omega^k M$, so that the map $\widetilde{\omega}$ has the appropriate degree. We conclude the first statement by Proposition \ref{maincl} 
$$\widetilde{\omega} \hat{\nu_2} (x) = \iota_{\nu_n(x)} \omega = d \iota_x \omega = d \widetilde{\omega} (x).$$
This proves both that $\widetilde{\omega}$ is a cochain map and that $\widetilde{\omega}(x)$ is a hamiltonian form.
\qed\\

This result gives a justification of some of the conventions that we have been adopting along the chapters. First, it emphasizes the fact that we selected the correct sign for $\hat{\nu_2}$. Second, we get a justification of why we defined the contraction operation the way we did. $\hat{\iota}_{\nu_n(x)} = (-1)^{|x|} d \hat{\iota}_{x}$ does not give a cochain map.

We have already seen that multi\=/symplectic and hamiltonian multi\=/vector fields form a Lie sub\=/algebra of the Lie algebra of multi\=/vector fields furnished with the Schouten bracket. We are not able to apply Theorem \ref{main} to neither multi\=/symplectic nor hamiltonian multi\=/vector fields since we have not checked that any of those Lie sub\=/algebras are closed under the product, and thus they do not need to define a Gerstenhaber algebra. Nevertheless, both multi\=/symplectic and hamiltonian multi\=/vector fields are \mbox{$L_{\infty}$}\=/sub\=/algebras of the \Linf of multi\=/vector fields.

\begin{pp}\label{sub}
Let $(M, \omega)$ be a pre\=/multi\=/symplectic manifold. Then $\mathfrak{X}_{\textrm{sym}}^{ \bullet}(M)$ and \ham are closed under $\nu_n$ for all $n\geqslant 1$ and hence they define \mbox{$L_{\infty}$}\=/sub\=/algebras of the \Linf of multi\=/vector fields on $M$.
\end{pp}

\dem
From Proposition \ref{maincl} we know that $\iota_{\nu(x)}\omega = d \iota_x \omega$ for all $x$ which are exterior power of multi\=/symplectic multi\=/vector fields. This shows that $\nu(x)$ is a hamiltonian multi\=/vector field. In particular $\nu(x)$ is multi\=/symplectic. If the multi\=/vector fields were hamiltonian, then they are multi\=/symplectic and hence $\nu(x)$ is hamiltonian.
\qed


\color{BlueViolet}
\section{\Linfsb on Hamiltonian Forms}
\color{black}

In symplectic geometry there is a Lie bracket, the Poisson bracket defined on hamiltonian forms (only hamiltonian functions exist). If $(M, \omega)$ is a pre\=/symplectic manifold, the bracket is given by $\{f, g\} = \hat{\iota}_{x_f x_g}\omega = - \iota_{x_f x_g} \omega$ where $(x_f, f)$ and $(x_g, g)$ are hamiltonian pairs. Due to the fact that $\omega$ is anti\=/symmetric, the bracket is anti\=/symmetric as well. 

If $(M, \omega)$ is now a pre\=/$m$\=/symplectic manifold, $\omega$ is (graded\=/) symmetric making $\{\alpha, \beta\}:=-\iota_{x_{\alpha} x_{\beta}} \omega$ on $\smallsetminus[m]\Omega_{ham}(M)$ symmetric (where $(x_{\alpha}, \alpha)$ and $(x_{\beta}, \beta)$ are hamiltonian pairs). Even more, $\{-,-\}$ respect the hamiltonian character, since $d \{\alpha, \beta \}= -d \iota_{x_{\alpha} x_{\beta}} \omega = \iota_{-\nu(x_{\alpha} x_{\beta})} \omega$, so that $(-\nu(x_{\alpha}, x_{\beta}), \{\alpha, \beta\})$ is a hamiltonian pair. This structure, by contrast, does not satisfy the $3$\=/Jacobi identity. That is, given $(x_i, \alpha_i)$ hamiltonian pairs for $i \in \, \{1,2,3\}$:
\begin{eqnarray*}
	\sum_{\sigma \in \mathcal{S}h_1^2} \epsilon(\sigma) \{\{\alpha_{\sigma(1)}, \alpha_{\sigma(2)}\}, \alpha_{\sigma(3)}\}
	&=& \sum_{\sigma \in \mathcal{S}h_1^2} \epsilon(\sigma) -\iota_{-\nu(x_{\sigma(1)} x_{\sigma(2)}) x_{\sigma(3)}} \omega \\ 
	= \iota_{\sum_{\sigma \in \mathcal{S}h_1^2} \epsilon(\sigma) \nu(x_{\sigma(1)} x_{\sigma(2)}) x_{\sigma(3)}} \omega
	&=& \iota_{\nu_3(x_1 x_2 x_3)} \omega \, = \, d \iota_{x_1 x_2 x_3} \omega.
\end{eqnarray*}

This topic has already been discussed by Rogers in \cite{R}. Rogers developed an \Linf on forms that only considers hamiltonian forms on top degree (bottom degree in our convention of $\smallsetminus[m]\Omega(M)$).

In this section, we extend the \Linf of Rogers to include lower degree hamiltonian forms (those associated with higher hamiltonian multi-vector fields). In the first subsection we study the naive analogue to Rogers' result, where we restrict all forms to be hamiltonian and extend all brackets to be defined everywhere. That naive family of brackets does not define an \Linf structure on $\smallsetminus[m]\Omega_{\textrm{ham}}(M)$.

On the following subsection we do produce a higher analogue including hamiltonian forms of all degrees. Later we relate that \Linf to the one of hamiltonian multi-vector fields.

\subsection{Naive family of brackets}

 This subsection tries to justify why the naive \Linf structure on hamiltonian forms does not work and why Rogers' \Linf is the closest one to this naive structure. In the next section we give an alternative construction on an \Linf that on one hand includes hamiltonian forms of all degrees, but on the other hand is much bigger and less intuitive. 

The naive structure to be considered has as graded vector space $\smallsetminus[m]\Omega_{\textrm{ham}}(M)$ and as maps $l_{n\geqslant 2}(\alpha_1 \otimes \cdots \otimes \alpha_n) := -\{\alpha_1, \ldots, \alpha_n\}_n=\iota_{x_1 \cdots x_n} \omega$ where $(x_i, \alpha_i)$ are hamiltonian pairs for all $i \in \, \{1, \ldots, n\}$.  As a first remark, observe that the definition of $l_n$ is independent on the hamiltonian multi\=/vector field chosen. If $\iota_{x_n}\omega=\iota_{y_n}\omega=d\alpha_{n}$ then 
$$\iota_{x_1 \cdots x_n} \omega = \iota_{x_1 \cdots x_{n-1}} \iota_{x_n} \omega = \iota_{x_1 \cdots x_{n-1}} \iota_{y_n} \omega = \iota_{x_1 \cdots x_{n-1} y_n} \omega.$$

If we take $l_1 \alpha_1 = \iota_{x_1} \omega = d \alpha_1$ we will get the usual de Rham differential. This is not a good idea as we will see in the upcoming paragraphs. As we want the pair $(\nu_n, l_n)$ to define a hamiltonian pair we can take $l_1 := -d$ since $d \circ l_1 \alpha_1 = -d \circ d \alpha_1 = 0 = \iota_{0} \omega = \iota_{\nu_1 x_1} \omega$. That is why $l_{n} \ins l_1 = 0$, which is clear for $n=1$ and for a larger $n$ as we have said $(0, d \alpha)$ is a hamiltonian pair and because $l_n$ does not depend on the hamiltonian multi\=/vector field chosen for $\alpha$ we can take $x = 0$. 

It is also clear that the maps $l_n$ are morphisms of degree $-1$. The degree of the forms in $\smallsetminus[m]\Omega_{\textrm{ham}}(M)$ is given by the degree of the associated hamiltonian multi\=/vector field. Since the associated hamiltonian multi\=/vector field to $\iota_{x_1 \cdots x_n}\omega$ is $\nu_n(x_1, \ldots, x_n)$ and $\nu_n$ is of degree $-1$ we conclude that $l_n$ is also of degree $-1$ for all $n \geqslant 1$. 

Despite that all, $(\smallsetminus[m]\Omega_{\textrm{ham}}(M), \{l_n\}_{n\in\mathbb{N}})$ fails to be an \mbox{$L_{\infty}$}\=/algebra. Let $(x_i, \alpha_i)$ be a hamiltonian pair for all $i\geqslant 1$. $J(1)=d \circ d = 0$ trivially holds. We now observe that given $n\geqslant 3$ we have that $l_{n-1} \ins l_2 + l_1 \ins l_n = 0$. This is because 
\begin{eqnarray*}
(l_{n-1} \ins l_2)(\alpha_1 \otimes \cdots \otimes \alpha_n) &=& \sum_{\sigma \in \mathcal{S}h_{n-2}^2} \epsilon(\sigma) l_{n-1}(l_2(x_{\sigma(1)} \otimes x_{\sigma(2)})\otimes \alpha_{\sigma(3)} \otimes \cdots \otimes \alpha_{\sigma(n)}) \\
&=& \sum_{\sigma \in \mathcal{S}h_{n-2}^2} \epsilon(\sigma) \iota_{\nu(x_{\sigma(1)} x_{\sigma(2)}) x_{\sigma(3)} \cdots x_{\sigma(n)}} \omega \\
&=&  \iota_{\sum_{\sigma \in \mathcal{S}h_{n-2}^2} \epsilon(\sigma)\nu(x_{\sigma(1)} x_{\sigma(2)}) x_{\sigma(3)} \cdots x_{\sigma(n)}} \omega \\
&=& \iota_{\nu_n(x_1, \ldots, x_n) } \omega = d \iota_{x_1 \cdots x_n} \omega \\
&=& -(l_1 \circ l_n) (\alpha_1 \otimes \cdots \otimes \alpha_n) = -(l_1 \ins l_n) (\alpha_1 \otimes \cdots \otimes \alpha_n).
\end{eqnarray*}

We find a problem when looking at $2$\=/Jacobi:
$$J(2)(\alpha_1 \otimes \alpha_2)=(l_1 \ins l_2)(\alpha_1 \otimes \alpha_2) + (l_2 \ins l_1)(\alpha_1 \otimes \alpha_2) =
-d \iota_{x_1 x_2} \omega = -\iota_{\nu(x_1 x_2)} \omega.$$
This is in general non\=/zero, unless the last expression vanishes by degree reasons. $l_i \ins l_j$ lowers the degree by $2$. $|\alpha_1 \otimes \alpha_2| \geqslant 2$, but if that degree is precisely $2$, $|(l_i \ins l_j)(\alpha_1 \otimes \alpha_2)| = 0$, and then $(l_i \ins l_j)(\alpha_1 \otimes \alpha_2) \in \, (\smallsetminus[m]\Omega_{\textrm{ham}}(M))_0 = \{0\}$. That happens when both $\alpha_1$ and $\alpha_2$ are $(m-1)$ hamiltonian forms, that is, when the corresponding multi\=/vector fields are ordinary hamiltonian vector fields.

In order to fix this we re\=/define $l_2$ to be:

\[
 l_2(\alpha_1 \otimes \alpha_2) =
  \begin{cases} 
      \hfill \iota_{x_1 x_2} \omega    \hfill & \text{ if $|\alpha_1| = |\alpha_2| = 1$ } \\
      \hfill 0 \hfill & \text{ else. } \\
  \end{cases}
\]

In this way $2$\=/Jacobi holds.

In fact, if we keep on letting $n$ grow, at each step, for every $n\geqslant 2$; $J(n)$ does not vanish unless $l_i = 0$ outside $\left((\smallsetminus[m] \Omega_{\textrm{ham}}(M))_1\right)^{\otimes i}$ for all $i \in \, \{2, \ldots, n\}$. That makes us redefine $l_i$ in this way. We are going to prove that assertion by induction. The inductive step has already being computed ($n=2$). We assume that the statement is true for $n (\geqslant 2)$ and we are going to show it for $n+1 (\geqslant 3)$.

$$J(n+1)=\sum_{i+j=n+2} l_j \ins l_i = l_{n+1} \ins l_1 + l_{n} \ins l_2 + \left( \sum_{\begin{smallmatrix}i+j=n+2 \\ 3 \leqslant i \leqslant n \end{smallmatrix}} l_j \ins l_i\right) + l_1 \ins l_{n+1}.$$

We have previously remarked that $l_{n+1} \ins l_1 = 0$. For every $3 \leqslant i \leqslant n$ we have that $|l_i(\alpha_{1} \otimes \ldots \otimes \alpha_i)| \geqslant i-1 \geqslant 2$. Since the statement holds for $n$ we get that $l_j = 0$ outside $\left( (\smallsetminus[m] \Omega_{\textrm{ham}}(M))_1\right)^{\otimes j}$ for all $j \in \, \{2, \ldots, n\}$ so that
$\sum_{i+j=n+2}^{3 \leqslant i \leqslant n} l_j \ins l_i = 0$.

The only two summands remaining are $l_n \ins l_2 + l_1 \ins l_{n+1}$. If all the inputs are in $(\smallsetminus[m] \Omega_{\textrm{ham}}(M))_1$ that sum vanishes as we have seen before. Outside $\left((\smallsetminus[m] \Omega_{\textrm{ham}}(M))_1\right)^{\otimes n+1}$ we have that $l_n \ins l_2 = 0$, but 
$$(l_1 \ins l_{n+1})(\alpha_1 \otimes \cdots \otimes \alpha_{n+1}) = -d \iota_{x_1 \cdots x_{n+1}}\omega = -\iota_{\nu_{n+1}(x_1, \ldots, x_{n+1})} \omega.$$

This does not vanish in general, that is why we fix it by setting:
\[
 l_{n+1}(\alpha_1 \otimes \cdots \otimes \alpha_{n+1}) =
  \begin{cases} 
      \hfill \iota_{x_1 \cdots x_{n+1}} \omega    \hfill & \text{ if $|\alpha_i| = 1$ for every $1 \leqslant i \leqslant n+1$ } \\
      \hfill 0 \hfill & \text{ else. } \\
  \end{cases}
\]

As a consequence of this discussion we have the following result:

\begin{pp}\label{ROG}\index{Differential form!$L_{\infty}$\=/algebra of}
Let $(M, \omega)$ be a pre\=/$m$\=/symplectic manifold. $(\smallsetminus[m]\Omega_{\textrm{ham}}(M), \{l_n\})$ is an \Linf where $l_1 = -d$ and for every $n \geqslant 2$
\[
 l_{n}(\alpha_1 \otimes \cdots \otimes \alpha_{n}) =
  \begin{cases} 
      \hfill \iota_{x_1 \cdots x_{n}} \omega    \hfill & \text{ if $|\alpha_i| = 1$ for every $1 \leqslant i \leqslant n$ } \\
      \hfill 0 \hfill & \text{ else. } \\
  \end{cases}
\]
\end{pp}

The proof has been given in the previous paragraphs.

\begin{rk}\label{rkrog}
\begin{enumerate}
\item Since the only non\=/zero bracket defined on $(\smallsetminus[m] \Omega_{\textrm{ham}}(M))_{n\geqslant 2}$ is $l_1 = -d$, we do not need the input forms to be hamiltonian in order to have a well defined output, and hence we can extend the underlying graded vector space $L$ of the \Linf given on Proposition \ref{ROG} to be:
\[
 L_{n} =
  \begin{cases} 
      \hfill \Omega_{ham}^{m-1}(M)   \hfill & \text{ if $n=1$ } \\
      \hfill \Omega^{m-n}(M) \hfill & \text{ if $n \in \, \{2, \ldots, m\}$ } \\
			\hfill \{0\} \hfill & \text{ else. }
  \end{cases}
\]

This \Linf is the one given by Rogers in \cite{R} and the one used to define a homotopy momentum map in \cite{FRZ}. Actually, the \Linf in \cite{R} differs from this one in a sign: while $l_1$ is $-d$, the higher brackets ($n \geqslant 2$) are defined to be $-\iota_{x_1 \cdots x_n}\omega$. This occurs because of a different convention in defining what a hamiltonian pair is. This \Linf will be useful at defining what a momentum map is from a Lie algebra. In the following section we will develop an \Linf that gives rise to a different idea of momentum map from a higher \Linf acting on a manifold.

\item It is important to have in mind that this \Linf is a generalization of the \Linf of functions on a symplectic manifolds with opposite sign. We have done it this way because, as we will see in the next section, we get an \Linfm to multi\=/vector fields and not an anti\=/morphism: just like in the symplectic case. This gives a clue of how to remember the signs: the usual differential and the usual contraction for two elements ($d$ and $\hat{\iota}_{xy}\omega$) have the opposite signs to $l_1 = -d$ and $l_2(\alpha, \beta) = -\hat{\iota}_{xy}\omega = -\iota_{xy}\omega$ where $(x,\alpha)$ and $(y, \beta)$ are hamiltonian pairs.
\end{enumerate}
\end{rk}


\subsection{Higher \Linfb structures on differential forms}

In the last section we showed that the naive set of brackets on differential forms on a pre\=/multi\=/symplectic manifold does not constitute an \mbox{$L_{\infty}$}\=/algebra. We followed \cite{R} to fix the problems arising in the proof of the Jacobi identity by forcing some brackets to be zero in order to get the \Linf shown in Proposition \ref{ROG}. The new set of brackets gives an \Linf structure to differential forms, but as we observed in Remark \ref{rkrog} only top degree hamiltonian forms are considered.

That process is not the only one solving the problem, here we will expose a different idea. Instead of defining the brackets to be zero we give several copies of the same spaces, given a bi\=/grading on the underlying graded vector space.  We will define a new family of brackets which are defined in different ways depending on the bi\=/grading.

\begin{df}\label{braforms}\index{Differential form!higher $L_{\infty}$\=/algebra of}
Let $(M, \omega)$ be a pre\=/$m$\=/symplectic manifold. We define $\Xi(M) = \Xi(M,\omega)$ to be the total space of the bicomplex $\Xi_{i}^j(M)$:
$$\Xi_k(M) := \bigoplus_{i+j=k} \Xi_i^j(M) \textrm{ for every positive } k \in \mathbb{Z} \textrm{ where}$$
$$\Xi_{i=1}^j(M):= \mathsf{tr}_{m-i}(\Omega^{m-1-j}_{\textrm{ham}}(M)) \textrm{ and  } \Xi_{i \geqslant 2}^{j}(M):=\mathsf{tr}_{m-i}(\Omega^{m-i-j}(M))$$ 
(recall that $\mathsf{tr}_{m-i}$ means truncation in degrees $\{0, \ldots , m-i \}$).
We endow $\Xi(M)$ with a set of brackets given by the maps $\{l_i\}_{i\in\mathbb{N}}$. Every map is not only going to be of total degree $-1$ but more precisely of degree $(-1,0)$. Each map $l_k$ is determined by its behavior on 
$$(l_k)_{i_1 \cdots i_k}^{j_1 \cdots j_k} \colon \Xi_{i_1}^{j_1}(M) \otimes \cdots \otimes \Xi_{i_k}^{j_k}(M) \longrightarrow \Xi(M).$$
The only non\=/zero maps are
\begin{eqnarray*}
(l_1)_{i \geqslant 2}^j (\alpha) &:=& -d (\alpha) \textrm{ and } \\
(l_{k \geqslant 2})_{1 \stackrel{\scriptscriptstyle{k}}{\cdots} 1}^{j_1 \cdots j_k} (\alpha_1 \otimes \cdots \otimes \alpha_k) &:=&  \iota_{x_1 \cdots x_k} \omega.
\end{eqnarray*}
Where $(x_i, \alpha_i)$ are hamiltonian pairs for every $i$.
\end{df}

Usually we do not write the pre-multi-symplectic form $\omega$ and simply write $\Xi(M)$. The graded vector space given in Definition \ref{braforms} and $l_1 =-d$ can be pictured like this ($i$ is represented in the horizontal direction and $j$ on the vertical one):

\begin{center}
\begin{tikzpicture}[description/.style={fill=white,inner sep=2pt},descr/.style={fill=white,inner sep=2.5pt}]
\matrix (m) [matrix of math nodes, row sep=1.5em,
column sep=2.2em, text height=1.5ex, text depth=0.25ex]
{\color{BlueViolet}m-1&\Omega_{\textrm{ham}}^{0}(M)&&&&\\
\color{BlueViolet}m-2&\Omega_{\textrm{ham}}^{1}(M)&\Omega^{0}(M)&&&\\
\color{BlueViolet} \ldots &\cdots&\cdots&\cdots&&\\
\color{BlueViolet} 1& \Omega_{\textrm{ham}}^{m-2}(M)&\Omega^{m-3}(M)& \cdots&\Omega^{0}(M)&&\\
\color{BlueViolet}0& \Omega_{\textrm{ham}}^{m-1}(M)&\Omega^{m-2}(M)&\cdots&\Omega^{1}(M)&\Omega^{0}(M)\\
&\color{BlueViolet}1&\color{BlueViolet}2&\color{BlueViolet}\ldots&\color{BlueViolet}m-1&\color{BlueViolet}m\\};
\path[->,font=\scriptsize]
(m-2-3) edge node[auto] {$-d$} (m-2-2)
(m-4-5) edge node[auto] {$-d$} (m-4-4)
(m-4-4) edge node[auto] {$-d$} (m-4-3)
(m-4-3) edge node[auto] {$-d$} (m-4-2)
(m-5-6) edge node[auto] {$-d$} (m-5-5)
(m-5-5) edge node[auto] {$-d$} (m-5-4)
(m-5-4) edge node[auto] {$-d$} (m-5-3)
(m-5-3) edge node[auto] {$-d$} (m-5-2);
\end{tikzpicture}
\end{center}

It is important to remark that the lower line is closed under the brackets. This is clear because all the brackets have degree $0$ on the vertical direction. 

\begin{tm}\label{mainforms}
Let $(M, \omega)$ be a pre\=/$m$\=/symplectic manifold. Then $(\Xi(M),\{l_i\})$ is an \mbox{$L_{\infty}$}\=/algebra.
\end{tm} 

Observe that the first column is precisely $\smallsetminus [m] \Omega_{\textrm ham} (M)$, therefor the theorem could be thought of as a way of resolving the space of hamiltonian forms into an \Linf (which is minimal due to the discussion above). Before proving the Theorem we would like to observe certain things about this construction. 

\begin{rk}\label{rk2}
\begin{enumerate}

\item A priori we could have defined $\Xi_{i \neq 1}^j(M)$ to be $\Omega^{n-(i+j)}(M)$ always without truncating the complex (the same for $\Xi_{1}^j(M)$). Nevertheless we can view $\Xi_i(M) := \bigoplus_{j=0}^{m-i} \Xi_i^j(M)$ as target space of $l_{i+1}$. This justifies the choice $1 \leqslant i$. 
\item Imagine we only want to work with hamiltonian forms of degrees concentrated in $\{m-1, \ldots, m-r\}$ for some positive integer $r$. In this way we will be willing to bound the graded vector space on $j$ from above without going all the way to $n-1$. Since $l_2$ takes two elements in $\Xi_1 (M)$ and gives an element in $\Xi_1(M)$, if we start defining $\Xi_1^j(M)$ to be zero unless $j\in \, \{1, \ldots, r\}$ after $1$ application of $l_2$ we could get an element in $\Xi_1^{2r -1}(M)$, this makes us define $\Xi_1^j(M)$ to be non\=/zero for $j\in \, \{1, \ldots, 2r-1\}$. After $k$ iterations of $l_2$ we could be getting an element in $\Xi_1^{2^k(r-1) +1}(M)$ what make us define $\Xi_1^j(M)$ to be non\=/zero for $j\in \, \{1, \ldots, 2^k(r-1) +1\}$. The series $a_k:=2^k(r-1) +1$ is unbounded unless $r=1$ when it is constant $a_k = 1$. 

{\bf That tells us that the only way of bounding the complex even more is in the particular case where only top degree hamiltonian forms are required.}

\end{enumerate}
\end{rk}

In that case we define $(\Xi^0)_1(M) := \Xi^0_1(M) = \Omega_{\textrm{ham}}^{m-1}(M)$\index{Differential form!base $L_{\infty}$\=/algebra of}. Now $l_{i+1\geqslant 2}$ takes values only in $\Xi_{i}^0 (M)$ and we define $(\Xi^0)_{i}(M) := \Xi^0_{i}(M) = \Omega^{m-i}(M)$ for every $i\geqslant 2$. $\Xi^0(M)$ is no more than the base row $\Xi_{\bullet}^0(M)$:
\begin{center}
	\begin{tikzpicture}[description/.style={fill=white,inner sep=2pt},descr/.style={fill=white,inner sep=2.5pt}]
	\matrix (m) [matrix of math nodes, row sep=1.5em,
	column sep=2.2em, text height=1.5ex, text depth=0.25ex]
	{ \Omega_{\textrm{ham}}^{m-1}(M)&\Omega^{m-2}(M)&\cdots&\Omega^{1}(M)&\Omega^{0}(M)\\
		\color{BlueViolet}1&\color{BlueViolet}2&\color{BlueViolet}\ldots&\color{BlueViolet}m-1&\color{BlueViolet}m\\};
	\path[->,font=\scriptsize]
	(m-1-2) edge node[auto] {$-d$} (m-1-1)
	(m-1-5) edge node[auto] {$-d$} (m-1-4)
	(m-1-4) edge node[auto] {$-d$} (m-1-3)
	(m-1-3) edge node[auto] {$-d$} (m-1-2);
	\end{tikzpicture}
\end{center}

A direct consequence of the last statement is that the \Linf of Rogers \cite{R} is an \mbox{$L_{\infty}$}\=/subalgebra of the one in Theorem \ref{mainforms}.

\begin{cl}
	the \Linf of Rogers \cite{R} is an \mbox{$L_{\infty}$}\=/subalgebra of the one in Theorem \ref{mainforms}.
\end{cl}

\noindent{\bf Proof of Theorem \ref{mainforms}:}
The graded vector space $\Xi(M)$ is concentrated in positive degrees, in fact it is concentrated in degrees $\{1, \ldots, m\}$. The maps $\{l_k\}$ are linear (since $\widetilde{\omega}$ is linear), well defined (as we saw in the previous section, the value of $l_k$ does not depend on the paired hamiltonian multi\=/vector field chosen), symmetric (clear for $k=1$, for a larger $k$ both the contraction operator and $\omega$ are graded symmetric making $l_k$ symmetric) and of degree $-1$ (it is even of degree $(-1,0)$ taking into account the bi\=/grading).

Only the Jacobi identity remains to be checked. We are going to call the grading given by the lower\=/index, horizontal grading or simply. In this way $\alpha \in \Xi_i^j(M)$ has {\bf horizontal degree $i$}.

As in the previous theorem, we observe that $l_n \ins l_1 = 0$ for every $n$ since $0$ is a hamiltonian multi\=/vector field for an exact form ($d \circ d = 0 = \iota_{0} \omega$). This already shows that $J(1)=0$. It is also true that $l_1 \ins l_2 = 0$, just by degree reasons: the image of $l_2$ is concentrated in horizontal degree $\{1\}$, where $l_1$ is zero. This shows that $J(2)=l_1 \ins l_2 + l_2 \ins l_1 = 0+0 = 0$.

For a larger $n$, $n \geq 3$ we have as in the previous section:
$$J(n)=\sum_{p+q=n+1} l_q \ins l_p = l_{n} \ins l_1 + l_{n-1} \ins l_2 + \left( \sum_{\begin{smallmatrix}p+q=n+1 \\ 3 \leqslant p \leqslant n-1 \end{smallmatrix}} l_q \ins l_p\right) + l_1 \ins l_{n}.$$

We study that expression term by term. $l_n \ins l_1 = 0$ as we have already said. The image of $l_{3 \leqslant p \leqslant n-1}$ is concentrated in horizontal degree $p-1 \geqslant 2$, the corresponding $l_{2 \leqslant q = n+1-p \leqslant n-2}$ is zero on elements of horizontal degree $p-1$. This implies that $$ \sum_{\begin{smallmatrix}p+q=n+1 \\ 3 \leqslant p \leqslant n-1 \end{smallmatrix}} l_q \ins l_p = 0.$$

The only two remaining terms are $l_{n-1} \ins l_2$ and $l_1 \ins l_{n}$. We need to be careful with the different copies of the same space on $\Xi(M)$ to see if the result (true for $\Xi^0(M)$) also holds. We fix $\alpha_k \in \Xi_{i_k}^{j_k}(M)$ for every $k \in \, \{1, \ldots, n\}$. If any $i_k \neq 1$ then both $l_{n-1} \ins l_2$ and $l_1 \ins l_{n}$ vanish on $(\alpha_1 \otimes \cdots \otimes \alpha_n)$. If every $i_k = 1$ let $(x_k, \alpha_k)$ be hamiltonian pairs for every $k$. Then 
$$(l_1 \ins l_{n})(\alpha_1 \otimes \cdots \otimes \alpha_n)=-d\iota_{x_1 \cdots x_n}\omega \in \, \Xi_{n-2}^{(\sum_{k=1}^{n} j_k)}(M) \textrm{ and}$$
$$(l_{n-1} \ins l_2)(\alpha_1 \otimes \cdots \otimes \alpha_n)=\iota_{\nu_n(x_1, \ldots, x_n)}\omega = d\iota_{x_1 \cdots x_n}\omega \in \, \Xi_{n-2}^{(\sum_{k=1}^{n} j_k)}(M).$$
We have concluded the  proof since
$$J(n)=\sum_{p+q=n+1} l_q \ins l_p = 0 + l_{1} \ins l_n + 0 - l_1 \ins l_{n} = 0.$$
\qed


\subsection{\Linfb structures on hamiltonian pairs}

Theorem \ref{mainforms} together with Remark \ref{rk2} give two \Linf structures on differential forms on a multi\=/symplectic manifold based on the behavior of hamiltonian forms. Those \Linfs and the one on hamiltonian multi\=/vector fields are closely related due to Proposition \ref{maincl}.

\begin{lm}\label{pair}
Let $(M, \omega)$ be a multi\=/symplectic manifold. For every positive integer $n$ and every family of hamiltonian pairs $(x_k, \alpha_k)$ where $\alpha_k \in \, \Xi_1(M)$ for every $k \in \, \{1, \ldots, n\}$ we have that $(\nu_n(x_1, \ldots, x_n), l_n(\alpha_1 \otimes \cdots \otimes \alpha_n))$ is a hamiltonian pair. 
\end{lm}

Lemma \ref{pair} suggests an \Linf structure on hamiltonian pairs. 

\begin{df}\label{brapairs}\index{Hamiltonian!pair!$L_{\infty}$\=/algebra of}
Let $(M, \omega)$ be a pre\=/$m$\=/symplectic manifold. We define $\widetilde{\Xi}(M)$ to be the total space of the bicomplex $\widetilde{\Xi}_{i}^j(M)$:
$$\widetilde{\Xi}_k(M) := \bigoplus_{i+j=k} \widetilde{\Xi}_i^j(M) \textrm{ for every positive } k \in \mathbb{Z} \textrm{ where}$$
$$\widetilde{\Xi}_{i=1}^j(M):= \mathsf{tr}_{m-i}(\widetilde{\mathfrak{X}}^{j+1}(M)) \textrm{ and  }\widetilde{\Xi}_{i \geqslant 2}^{j}(M):=\mathsf{tr}_{m-i}(\Omega^{m-i-j}(M)).$$ 
We define a family of brackets $\{\tilde{l}_k\}$ on $\widetilde{\Xi}(M)$ of degree $(-1,0)$ where the only non\=/zero maps are
$$ (\tilde{l}_1)_{i \geqslant 2}^j  \textrm{ and } (\tilde{l}_{k \geqslant 2})_{1 \stackrel{\scriptscriptstyle{k}}{\cdots} 1}^{j_1 \cdots j_k}.$$
Which are given by 
\[
 (\tilde{l}_1)_i^j(\alpha) =
  \begin{cases} 
      \hfill (0, -d \alpha)  \hfill & \text{ if $i=2$ } \\
      \hfill -d \alpha \hfill & \text{ if $i \in \, \{3, \ldots, m\}$} \\
			\hfill 0 \hfill & \text{ else. }
  \end{cases}
\]
\[
 (\tilde{l}_{k\geqslant 2})((x_1, \alpha_1) \otimes \cdots \otimes (x_k, \alpha_k)) =
  \begin{cases} 
      \hfill (\nu(x_1 x_2), \iota_{x_1 x_2} \omega)  \hfill & \text{ if $k=2$ } \\
      \hfill \iota_{x_1 \cdots x_k} \omega & \text{ if $k \geqslant 3$.}
  \end{cases}
\]
for every $((x_1, \alpha_1) \otimes \cdots \otimes (x_k, \alpha_k)) \in \, \widetilde{\Xi}_1^{j_1}(M) \otimes \cdots \otimes \widetilde{\Xi}_1^{j_k}(M)$.
\end{df}

We can picture the graded vector space, together with $\tilde{l}_1$, in the following diagram. See the resemblance to the diagram describing the \Linf of hamiltonian forms studied before.

\begin{center}
\begin{tikzpicture}[description/.style={fill=white,inner sep=2pt},descr/.style={fill=white,inner sep=2.5pt}]
\matrix (m) [matrix of math nodes, row sep=1.5em,
column sep=2.2em, text height=1.5ex, text depth=0.25ex]
{\color{BlueViolet}m-1& \widetilde{\mathfrak{X}}^{m}(M)&&&&\\
\color{BlueViolet}m-2&\widetilde{\mathfrak{X}}^{m-1}(M)&\Omega^{0}(M)&&&\\
\color{BlueViolet} \ldots &\cdots&\cdots&\cdots&&\\
\color{BlueViolet} 1& \widetilde{\mathfrak{X}}^{2}(M)&\Omega^{m-3}(M)& \cdots&\Omega^{0}(M)&&\\
\color{BlueViolet}0& \widetilde{\mathfrak{X}}^{1}(M)&\Omega^{m-2}(M)&\cdots&\Omega^{1}(M)&\Omega^{0}(M)\\
&\color{BlueViolet}1&\color{BlueViolet}2&\color{BlueViolet}\ldots&\color{BlueViolet}m-1&\color{BlueViolet}m\\};
\path[->,font=\scriptsize]
(m-2-3) edge node[auto] {$\tilde{l}_1$} (m-2-2)
(m-4-5) edge node[auto] {$\tilde{l}_1$} (m-4-4)
(m-4-4) edge node[auto] {$\tilde{l}_1$} (m-4-3)
(m-4-3) edge node[auto] {$\tilde{l}_1$} (m-4-2)
(m-5-6) edge node[auto] {$\tilde{l}_1$} (m-5-5)
(m-5-5) edge node[auto] {$\tilde{l}_1$} (m-5-4)
(m-5-4) edge node[auto] {$\tilde{l}_1$} (m-5-3)
(m-5-3) edge node[auto] {$\tilde{l}_1$} (m-5-2);
\end{tikzpicture}
\end{center}

\begin{cl}\label{cpairs}
Let $(M, \omega)$ be a pre\=/$m$\=/symplectic manifold. Then $(\widetilde{\Xi}(M), \{\tilde{l}_k\})$ is an \mbox{$L_{\infty}$}\=/algebra.
\end{cl}

\dem
Every map is a product of two well defined symmetric linear maps of degree $-1$ by definition. The image of the maps lands on hamiltonian pairs by Lemma \ref{pair}. The Jacobi identity on the form part is satisfied due to Theorem \ref{mainforms}.

$J(1)= 0$ and $J(2)=0$ since $\tilde{l}_1$ is always zero in the multi\=/vector fields part. As usual we write for every $n\geqslant 3$
$$J(n)=\sum_{p+q=n+1} \tilde{l}_q \ins \tilde{l}_p = \tilde{l}_{n} \ins \tilde{l}_1 + \tilde{l}_{n-1} \ins \tilde{l}_2 + \left( \sum_{\begin{smallmatrix}p+q=n+1 \\ 3 \leqslant p \leqslant n-1 \end{smallmatrix}} \tilde{l}_q \ins \tilde{l}_p\right) + \tilde{l}_1 \ins \tilde{l}_{n}.$$

The corresponding equation in the multi\=/vector fields part behaves as follows. $\sum_{\begin{smallmatrix}p+q=n+1 \\ 3 \leqslant p \leqslant n-1 \end{smallmatrix}} \tilde{l}_q \ins \tilde{l}_p= 0$ by horizontal degree reasons. $\tilde{l}_{n} \ins \tilde{l}_1=0=\tilde{l}_1 \ins \tilde{l}_{n}$ because $\tilde{l}_1$ is always zero in the multi\=/vector field part. The only remaining term is $\tilde{l}_{n-1} \ins \tilde{l}_2$, which on the multi\=/vector field part is zero unless $n=3$. In that case $(\nu_{2}) \ins (\nu_2)=0$ because of the Jacobi identity for the Schouten bracket.
\qed\\

Just alike in the previous section, we want to point out the construction given in that Corollary also generalizes some construction used in \cite{FRZ}. That is done by trying to further bound the underlying vector space as we did for hamiltonian forms. The resulting construction is the following:

\begin{cl}
This \Linf called the Poisson Lie\=/m\=/algebra in \cite[Theorem 4.7]{FRZ} is an \mbox{$L_{\infty}$}\=/subalgebra of the one from Corollary \ref{cpairs}.	
\end{cl}

\begin{center}
	\begin{tikzpicture}[description/.style={fill=white,inner sep=2pt},descr/.style={fill=white,inner sep=2.5pt}]
	\matrix (m) [matrix of math nodes, row sep=1.5em,
	column sep=2.2em, text height=1.5ex, text depth=0.25ex]
	{ \widetilde{\mathfrak{X}}^{1}(M)&\Omega^{m-2}(M)&\cdots&\Omega^{1}(M)&\Omega^{0}(M)\\
		\color{BlueViolet}1&\color{BlueViolet}2&\color{BlueViolet}\ldots&\color{BlueViolet}m-1&\color{BlueViolet}m\\};
	\path[->,font=\scriptsize]
	(m-1-2) edge node[auto] {$-d$} (m-1-1)
	(m-1-5) edge node[auto] {$-d$} (m-1-4)
	(m-1-4) edge node[auto] {$-d$} (m-1-3)
	(m-1-3) edge node[auto] {$-d$} (m-1-2);
	\end{tikzpicture}
\end{center}

We denote it by $\widetilde{\Xi}^0 (M)$ and the inclusion will be denoted by $\widetilde{\textrm{i}}^0 \colon \widetilde{\Xi}^0 (M) \rightarrow \widetilde{\Xi} (M)$ similarly to the inclusion for the \Linfs of Hamiltonian forms: $\textrm{i}^0 \colon \Xi^0 (M) \rightarrow \Xi (M)$.

%
%
%


\subsection{Relating the \Linfsb of multi\=/vector fields, forms and pairs}


The \Linf structures on hamiltonian multi\=/vector fields, those on hamiltonian pairs and those on hamiltonian forms are very close. Given any hamiltonian pair, we can {\bf project onto the multi\=/vector field} part (in horizontal degree $1$ and $0$ elsewhere).
 
$$\pi \colon (\smallsetminus(\widetilde{\Xi} (M)), \smallsetminus(\tilde{l}_1)) \longrightarrow (\smallsetminus(\ham), \smallsetminus(\nu_1)=0) \textrm{ and }$$
\noindent{ {\bf project  onto the form part} at horizontal degree $1$ and the identity elsewhere.
$$\phi \colon (\smallsetminus(\widetilde{\Xi} (M)), \smallsetminus(\tilde{l}_1)) \longrightarrow (\smallsetminus({\Xi} (M)), \smallsetminus(l_1)).$$	
Both are cochain maps by construction. For example $\pi \tilde{l}_1(\smallsetminus(x, \alpha)) = \pi(\smallsetminus(0,d\alpha)) = 0 = (0 \circ \pi)(\smallsetminus(x,\alpha))$. The same occurs if we replace $\pi$ by $\pi \circ \widetilde{\textrm{i}}^0$, or $\phi$ by $\phi \circ \widetilde{\textrm{i}}^0$. The later one actually factors through the \Linf of Rogers: 
$$\phi^{\prime} \colon (\smallsetminus(\widetilde{\Xi}^0 (M)), \smallsetminus(\tilde{l}_1)) \longrightarrow (\smallsetminus({\Xi}^0 (M)), \smallsetminus(l_1)).$$

These cochain maps lift to strict \mbox{$L_{\infty}$}\=/morphisms. This fact is exposed in the following lemma, which is a generalization (for $r=1$) of \cite[Proposition 4.8]{FRZ}.

\begin{lm}
Let $(M, \omega)$ be a pre\=/multi\=/symplectic manifold. The cochain maps 
$$\pi \colon (\smallsetminus(\widetilde{\Xi}) (M), \smallsetminus(\tilde{l}_1)) \longrightarrow (\smallsetminus(\ham), \smallsetminus(\nu_1)=0),$$
$$\pi^{\prime} = \pi \circ \widetilde{\textrm{i}}^0 \colon (\smallsetminus(\widetilde{\Xi}^0) (M), \smallsetminus(\tilde{l}_1)) \longrightarrow (\smallsetminus(\ham), \smallsetminus(\nu_1)=0),$$
$$\phi \colon (\smallsetminus(\widetilde{\Xi} (M)), \smallsetminus(\tilde{l}_1)) \longrightarrow (\smallsetminus({\Xi} (M)), \smallsetminus(l_1)) \textrm{ and}$$
$$\phi^{\prime} \colon (\smallsetminus(\widetilde{\Xi}^0 (M)), \smallsetminus(\tilde{l}_1)) \longrightarrow (\smallsetminus({\Xi}^0 (M)), \smallsetminus(l_1)).$$
lift to strict \mbox{$L_{\infty}$}\=/morphisms given by
$$\pi_1=\smallsetminus(\pi) \colon (\widetilde{\Xi} (M), \{\tilde{l}_n\}) \longrightarrow (\ham, \nu_2),$$
$$\pi_1^{\prime}=\smallsetminus(\pi^{\prime}) \colon (\widetilde{\Xi}^0 (M), \{\tilde{l}_n\}) \longrightarrow (\ham, \nu_2),$$
$$\phi_1=\smallsetminus(\phi) \colon (\widetilde{\Xi} (M), \{\tilde{l}_n\}) \longrightarrow ({\Xi} (M), \{\l_n\}) \textrm{ and}$$
$$\phi_1^{\prime}=\smallsetminus(\phi^{\prime}) \colon (\widetilde{\Xi}^0 (M), \{\tilde{l}_n\}) \longrightarrow ({\Xi}^0 (M), \{\l_n\}).$$
Even more, if $(M, \omega)$ is symplectic, then $\phi_1^\prime$ is a strict \mbox{$L_{\infty}$}\=/quasi\=/isomorphism.
\end{lm}

\dem
It is enough to show it for $\pi$ and $\phi$ since $\widetilde{\textrm{i}}^0$ and $\textrm{i}^0$ clearly lift to strict \mbox{$L_{\infty}$}\=/morphisms.
We begin showing the strict \Linfm at hamiltonian multi\=/vector fields. Let $\{r_n\}$ be the family of brackets in $(\ham, \nu_2)$ ($r_2 = \nu_2$ and $r_{n\neq 2}=0$). We have to check for every positive $n$ that $\pi_1 \circ \tilde{l}_n = r_n \circ (\pi_1 \otimes \stackrel{\scriptscriptstyle{n}}{\cdots} \otimes \pi_1)$. If $n=1$ the statement is just equivalent to the fact that $\smallsetminus(\pi)$ is a cochain map. For $n \geqslant 2$, $\tilde{l}_n$ is non\=/zero only when all the inputs are in $(\widetilde{\Xi})_1^{\bullet} (M)$. Given $n$ elements $\{(x_i, \alpha_i) \in \, \widetilde{\Xi}_1^{j_i}(M)\}_{i=1}^n$ we have that $\tilde{l}_n ((x_1, \alpha_1) \otimes \cdots \otimes (x_n, \alpha_n)) \in \, \widetilde{\Xi}_{n-1}^{\sum_{i=1}^{n}j_i }(M)$, $\pi_1$ on such term is zero unless $n = 2$: precisely the case in which $r_n(\pi_1 \otimes \stackrel{\scriptscriptstyle{n}}{\cdots} \otimes \pi_1)$ is non\=/zero. In that case
$$ (\pi_1 \circ \tilde{l}_2)((x_1, \alpha_1) \otimes (x_2, \alpha_2)) = \nu_2(x_1 \wedge x_2) = (\nu_2 \circ (\pi_1 \otimes \pi_1))((x_1, \alpha_1) \otimes (x_2, \alpha_2)).$$
The part of the statement regarding hamiltonian forms holds trivially.
\begin{eqnarray*}
(\phi_1 \circ \tilde{l}_n)((x_1, \alpha_1) \otimes \cdots \otimes (x_n, \alpha_n)) = l_n(x_1, \ldots, x_n) && \\
= (l_n \circ (\phi_1 \otimes \stackrel{\scriptscriptstyle{n}}{\cdots} \otimes \phi_1))((x_1, \alpha_1) \otimes \cdots \otimes (x_n, \alpha_n)). &&
\end{eqnarray*}
Even more, $\phi_1^\prime$ is given by the identity every\=/where except on degree $1$ when it is equal to $(\widetilde{\omega}, \textrm{id}) \colon \widetilde{\mathfrak{X}}(M) \longrightarrow \Omega^{m-1}(M)$. If $(M,\omega)$ is multi\=/symplectic, then that map is an isomorphism, turning $\phi_1^{\prime}$ into a quasi\=/isomorphism. 
\qed

We could have defined $\pi$ towards multi\=/vector fields and then just observe that the image of $\pi$ is indeed always hamiltonian.

\subsection{Synchronized morphisms}

At this moment, we are ready to introduce a kind of map that will be the key to define higher momentum maps. We will be talking about pairs of maps $({\sf a}, \tilde{f})$ from an \Linf $L$ to $\widetilde{\Xi}(M)$ and $\widetilde{\Xi}(M)$ making the following diagram in the category of \Linfs commute:
 
\begin{center}
	\begin{tikzpicture}[description/.style={fill=white,inner sep=2pt}]
	\matrix (m) [matrix of math nodes, row sep=1.5em,
	column sep=3.5em, text height=1.5ex, text depth=0.25ex]
	{ & \widetilde{\Xi}(M) \\
		L & \mathfrak{X}^{\bullet}(M) \\};
	\path[->,font=\scriptsize]
	(m-2-1) edge node[auto] {${\sf a}$} (m-2-2)
	(m-1-2) edge node[auto] {$\pi_1$} (m-2-2)
	(m-2-1) edge node[auto] {$\tilde{f}$} (m-1-2);
	\end{tikzpicture}
\end{center}

For obvious reasons we will say that $\tilde{f}$ {\bf lifts} ${\sf a}$.

Since the target space of $\tilde{f}$ has several copies of the same underlying vector spaces (the spaces of differential forms of a given degree), it will become very important to be able to control to which copy are things being mapped to. Let us be explicit, when talking about maps $\{\tilde{f}_n \colon \wedge^n L \rightarrow \Xi(M) \textrm{ or } \widetilde{\Xi}(M)\}$ we would like to use the $i$-th column to be the target of $\tilde{f}_i (\wedge^i L)$. In this way, the image of $L$ via $\tilde{f}_1$ is fully contained in the first column: that of Hamiltonian forms or pairs. Such maps are going to be called synchronized with the horizontal grading, or simply {\bf synchronized}:

$$\tilde{f} \textrm{ is synchronized iff } \tilde{f}_i (\wedge^i L) \subset \widetilde{\Xi}_i (M).$$

The definition can be given in general for an \Linf with an extra grading:

\begin{df}[Synchronized\index{Linfinity morphism @ $L_{\infty}$\=/morphism!synchronized} \Linfm]
Let $(\{V_n\}, \{v_n\})$ be an \mbox{$L_{\infty}$}-algebra whose underlying total space $V^{\oplus}$ is given another grading $(\varphi, U)$. An \Linf morphism $f \colon L \longrightarrow V$ is said to be synchronized with the grading $(\varphi, U)$ if $f_n((L)^{\otimes n}) \subset U_n$ for every positive $n$.
\end{df}

The following Lemma will be key in understanding the connection between the different notions of moment maps that will be given in the next sections:

\begin{lm}\label{lmlift}
Let $(L, \{l_n\})$ be an \Linf concentrated in degrees $\{1, \ldots, r\}$. Let ${\sf a} \colon (L,\{l_n\}) \longrightarrow (\mathfrak{X}^{\bullet}(M), \nu_2)$ be an \Linfm. There is a one\=/to\=/one correspondence between
\begin{enumerate}
	\item synchronized \mbox{$L_{\infty}$}\=/morphisms $\tilde{f}$ from $L$ to $\widetilde{\Xi}(M)$ lifting ${\sf a}$, and
	\item synchronized \mbox{$L_{\infty}$}\=/morphisms $f$ from $L$ to $\Xi(M)$ with the property that $\widetilde{\omega}\circ {\sf a}_1 = d \circ f_1$. 
\end{enumerate}

If $r=1$ we have the same one\=/to\=/one correspondence between 
\begin{enumerate}
	\item \mbox{$L_{\infty}$}\=/morphisms $\tilde{f}$ from $L$ to $\widetilde{\Xi}^0(M)$ lifting ${\sf a}$, and
	\item \mbox{$L_{\infty}$}\=/morphisms $f$ from $L$ to $\Xi^0(M)$ with the property that $\widetilde{\omega}\circ {\sf a}_1 = d \circ f_1$. 
\end{enumerate}
\end{lm}

\dem
Let $\tilde{f}$ be a synchronized lift of ${\sf a}$. We take $f := \phi_1 \circ \tilde{f}$ which is an \Linfm from $L$ to $\Xi(M)$. Moreover, it is synchronized since $\phi_1((\widetilde{\Xi})_n^{\bullet}(M)) \subset (\Xi)_n^{\bullet}(M)$ for every $n$ by definition of $\phi_1$. $\tilde{f}_1$ takes values in $(\widetilde{\Xi})_1^{\bullet}(M)$ and the lifting condition gives that $\tilde{f}_1 = ({\sf a}_1, f_1)$ and since it is hamiltonian pair $\widetilde{\omega}\circ {\sf a}_1 =\iota_{{\sf a}_1} \omega = d \circ f_1$.

Conversely, if $f$ is a synchronized \Linfm from $L$ to $\Xi(M)$ such that $\widetilde{\omega}\circ {\sf a}_1 = d \circ f_1$. We fix $\tilde{f}_{j\geqslant 2} = f_{j}$ and $\tilde{f}_1 = ({\sf a}_1, f_1)$. Clearly $\tilde{f}$ takes values in $\widetilde{\Xi}(M)$ because $\widetilde{\omega}\circ {\sf a}_1 = d \circ f_1$. It is an \Linfm because both the form and the multi\=/vector field part are $L_{\infty}$\=/morphisms and because $\pi_1$ is strict. It is a lift and it is synchronized by construction.
At the case $r=1$ the only thing to remark is that \mbox{$L_{\infty}$}\=/morphisms and synchronized \mbox{$L_{\infty}$}\=/morphisms are the same. The result follows from the previous lines by the same construction, noticing again that underlying vector spaces of $\Xi^0(M)$ and $\widetilde{\Xi}^0(M)$ are the base arrows of $\Xi(M)$ and $\widetilde{\Xi}(M)$ respectively.
\qed

\begin{cl}\label{strict}
	Let $(L, \{l_n\})$ be an \Linf concentrated in degrees $\{1, \ldots, r\}$. Let ${\sf a} \colon (L,\{l_n\}) \longrightarrow (\mathfrak{X}^{\bullet}(M), \nu_2)$ be an \Linfm such that there exists a synchronized \mbox{$L_{\infty}$}\=/morphisms $f$ from $L$ to $\Xi(M)$ with the property that $\widetilde{\omega}\circ {\sf a}_1 = d \circ f_1$. 
	Then ${\sf a}$ is strict.
\end{cl}

\dem
From Lemma \ref{lmlift} we know that there exists a lift $\widetilde{f}$ of ${\sf a}$. For every $n \geqslant 2$, ${\sf a}_n = \pi_1 \circ \widetilde{f}_n$. But $\widetilde{f}_n (\wedge^n L) \subset \widetilde{\Xi}_n (M)$ and $\pi_i(\widetilde{\Xi}_n (M)) = 0$ for $n \geqslant 2$.
\qed

As a consequence of this Corollary, $(l_n){\sf a}(l_1)$
$$ {\sf a}_1 \circ l_2 = \nu_2 ({\sf a}_1, {\sf a}_1) \textrm{ and }.$$
$$ {\sf a}_1 \circ l_n = 0 \textrm{ for all } n \geqslant 3.$$
In particular if $L$ is a Lie algebra, ${\sf a} = {\sf a}_1$ is a honest Lie algebra morphism (whenever shifted appropriately).


\color{BlueViolet}
\section{Transfer: cochain homotopies and \mbox{\boldmath ${L_{\infty}}$}\=/morphisms}
\color{black}

The second big result in this paper is about the equivalence between different notions of higher momentum maps in multi\=/symplectic geometry. As hinted in the previous section, this will be done by comparing \Linf lifts and maps with certain cohomological properties. As a matter of fact, the result can be extrapolated outside if the world of symplectic geometry and be stated in a homological algebra fashion. That is precisely what will be done in this chapter. The main idea is to give a general tool to relate \Linf morphisms and cochain null\=/homotopies.

We begin by fixing some notation. We will be following \cite{W} and \cite{DS}. We consider cochain complexes\index{Cochain!complex} of vector spaces $C=\{C^n\}_{n\in \mathbb{Z}}$ with linear maps as co\=/differentials $d^n \colon C^{n} \longrightarrow C^{n+1}$. An element in the image of $d$ is called a co\=/boundary and one in the kernel of $d$ is called a co\=/cycle. Cochain null\=/homotopies can be defined in terms of the cone object:

\begin{df}[Cone]\index{Cochain!cone}
Let $(C, d)$ be a cochain complex. The cone cochain complex $(\textrm{cone}(C), d_{\textrm{cone}})$ is given by $\textrm{cone}(C)^n := C^{n-1} \oplus C^{n}$ for every $n \in \mathbb{Z}$ and $d_{\textrm{cone}}^n(c_1, c_2) := (c_2 -d^{n-1} c_1, d^n c_2)$ for every $(c_1, c_2) \in \, \textrm{cone}(C)^{n}$.
\end{df}

\begin{df}[Cochain null\=/homotopy]\index{Cochain!null\=/homotopy}
	A cochain null\=/homotopy between the cochain complexes $(C, d_{C})$ and $(D, d_{D})$ is a cochain map $(h,f)$ from $(C, d_{C})$ to $(\textrm{cone}(D), d_{\textrm{cone}(D)})$. $f$ is said to be null\=/homotopic by $h$, something which is denoted by $f \raisebox{-0.2cm}{ $\stackrel{\simeq}{\scriptscriptstyle{h}}$ } 0$. 
\end{df}

A cochain homotopy defining that homotopy $\varphi=(h,0,f)$ is the same as a cochain map 
$$\varphi^{\prime}=(h,f)\colon (C, d_{C}) \longrightarrow (\textrm{cone}(D), d_{\textrm{cone}(D)})$$ because the cochain maps conditions are alike:
$$(f -d_{D} h, d_D f) = d_{\textrm{cyl}(D)} \circ \varphi^{\prime} =  \varphi \circ d_{C} = (h d_C, 0, f d_C) \textrm{ and}$$
$$(f -d_{D} h, 0, d_D f) = d_{\textrm{cyl}(D)} \circ \varphi =  \varphi \circ d_{C} = (h d_C, 0, f d_C).$$ 

A cochain null\=/homotopy is then the same as a pair $(h,f)$ where $f$ is a cochain map $f \colon (C, d_{C}) \longrightarrow (D, d_{D})$ and $h \in \, \underline{\textrm{Hom}}_{{\sf grVec}}(C, D)_{-1}$ such that
$$ f = h \circ d_C + d_D \circ h.$$


\subsection{Transfer theorem for cochain homotopies from an \Linfb}\label{52}

In this section we are going to try to reproduce in a more general setting what we have done in the previous chapter concerning differential forms: we were able to transfer the \Linf structure from $\mathfrak{X}_{\mathrm{ham}}(M)$ to $\Omega^{\bullet}(M)$.

\begin{rk}\label{tfrk} Given a positive integer $m$ are going to use the following notations:
\begin{enumerate}
\item $(D, d)$ is a cochain complex concentrated in degrees $\{0, \ldots, m\}$.
\item $(\mathfrak{g}, \nu)$ is going to represent a graded Lie\=/algebra and $(C(\mathfrak{g}), d_{\textrm{CE}})$ the associated cochain complex given by Proposition \ref{cochain}. The graded vector space is $C(\mathfrak{g})^n = (S^{\bullet}\mathfrak{g})_{m+1-n}$ for $n \in \, \{0, \ldots, m\}$ and $0$ otherwise.
\end{enumerate}
\end{rk}

It is possible to think of the truncation to degrees $\{0, \ldots, m\}$ of a larger cochain complex. To be precise, $(C(\mathfrak{g}), d_{\textrm{CE}})$ is the cochain complex from Proposition \ref{cochain} but shifted by $m+1$ and truncated in degrees $\{0, \ldots, m\}$.

\begin{dfpp}[Hamiltonian coboundaries]
	For any cochain map $f \colon (C(\mathfrak{g}), d_{\textrm{CE}} ) \longrightarrow (D,d)$ we define the sub\=/cochain complex of hamiltonian elements $D^{f} \subset D$ given by 
	$$D^{f} := \{ \alpha \in D \colon d \alpha \in \mathrm{im}(f)\}.$$
\end{dfpp}

Observe that $D^{f}$ is a sub\=/cochain complex since for every $\alpha \in \, D^{f}$ we have that $d (d \alpha) = 0 \in \mathrm{im}(f)$. Observe also that $\mathrm{im}(f) \subset D^{f}$ since $d(f(x)) = f(c(x)) \in \mathrm{im}(f)$. Similarly we can define $f$-hamiltonian pairs:

\begin{df}[Hamiltonian pairs]
	Let $f \colon C(\mathfrak{g}) \longrightarrow D$ be in the hypothesis of Theorem \ref{maintf}. The graded vector space of pairs $\widetilde{D^{f}}$ is the graded vector space concentrated in degrees $\{0, \ldots , m-1\}$ given degree\=/wise by:
	$$\left( \widetilde{D^{f}} \right)^i := \{(x,\alpha) \in \, C(\mathfrak{g})^i \times \left(D^f \right)^{m-i} \colon d \alpha = f (x) \}.$$
\end{df}

The other key ingredient in the last chapter was the algebra multiplication $\otimes$ on $S^{\bullet}\mathfrak{g}$. This gives rise to a $0$\=/degree map $\mu^{k} \colon (S^{\bullet}\mathfrak{g})^{\otimes k} \longrightarrow S^{\bullet}\mathfrak{g}$ as we have seen in the first chapter. Its decalage, $\dec(\mu^{k})$ is of degree $(1-k)(m+1)$.

The following theorem is the generalization of Theorem \ref{mainforms} to the more general setting given by Remark \ref{tfrk}.

\begin{tm}[Transfer Theorem]\label{maintf}\index{Transfer!theorem}\label{mainTransfer}
Let $f \colon (C(\mathfrak{g}), d_{\textrm{CE}} ) \longrightarrow (D,d)$ satisfying the following properties:
\begin{enumerate}
\item For every $x$ and $y$ in $\mathfrak{g}$ such that $f(x) = f(y)$ it follows that $f(x,-) = f(y,-)$.
\item There exists a morphism $p \in \underline{\textrm{Hom}}_{{\sf grVec}}(D^{f},\mathfrak{g}\subset C(\mathfrak{g}))_1$ such that $f \circ p = d$.
\end{enumerate}
We define $Q$ to be the total graded vector space of the bi\=/graded space $Q_i^j$:
$$Q_k := \bigoplus_{i+j=k} Q_i^j \textrm{ for every positive } k \in \, \mathbb{Z} \textrm{ where}$$

$$Q_{i=1}^j := \mathsf{tr}_{m-1} \left( D^{f} \phantom{.}^{m-(1+j)} \right) \textrm{ and } Q_{i\geqslant 2}^j := \mathsf{tr}_{m-i} \left( D^{m-(i+j)} \right).$$
We endow $Q$ with a family of brackets of bi\=/degree $(-1, 0)$ given by the maps $\{q_n \colon Q^{\otimes n} \longrightarrow Q\}$. The only non\=/zero maps are:
\begin{eqnarray*}
\left(q_1 \right)_{i \geqslant 2}^j := -d \textrm{ and}\\
\left( q_{k \geqslant 2} \right)_{1 \stackrel{\scriptscriptstyle{k}}{\cdots} 1}^{j_1 \cdots j_k} := f \circ \dec(\mu^k) \circ p^{\otimes k}.
\end{eqnarray*}
Then $(Q, \{q_n\})$ is an \Linf and we modify $p$ to get a strict \Linfm from $Q$ to $\mathfrak{g}$ by extending it to be zero in $Q_{i\geqslant 2}^{\bullet}$.
\end{tm}

The double chain complex $Q$ can be represented in the following diagram:

\begin{center}
	\begin{tikzpicture}[description/.style={fill=white,inner sep=2pt},descr/.style={fill=white,inner sep=2.5pt}]
	\matrix (m) [matrix of math nodes, row sep=1.5em,
	column sep=2.2em, text height=1.5ex, text depth=0.25ex]
	{\color{BlueViolet}m-1& \left(D^{f}\right)^{0}&&&&\\
		\color{BlueViolet}m-2&\left(D^{f}\right)^{1}&D^{0}&&&\\
		\color{BlueViolet} \ldots &\cdots&\cdots&\cdots&&\\
		\color{BlueViolet} 1& \left(D^{f}\right)^{m-2}&D^{m-3}& \cdots&D^{0}&&\\
		\color{BlueViolet}0& \left(D^{f}\right)^{m-1}&D^{m-2}&\cdots&D^{1}&D^{0}\\
		&\color{BlueViolet}1&\color{BlueViolet}2&\color{BlueViolet}\ldots&\color{BlueViolet}m-1&\color{BlueViolet}m\\};
	\path[->,font=\scriptsize]
	(m-2-3) edge node[auto] {$-d$} (m-2-2)
	(m-4-5) edge node[auto] {$-d$} (m-4-4)
	(m-4-4) edge node[auto] {$-d$} (m-4-3)
	(m-4-3) edge node[auto] {$-d$} (m-4-2)
	(m-5-6) edge node[auto] {$-d$} (m-5-5)
	(m-5-5) edge node[auto] {$-d$} (m-5-4)
	(m-5-4) edge node[auto] {$-d$} (m-5-3)
	(m-5-3) edge node[auto] {$-d$} (m-5-2);
	\end{tikzpicture}
\end{center}

The proof is included in the next subsection. The similar result for pairs follows from the Theorem:

\begin{cl}[Transfer Theorem for pairs]\label{tfp}\index{Transfer!theorem!for pairs}
Let $f \colon (C(\mathfrak{g}), d_{\textrm{CE}} ) \longrightarrow (D,d)$ be a cochain map in the hypothesis of Theorem \ref{maintf}. We define $\widetilde{Q}$ to be the total graded vector space of the bi\=/graded space $\widetilde{Q}_i^j$:
$$\widetilde{Q}_k := \bigoplus_{i+j=k} \widetilde{Q}_i^j \textrm{ for every positive } k \in \, \mathbb{Z} \textrm{ where}$$
$$\widetilde{Q}_{i=1}^j := \widetilde{D^f}_{j+1} \textrm{ and } \widetilde{Q}_{i\geqslant 2}^j := Q_{i}^j.$$
We endow $\widetilde{Q}$ with a family of brackets of bi\=/degree $(-1, 0)$ given by the maps $\{\widetilde{q}_n \colon \widetilde{Q}^{\otimes n} \longrightarrow \widetilde{Q}\}$. The only non\=/zero maps are:
\[
 (\tilde{q}_1)_i^j(\alpha) =
  \begin{cases} 
      \hfill (0, -d \alpha)  \hfill & \text{ if $i=2$ } \\
      \hfill -d \alpha \hfill & \text{ if $i \in \, \{3, \ldots, m\}$} \\
			\hfill 0 \hfill & \text{ else. }
  \end{cases}
\]
\[
 (\tilde{q}_{k\geqslant 2})((x_1, \alpha_1) \otimes \cdots \otimes (x_k, \alpha_k)) =
  \begin{cases} 
      \hfill (\nu(x_1 x_2), q_2(\alpha_1,\alpha_2))  \hfill & \text{ if $k=2$ } \\
      \hfill q_k(\alpha_1, \ldots, \alpha_k) & \text{ if $k \geqslant 3$.}
  \end{cases}
\]
Then $(\widetilde{Q}, \{\widetilde{q}_n\})$ is an \Linf and we can modify $p$ to get a strict \Linfm from $\widetilde{Q}$ to $\mathfrak{g}$ by extending it to be zero in $\widetilde{Q}_{i\geqslant 2}^{\bullet}$.
\end{cl}

%

There is a strict \Linfm $\phi \colon \widetilde{Q} \longrightarrow Q$ given by the projection on the $D$ factor. The bottom row of both \Linfs are \mbox{$L_{\infty}$}\=/subalgebras.

\begin{pp}\label{cltf}
	In the hypothesis of Theorem \ref{maintf} $\left(\left(Q^0\right)_{\bullet} := Q_{\bullet}^0, \, \{{q_n}_{|Q^0}\}\right)$ is an \mbox{$L_{\infty}$}\=/algebra. Similarly in hypothesis of Corollary \ref{tfp} $\left( \left(\widetilde{Q}^0\right)_{\bullet} := \widetilde{Q}_{\bullet}^0, \, \{{\widetilde{q}_n}|_{\widetilde{Q}^0}\}\right)$ is an \mbox{$L_{\infty}$}\=/algebra.
\end{pp}

%
%

\subsubsection{Proof of the theorem}

\noindent{\bf Proof of Theorem \ref{maintf}:}
$Q$ is a graded vector space concentrated in positive degrees by definition. The brackets are linear and they are claimed to be of bi\=/degree (-1,0). Let us look at the degree of the maps.

Let us suppose that $Q_{i\geqslant2}^{j} \neq 0$, that means that $Q_{i\geqslant 2}^{j} = D^{m-(i+j)}$ and hence $\left(q_1\right)_{i\geqslant 2}^{j} = -d$ goes to $(D^{f})^{m-(i+j)+1} \subset Q_{i-1}^j$. The fact that it ends up in $D^{f}$ and not only $D$ is because $d(d \alpha) = 0 \in \textrm{im}(f)$.

The maps $q_{k\geqslant 2}$ are also of bi\=/degree $(-1,0)$. We fix $k \geqslant 2$ and $t$ a positive integer. We start at $Q_1^{j_1} \otimes \cdots \otimes Q_1^{j_k} \subset Q_t$ (where $k + \sum_{i=1}^k j_i = t$):

\begin{center}
	\begin{tikzpicture}[description/.style={fill=white,inner sep=2pt},descr/.style={fill=white,inner sep=2.5pt}]
	\matrix (m) [matrix of math nodes, row sep=1.5em,
	column sep=2.5em, text height=1.5ex, text depth=0.25ex]
	{\left(D^{f}\right)^{m-1-j_1} \otimes \cdots \otimes \left(D^{f}\right)^{m-1-j_1}&Q_1^{j_1} \otimes \cdots \otimes Q_1^{j_k}&Q_{k-1}^{\sum_{i=1}^k j_i}\\
		C(\mathfrak{g})^{m-j_1} \otimes \cdots \otimes C(\mathfrak{g})^{m-j_k}&C(\mathfrak{g})^{(1-k)(m+1) + km -\sum_{i=1}^k j_i} & \left(D^{f}\right)^{m-(t-1)}
		\\};
	\path[->,font=\scriptsize]
	(m-1-1) edge node[auto] {$=$} (m-1-2)
	(m-1-2) edge node[auto] {$q_k$} (m-1-3)
	(m-1-1) edge node[auto] {$p^{\otimes k}$} (m-2-1)
	(m-2-1) edge node[auto] {$\dec(\mu^k)$} (m-2-2)
	(m-2-2) edge node[auto] {$f$} (m-2-3)
	(m-2-3) edge node[auto] {$\subset$} (m-1-3);the following equality holds
	\end{tikzpicture}
\end{center}

Observe that $p \colon \left(D^f\right)^{m-1-j} \longrightarrow C^{m-j}$ can be understood as a map from $Q_1^j$ to $\left(S^{\bullet}\mathfrak{g}\right)_{1+j}$. This shows that {\bf the Koszul sign in $Q$ and in $S^{\bullet}\mathfrak{g}$} after $p$ {\bf are the same}. Since $\mu^k$ is graded symmetric in $S^{\bullet}\mathfrak{g}$ so it is $q_{k\geqslant 2}$ in $Q$.

We need to show that the Jacobi identity holds. We first observe that $q_{k} \circ_1 q_1 = 0$ this is clear for $k=1$ since $D$ was a cochain map. For a larger $k$ we have that $f \circ p \circ d = d \circ d = 0 = f(0)$ so that $f \circ_1 (p \circ (-d)) = f \circ_1 0 = 0$. By linearity of $\mu$ and $f$ we conclude that $q_{k} \circ_1 q_1 = 0$. This shows that $J(1)=J(2)=0$ ($q_1 \circ q_2 =0$ by horizontal degree reasons). We fix $n\geqslant 3$ and write as usual:
$$J(n)=\sum_{i+j=n+1} q_j \ins q_i = q_{n} \ins q_1 + q_{n-1} \ins q_2 + \left( \sum_{\begin{smallmatrix}i+j=n+1 \\ 3 \leqslant i \leqslant n-1 \end{smallmatrix}} q_j \ins q_i\right) + q_1 \ins q_{n}.$$

We have that $\sum_{\begin{smallmatrix}i+j=n+1 \\ 3 \leqslant i \leqslant n-1 \end{smallmatrix}} q_j \ins q_i = 0$ by horizontal degree reasons as in the results about hamiltonian forms and pairs. $J(n)$ reduces to
$$J(n)= q_{n-1} \ins q_2 + q_1 \ins q_{n}.$$
Any of those summands is zero unless all the inputs are in horizontal degree $1$. Let us have a look on the term $q_{n-1} \ins q_2$ in that case. For clearness of the proof we are going to denote $q_{k\geqslant 2}$ just by $f(p \stackrel{\scriptscriptstyle{k}}{\cdots} p)$ (all the $\smallsetminus[m+1](\mu^k)$ are here omitted). Now we have that $f \circ p \circ f = d \circ f = f \circ c$ so that $f \circ_1 (p \circ f) = f \circ_1 (c)$. Let us fix $y=y_1 \otimes \cdots \otimes y_n \in (Q_{1}^{\bullet})^{\otimes n}$.
\begin{eqnarray*}
	(q_{n-1} \ins q_2)(y) &=& \sum_{\sigma \in \mathcal{S}h_{n-2}^2} \epsilon(\sigma) f\left( p\left(f(p (y_{\sigma(1)}) p(y_{\sigma(2)}))\right) p(y_{\sigma(3)}) \cdots p(y_{\sigma(n)})\right) \\
	&=& \sum_{\sigma \in \mathcal{S}h_{n-2}^2} \epsilon(\sigma) f(c\left(p y_{\sigma(1)} p y_{\sigma(2)})\right) p(y_{\sigma(3)}) \cdots p(y_{\sigma(n)}))\\
	&=& \left( (f \circ c) \circ p^{\otimes n} \right) (y)= \left( (d \circ f) \circ p^{\otimes n} \right) (y) = -(q_1 \ins q_n)(y).
\end{eqnarray*}

Here it is very important to observe that $c = \smallsetminus[m+1](\hat{\nu})$ so that $c$ on two elements is just $\smallsetminus[m+1](\nu)$ and the sum over the shuffles gives again $c$.

For the proof of the strict \Linfm definded by $p$ we first need to see that it defines a symmetric linear map of degree $0$ and then we will modify $p$ so that the \Linfm condition is satisfied. Recall that $p \colon \left(D^f\right)^{m-1-j} \longrightarrow C(\mathfrak{g})^{m-j}$ can be understood as a map from $Q_1^j$ to $\left(S^{\bullet}\mathfrak{g}\right)_{1+j}$ which shows that $p \colon Q \longrightarrow \mathfrak{g}$ is a symmetric linear map of degree $0$.

It is important to notice that the definition of $q$ is independent of the choice of $p$ as far as $f \circ p = d$. Suppose we are given $p^{\prime} \in \underline{\textrm{Hom}}_{{\sf grVec}}(D^{f},C(\mathfrak{g}))_1$ such that $f \circ p^{\prime} = d = f \circ p$. This means that $f \circ_1 p^{\prime} = f \circ_1 p$ by the first hypothesis and hence the brackets defined by $p^{\prime}$ and $p$ coincide. That is why we are going to force $p$ to be $0$ on exact froms $d \circ d \alpha = 0 = f (0) = f (p \alpha)$.

Now it is clear that $p \circ q_n$ is zero for $n \neq 2$, precisely the case in which $\nu_n \neq 0$.  In the case $n=2$ we have that $f \circ \nu(p\alpha_1, p\alpha_2) = d \circ f(p\alpha_1, p\alpha_2)$ so that we can define $p$ such that $p(f(p\alpha_1, p\alpha_2)) = \nu(p\alpha_1, p\alpha_2)$, that is $p \circ q_2 = \nu \circ (p \otimes p)$ proving that $p$ is a stric \Linf and concluding the proof.
\qed\\


\subsection{Transfer and morphisms}\label{53}

In the previous section we transfered the \Linf structure from a Lie$[1]$\=/algebra $(\mathfrak{g}, \nu)$ to a cochain complex $D$ thanks to cochain homotopy $f$ satisfying certain conditions.

In this subsection we study the relation between \Linf morphisms from an \Linf $L$ towards $Q$ (synchronized ones) and cochain maps (homotopies) from the Chevalley\=/Eilenberg cochain complex associated to $L$ to $D$. As we did for $\mathfrak{g}$, we fix a notation for the shifted and truncated Chevalley\=/Eilenberg of $L$: $(C(L), e)$ is shifted by $m+1$ ($C(L)^n = (S^{\bullet} L)_{m+1-n}$ and $e = \smallsetminus[m+1](\sum_{n}\hat{l}_n)$). 

We begin by establishing a relationship between synchronized lifts of a map to $\mathfrak{g}$ and synchronized \Linf maps to $Q$:

\begin{lm}\label{lmtf}
In the notation of Corollary \ref{tfp} let $\mathsf{a} \colon (L,\{l_n\}) \longrightarrow \mathsf{tr}_{\{1,\ldots,m\}}(\mathfrak{g}, \nu)$ be an \mbox{$L_{\infty}$}\=/morphism. 

There is a one\=/to\=/one correspondence between synchronized \mbox{$L_{\infty}$}\=/morphisms $\tilde{\mathsf{h}}$ from $L$ to $\widetilde{Q}$ lifting ${\sf a}$ and synchronized \mbox{$L_{\infty}$}\=/morphisms $\mathsf{h} = \phi \circ \tilde{\mathsf{h}} $ from $L$ to $Q$ with the property that $f \circ {\sf a}_1 =d \circ \mathsf{h}_1$. 

If $L$ is concentrated in degree $1$ we have the same one\=/to\=/one correspondence between \mbox{$L_{\infty}$}\=/morphisms $\tilde{\mathsf{h}}$ from $L$ to $\widetilde{Q}^0$ lifting ${\sf a}$ and \mbox{$L_{\infty}$}\=/morphisms $\mathsf{h} = \phi \circ \tilde{\mathsf{h}}$ from $L$ to $Q^0$ with the property that $f \circ {\sf a}_1 =d \circ \mathsf{h}_1$. 
\end{lm}

The condition $f \circ {\sf a}_1 =d \circ \mathsf{h}_1$ sais that $\mathsf{h}_1(L) \subset D^f$, in other words, that we can view the elements in $L$ to be hamiltonian. We are interested in lifts of the kind

\begin{center}\leavevmode
\xymatrix{ &&\widetilde{Q}\ar^{p}[d]\\
L\ar^{\mathsf{a}}[rr]\ar^{\mathsf{h}}@{.>}[urr]&&\mathfrak{g}}
\end{center}

\begin{pp}\label{strict2}
	Let $(L, \{l_n\})$ be an \Linf and let ${\sf a} \colon (L,\{l_n\}) \longrightarrow \mathsf{tr}_{\{1,\ldots,m\}}(\mathfrak{g}, \nu)$ be an \Linfm such that there exists a synchronized \mbox{$L_{\infty}$}\=/morphisms $\mathsf{h}$ from $L$ to $Q$ with the property that $f \circ {\sf a}_1 =d \circ \mathsf{h}_1$. 
	Then ${\sf a}$ is strict.
\end{pp}

\dem
The existence of a synchronized lift implies that $\mathsf{a}_{n\geq 2} = p \circ \mathsf{h}_{n \geq 2} = 0$ since $\mathsf{h}_n(L)$ is concentrated in horizontal degree $n$ where $p$ is defined to be zero. This fact has very relevant consequences: the \Linfm conditions on $\mathsf{a}$ are simply $\mathsf{a}_1 (l_2(x_1, x_2)) = \nu(\mathsf{a}_1 x_1, \mathsf{a}_1 x_2)$ for every $x_1 \otimes x_2 \in \, L^{\otimes 2}$ and $\mathsf{a}_1 \circ l_{n\neq 2} = 0$ (in particular $\mathsf{a}_1 \circ l_1 =0$). That shows the claim.
\qed

In general $\mathsf{a}(L)$ is concentrated in degrees $\{1, \ldots, m+1\}$ but the existence of a synchronized lift tells us that $\mathsf{h}(L)$ is concentrated in degrees $\{1, \ldots, m\}$ and hence  $\mathsf{a}(L)=\mathsf{a}_1(L)=(p \circ \mathsf{h})(L)$. Moreover, if $L$ is a honest Lie algebra, $\mathsf{a}$ is a Lie algebra morphism (up to a shift).

\subsubsection{From cochain maps to \Linf morphisms}

Recall that, with the notation for the shifted Chevalley\=/Eilenberg complex introduced at the begining of the subsection, the \Linfm $\mathsf{a}$ can be viewed as a cochain map $\hat{\mathsf{a}} \colon C(L) \longrightarrow C(\mathfrak{g})$. Composing with  $f$ we get a new cochain map $\mathsf{f} := f \circ \hat{\mathsf{a}} \colon C(L) \longrightarrow D$. In the case $\mathsf{a}$ is concentrated in $\mathsf{a}_1$ (as in the hypothesis of Proposition \ref{strict2}), the equation is explicitly given by: 
$$\mathsf{f}(x_1, \ldots, x_j) := f (\mathsf{a}_1 x_1 \ldots \mathsf{a}_1 x_j) \textrm{ for every } x_1 \cdots x_j \in (S^j L)_{m+1-i} \subset C(L)^i .$$

Knowing about the main transfer Theorem \ref{maintf}, we are interested in extending $\mathsf{f}$ to get a cochain homotopy $(\mathsf{h}, \mathsf{f})$. In this section we are trying to answer when such a cochain homotopy gives rise to an \Linf algebra morphism towards the transferred \Linf $Q$ (or more in general to $\widetilde{Q}$ in order to introduce compatibility relations with $\mathsf{a}$).

\begin{tm}[Transfer result for morphisms]\label{tfmp}\index{Transfer!morphisms}
	Let $(\mathsf{h}, \mathsf{f}) \colon C(L) \longrightarrow D$ be a cochain null\=/homotopy such that $(\mathsf{h} \circ e) (L) = 0$ and that $\mathsf{f}=f(\mathsf{a}_1, \ldots, \mathsf{a}_1)$. We use the notation $\mathsf{h}_j^i \colon (S^j L)_{m+1-i} \longrightarrow D^{i-1}$.
	
	We can view $\hat{\mathsf{h}}_j^i := \mathsf{h}_j^{m+1-i} \colon (L^{\otimes j})_i \longrightarrow D^{m-i} \subset Q_j^{i-j} \subset Q_i$.
	Then $\hat{\mathsf{h}} = \{\hat{\mathsf{h}}_j\}$ is a synchronized \Linfm from $L$ to $Q$ and it defines an \Linfm to $\widetilde{Q}$.
\end{tm}

\dem The last statement follows from the first. If $\hat{\mathsf{h}}$ is a synchronized \Linfm and $(\mathsf{h} \circ e) (L) = 0$ then $\mathsf{f} \circ {\sf a}_1 =d \circ \mathsf{h}_1$ by the homotopy condition. Applying Lemma \ref{lmtf} we get an \Linfm to $\widetilde{Q}$. 

$\hat{\mathsf{h}}$ is linear and symmetric by assumption (it was a cochain map, hence linear, and it was defined on the symmetric power, hence symmetric). It is also of degree zero and synchronized since $\hat{\mathsf{h}}_j^i \colon (L^{\otimes j})_i \longrightarrow Q_j^{i-j}$.

We would have proved the theorem as soon as we will prove that for every $n \in \mathbb{Z}$, $n \geqslant 1$ and every $x = x_1 \otimes \cdots \otimes x_n$ in $L^{\otimes n}$:
\beq\label{gamma}
\sum_{i+j =n +1} (\hat{\mathsf{h}}_j \ins l_i) (x) = \sum_{a_1 + \ldots a_b = n}^{a_k \geqslant 1} \sum_{\sigma \in \mathcal{S}{h}(a_1,\ldots, a_b)} \frac{1}{b!} (q_b \circ (\hat{\mathsf{h}}_{a_1} \otimes \cdots \otimes \hat{\mathsf{h}}_{a_b})) \, (\sigma \cdot x).
\eeq

By the definition of the co\=/differential in $C(L)$ the left hand side of the previous equation is just $\mathsf{h} \circ e$. On the right hand side we have to examine differently the cases $b=1$ and $b \geqslant 2$. So, firstly we have a look on the case where only $b=1$ is possible, this is when $n=1$:

$$\sum_{i+j =2} (\hat{\mathsf{h}}_j \ins l_i) (x) = \hat{\mathsf{h}}_1 \circ l_1 = 0 = q_1 \circ \hat{\mathsf{h}}_1 = \sum_{\begin{smallmatrix} a_1 + \ldots a_b = 1 \\ \sigma \in \mathcal{S}{h}(a_1,\ldots, a_b) \end{smallmatrix}}^{a_k \geqslant 1} \frac{1}{b!} (q_b \circ (\hat{\mathsf{h}}_{a_1} \otimes \cdots \otimes \hat{\mathsf{h}}_{a_b})) \, (\sigma \cdot x).$$

The zero on the left hand side is by assumption, as we have said on the first paragraph of this proof, and on the right hand side by horizontal degree reasons.

If $n\geqslant 2$ we have to deal with $b\geqslant 2$. $q_b$ only does not vanish if all the inputs have horizontal degree equal to $1$. $\hat{\mathsf{h}}_{a_k}$ takes values in horizontal degree $a_k$ since it is synchronized. The only case where $q_b \circ (\hat{\mathsf{h}}_{a_1} \otimes \cdots \otimes \hat{\mathsf{h}}_{a_b})$ is not identically zero is when $b=n$, in that case:

\begin{eqnarray*}
	&& \sum_{\sigma \in \mathcal{S}{h}_0^n} \frac{1}{n!} (q_n \circ (\hat{\mathsf{h}}_{1} \otimes \cdots \otimes \hat{\mathsf{h}}_{1})) \, (\sigma \cdot x) = \\
	&=& (q_n \circ (\mathsf{h}_1)^{\otimes n}) x = \left( f \circ (p \circ \mathsf{h}_1)^{\otimes n}\right) (x) = \mathsf{f}(x).
\end{eqnarray*}

Observe that because $\mathsf{f} \circ {\sf a}_1 =d \circ \mathsf{h}_1$ we can take $p$ such that $p \circ \mathsf{h}_1 = \mathsf{a}_1$. This shows that $\left( f \circ (p \circ \mathsf{h})^{\otimes n}\right) (x) = \mathsf{f}(x)$.

On the other hand when $b=1$ we have $(q_1 \circ \mathsf{h}_n)(x) = (-d \mathsf{h}) (x)$. At the case $n\geqslant 2$ the equation \ref{gamma} is:
$$\mathsf{h} \circ e = \mathsf{f} - d \circ \mathsf{h},$$
which holds since $(\mathsf{h},\mathsf{f})$ is a null\=/homotopy.
\qed\\

Observe that the condition $(\mathsf{h} \circ e) (L) = 0$ is necessary. In order to get an \Linfm lift to $\widetilde{Q}$ we need $\mathsf{a}$ to be just $\mathsf{a}_1$ as argued in Proposition \ref{strict2}, and that the homotopy $\mathsf{h}$ satisfies the condition $\mathsf{f} \circ {\sf a}_1 =d \circ \mathsf{h}_1$. {\it A priori} we only have $\mathsf{f} \circ {\sf a}_1 =d \circ \mathsf{h}_1 + \mathsf{h} \circ e$ so that $\mathsf{h}$ defines an \Linf lift if and only if $(\mathsf{h} \circ e) (L) = 0$.

That clearly holds if $l_1 = 0$ so that the transfer result for morphisms from a Lie$[1]$\=/algebra holds trivially.

\begin{cl}[Transfer result for morphisms for Lie\mbox{[1]}\=/algebras]\label{cltfmp}\index{Transfer!morphisms!for Lie$[1]$\=/algebras}
	Given a cochain null\=/homotopy $(\mathsf{h}, \mathsf{f}) \colon C(L) \longrightarrow D$ from an \Linf $L$ concentrated in degree $1$ such that $\mathsf{f}=f(\mathsf{a}_1, \ldots, \mathsf{a}_1)$. 
	
	We can view $\hat{\mathsf{h}}^i := \mathsf{h}_{i}^{m+1-i} \colon (L^{\otimes i})_{i} \longrightarrow D^{m-i} \subset Q_{i}^{0} = \left(Q^0\right)_i$.
	Then $\hat{\mathsf{h}} = \{\hat{\mathsf{h}}_j\}$ is a synchronized \Linfm from $L$ to $Q^0$ and it defines an \Linfm to $\widetilde{Q}^0$.
\end{cl}

%
%
%

\subsubsection{From \Linf morphisms to cochain maps}

We get converse theorems to the previous ones:
 
\begin{pp}[Converse transfer result for morphisms]\label{maincv}\index{Transfer!converse result}
Let $\tilde{\mathsf{h}} \colon L \longrightarrow \widetilde{Q}$ be a synchronized $L_{\infty}$\=/morphism lifting $\mathsf{a} \colon L \longrightarrow \mathfrak{g}$. Let $\mathsf{h} \colon L \longrightarrow Q$ be the associated \Linfm given by Lemma \ref{lmtf}. We view
$$\hat{\mathrm{h}}_j^{i} := \mathrm{h}_{j}^{m+1-i} \colon (S^j L)_{m+1-i} \subset C(L)^{i} \longrightarrow Q_j^{m+1-i-j} \subset D^{i-1}.$$
Then $(\hat{\mathsf{h}},\mathsf{f}:= f \circ \hat{\mathsf{a}})$ is a homotopy from $(C(L),e)$ to $(D,d)$.
\end{pp}



\dem Since $\mathsf{a}$ and $f$ are linear of degree $0$, $\mathsf{f}$ is also linear of degree $0$. 

For $n \geqslant 2$ the fact that $\mathsf{h}$ is a synchronized \Linfm translates into:
$$\sum_{i+j =n +1} (\mathsf{h}_j \ins l_i) = (q_1 \circ \mathsf{h}_n) (x) + (q_n \circ (\mathsf{h}_{1} \otimes \cdots \otimes \mathsf{h}_{1})) (x),$$
which, by the observation previous to this proof, is exactly
\begin{equation}\label{indepence}
\hat{\mathsf{h}}\circ e = -d \circ \hat{\mathsf{h}} +\mathsf{f}.
\end{equation}
That shows the cochain null\=/homotopy at every level except when only one input is considered. But in that case $f \circ {\sf a}_1 =d \circ \hat{\mathsf{h}}_1$ and the \Linfm condition for $n=1$ gives $\hat{\mathsf{h}} \circ e = 0$ by horizontal degree reasons. This shows that in this case $\mathsf{f} = d \circ \hat{\mathsf{h}} + \hat{\mathsf{h}} \circ e$ is also true.

It remains to be checked that $\mathsf{f}$ is a cochain map. But this follows from the null\=/homotopy condition:
$$ \mathsf{f} \circ e = (d \circ \hat{\mathsf{h}} + \hat{\mathsf{h}} \circ e) \circ e = d \circ \hat{\mathsf{h}} \circ e = d\circ (d \circ \hat{\mathsf{h}} + \hat{\mathsf{h}} \circ e) = d \circ \mathsf{f}.$$
\qed

Observe that it follows from equation \ref{indepence} that the definition of $\mathsf{f}$ is independent of $f$.

As usual we a simplification if the image of $\mathsf{h}$ is concentrated in horizontal degree $0$, that is when $h_{n\geqslant 2} = 0$. A sufficient condition for this to happen is that $L$ is concentrated in degree $1$.

\begin{cl}[Converse transfer result for morphisms for Lie\mbox{[1]}\=/algebras]\label{clcv}\index{Transfer!converse result!for Lie$[1]$\=/algebras}
Let $L$ be concentrated in degree $1$ and $\tilde{\mathsf{h}} \colon L \longrightarrow \widetilde{Q}^0$ be a synchronized $L_{\infty}$\=/morphism lifting $\mathsf{a} \colon L \longrightarrow \mathfrak{g}$. Let $\mathsf{h} \colon L \longrightarrow Q^0$ be the associated \Linfm given by Lemma \ref{lmtf}. We can view
$$\hat{\mathrm{h}}_i := \mathrm{h}^{m+1-i} \colon (S^\bullet L)_{m+1-i} = S^{m+1-i} L = C(L)^{i} \longrightarrow Q_{m+1-i}^{0} \subset D^{i-1}.$$
Then $(\hat{\mathsf{h}},\mathsf{f})$ is a homotopy from $(C(L),e)$ to $(D,d)$.
\end{cl}

\dem
We can compose $\tilde{\mathsf{h}}$ with the inclusion from $\widetilde{Q}^{0}$ to $\widetilde{Q}$ and apply Proposition \ref{maincv}.
\qed


\subsubsection{Comparison between both constructions}\label{54}

The immediate questions after the previous section are the following: 
\begin{enumerate}
\item Are Theorem \ref{tfmp} and Proposition \ref{maincv} inverse to each other?
\item Are Corollary \ref{cltfmp} and Corollary \ref{clcv} inverse to each other?
\end{enumerate}

If we start from an \Linf lift and we apply Proposition \ref{maincv} we see (in the proof) that the new homotopy $\mathsf{h}$ is such that $\mathsf{h} \circ e \, (L)= 0$. By Proposition \ref{strict2}, the fact that $\mathsf{a}$ factors through $\widetilde{Q}$ implies that $\mathsf{a}$ is strict. Therefore, $\mathsf{a}$ is simply $\mathsf{a}_1$ and $\mathsf{f} := f \circ \hat{\mathsf{a}} = f (\mathsf{a}_1 \ldots \mathsf{a}_1)$. We can conclude that the new cochain homotopy is in the hypothesis of Theorem \ref{tfmp} and hence that the hypothesis of that Theorem are minimal: we get all lifts starting from that kind of cochain homotopies.

\begin{cl}\label{central}
Let $\mathsf{a} \colon L \longrightarrow \mathsf{tr}_{\{1,\ldots,m\}}\mathfrak{g}$ be an \Linfm and $f \colon C(\mathfrak{g}) \longrightarrow D$ a cochain map  in the notation of the previous section. Then there is a one\=/to\=/one correspondence between:
\begin{enumerate}
	\item Synchronized \Linf lifts of $\mathsf{a}$ to $\widetilde{Q}$ the \Linf constructed in Theorem \ref{maintf} and
	\item Cochain null\=/homotopies $(\mathsf{h}, \mathsf{f}:=f\circ \hat{\mathsf{a}}) \colon C(L) \longrightarrow D$ such that $\mathsf{f} = f (\mathsf{a}_1, \ldots, \mathsf{a}_1)$ and $(\mathsf{h} \circ e)(L)=0$.
\end{enumerate}
\end{cl}

This result provides a positive answer to the first question above. Observe that the equivalence does not hold for all cochian homotopies, but only for thise that satisfy the conditions $\mathsf{f} = f (\mathsf{a}_1, \ldots, \mathsf{a}_1)$ (which does not depend on $\mathsf{h}$) and $(\mathsf{h} \circ e)(L)=0$ (which does depend on $\mathsf{h}$).

In the case of Lie$[1]$\=/algebras, the condition $(\mathsf{h} \circ e)(L)=0$ disappears as we saw in Corollary \ref{cltfmp}. 

\begin{cl}\label{central1}
Let $\mathsf{a} \colon L \longrightarrow \mathsf{tr}_{\{1, \ldots, m\}}\mathfrak{g}$ be an \Linfm from an \Linf concentrated in degree $1$ such that $\mathsf{f}:= f \circ \hat{\mathsf{a}} = f (\mathsf{a}_1, \ldots, \mathsf{a}_1) \colon L \longrightarrow D$ where $f \colon C(\mathfrak{g}) \longrightarrow D$ is the cochain map in the notation of the previous section. Then there is a one\=/to\=/one correspondence between:
\begin{enumerate}
	\item Synchronized \Linf lifts of $\mathsf{a}$ to $\widetilde{Q}^0$, the \Linf constructed in Theorem \ref{maintf}, and
	\item Cochain null\=/homotopies $(\mathsf{h}, \mathsf{f}) \colon C(L) \longrightarrow D$.
\end{enumerate}
\end{cl}

But for higher \Linfs the condition $(\mathsf{h} \circ e)(L)=0$ can fail. There are null\=/homotopies such that $(\mathsf{h} \circ e)(L)\neq 0$ so that the class of null\=/homotopies is larger than the ones representing \Linf lifts.

There are ways of forcing $\mathsf{h}$ to satisfy $(\mathsf{h} \circ e)(L)= 0$ withouth involving $\mathsf{h}$ itself. An example is that $L$ is concentrated in degree $1$ as in the previous Corollary. Another one is to require $l_1$ to be zero.

\begin{cl}\label{centrall1}
Let $\mathsf{a} \colon L \longrightarrow \mathsf{tr}_{\{1, \ldots, m\}}\mathfrak{g}$ be an \Linfm such that $l_1 = 0$ and that $\mathsf{f}:= f \circ \hat{\mathsf{a}} = f (\mathsf{a}_1, \ldots, \mathsf{a}_1) \colon L \longrightarrow D$ where $f \colon C(\mathfrak{L}) \longrightarrow D$ is the cochain map in the notation of the previous section. Then there is a one\=/to\=/one correspondence between:
\begin{enumerate}
	\item Synchronized \Linf lifts of $\mathsf{a}$ to $\widetilde{Q}$, the \Linf constructed in Theorem \ref{maintf}, and
	\item Cochain null\=/homotopies $(\mathsf{h}, \mathsf{f}) \colon C(L) \longrightarrow D$.
\end{enumerate}
\end{cl}


\color{BlueViolet}
\section{Momentum maps in multi\=/symplectic geometry}
\color{black}

The \Linf of Hamiltonian forms is a special example of the transferred \Linf $Q$ from last chapter. As a matter of fact, we can give a geometric interpretation of cochain null\=/homotopies and synchronized \Linf morpihsms towards $Q$, therefore using the results from last section into the world of multi\=/symplectic geometry.

From a purely multi\=/symplectic geometric view this chapter is devoted to study the possible higher analogues of co\=/momentum maps is multi\=/symplectic geometry.

The main motivation for this is the paper on homotopy moment maps \cite{FRZ}. It is possible to understand that paper in terms of cochain homotopies and then to extend the result to include higher Hamiltonian pairs. That observation is original of Blohmann (private communication). He observed that the results in \cite{FRZ} can be understood to be the following: momentum maps from a Lie algebra are cochain homotopy from the Chevalley\=/Eilenberg cochain complex to differential forms such that the generating function condition is satisfied. We elaborate on that idea: 

Let $(M, \omega)$ be a symplectic manifold and $\mathfrak{g}$ a Lie algebra acting on $M$ by symplectomorphisms (i.e. we have a Lie algebra morphism ${\sf a} \colon \mathfrak{g} \longrightarrow \mathfrak{X}(M)$). A {\bf co\=/momentum map} (see \cite[Section 18.1]{C} for example) is a morphism of Lie algebras $\mu \colon \mathfrak{g} \longrightarrow \mathcal{C}^{\infty}(M)$ such that $\iota_{\mathsf{a}x} \omega = d \mu(x)$ for all $x \in \, \mathfrak{g}$ (that condition is called the {\it generating function condition}). In particular, the image of $\mathsf{a}$ is a subset of hamiltonian vector fields:
$$\mathsf{a} \colon \mathfrak{g} \longrightarrow \mathfrak{X}_{\mathrm{ham}}^1(M).$$

We can understand $\mu$ in in two different ways. 
\begin{enumerate}
	\item From one side as an \Linfm from $\mathfrak{g}$ to $\Xi^0(M)$ (remark that $m=1$ in the case of symplectic geometry) which satisfies the extra condition $d \mu (x) = \iota_{{\sf a}(x)} \omega$. This is equivalent (due to Lemma \ref{lmlift}) to an \Linfm from $\mathfrak{g}$ to $\widetilde{\Xi}^0(M)$ lifting ${\sf a}$.
	\begin{center}
		\begin{tikzpicture}[description/.style={fill=white,inner sep=2pt}]
		\matrix (m) [matrix of math nodes, row sep=1.5em,
		column sep=3.5em, text height=1.5ex, text depth=0.25ex]
		{ & \widetilde{\mathfrak{X}}^{1}(M) \\
			\mathfrak{g} & \mathfrak{X}^{1}(M)\\};
		\path[->,font=\scriptsize]
		(m-2-1) edge node[auto] {${\sf a}$} (m-2-2)
		(m-1-2) edge node[auto] {$\pi_1$} (m-2-2)
		(m-2-1) edge node[auto] {$\tilde{h}$} (m-1-2);
		\end{tikzpicture}
	\end{center}
	\item Another point of view is to think of it as a cochain homotopy as follows:
	\begin{center}
		\begin{tikzpicture}[description/.style={fill=white,inner sep=2pt}]
		\matrix (m) [matrix of math nodes, row sep=1.5em,
		column sep=3.5em, text height=1.5ex, text depth=0.25ex]
		{\mathfrak{g} \wedge \mathfrak{g} & \mathfrak{g} \\
			\mathcal{C}^{\infty}(M) & \Omega^1(M) \\};
		\path[->,font=\scriptsize]
		(m-1-1) edge node[auto] {$[-,-]$} (m-1-2)
		(m-2-1) edge node[below] {$d$} (m-2-2)
		(m-1-2) edge node[auto] {$\mu$} (m-2-1)
		(m-1-1) edge node[left] {$\widetilde{\omega} \circ ({\sf a} \otimes {\sf a})$} (m-2-1)
		(m-1-2) edge node[auto] {$\widetilde{\omega} \circ {\sf a}$} (m-2-2);
		\end{tikzpicture}
	\end{center}
\end{enumerate}

Both definitions easily generalize to the set up of multi\=/symplectic geometry, in the case where we have a higher \Linf acting on a multi\=/symplectic manifold. As a consequence of last section, we will see in the following sections that the two previous interpretations are equivalent.

%
%
%
%


\subsection{\Linfb actions and symplectic co\=/momentum maps}\label{61}

Let $(M, \omega)$ denote a pre\=/m\=/symplectic manifold. Let $L$ be an \mbox{${L_{\infty}}$}\=/algebra.

We can think of defining an \Linf action of $L$ on $(M, \omega)$ as an \Linfm ${\sf a} \colon L \longrightarrow (\mathfrak{X}^{\bullet}(M), \nu_2)$, the problem is that this does not generalize the definition in symplectic geometry: we need to take a truncation of the target space:

\begin{df}[\Linfb action]\index{Linfinity algebra @ $L_{\infty}$\=/algebra!action}\label{hact}
Let $L$ be an \mbox{${L_{\infty}}$}\=/algebra. An \Linf action of $L$ on $(M, \omega)$ is an \Linfm 
$${\sf a} \colon L \longrightarrow \mathsf{tr}_{\{1, \ldots, m\}}\left((\mathfrak{X}^{\bullet}(M), \nu_2)\right).$$
\end{df}

Remember that in Chapter $2$ we gave a new \Linf structure on multi\=/vector fields that was not previously known nor used in the literature. We suggest a second definition.

\begin{df}[Strong \Linfb action]\index{Linfinity algebra @ $L_{\infty}$\=/algebra!strong action}
Let $L$ be an \mbox{${L_{\infty}}$}\=/algebra. A strong \Linf action of $L$ on $(M, \omega)$ is an \Linfm 
$${\sf a} \colon L \longrightarrow \mathsf{tr}_{\{1, \ldots, m\}}\left((\mathfrak{X}^{\bullet}(M), \{\nu_n\})\right).$$
\end{df}

We will use the notation $L \circlearrowleft (M,\omega)$ and $L \circlearrowleft_{\mathsf st} (M,\omega)$ to denote actions and strong actions of $L$ on the pre\=/multi\=/symplectic manifold $(M,\omega)$.

Observe that in the case of symplectic geometry the two notions agree:
$$\mathrm{tr}_{\{1\}}\left( (\mathfrak{X}^{\bullet}(M), \nu_2) \right) = \mathrm{tr}_{\{1\}}\left( (\mathfrak{X}^{\bullet}(M), \{ \nu_n \} ) \right) = (\mathfrak{X}^1(M), \nu_2).$$

An action (strong or not) $\mathsf{a} \colon L \circlearrowleft_{(\mathsf{st})} (M,\omega)$ on a pre\=/m\=/symplectic manifold is given by a family $\{\mathsf{a}_n \colon L^{\otimes n} \longrightarrow \mathsf{tr}_{\{1, \ldots, m\}}(\mathfrak{X}^{\bullet}(M))\}$. Observe that $\mathsf{a}_{n\geqslant(m+1)}$ is zero by degree reasons: $|\mathsf{a}_{n\geqslant(m+1)}(L^{\otimes n})|=|L^{\otimes n}| \geqslant n \geqslant m+1$ which means that $\mathsf{a}_{n\geqslant(m+1)}(L^{\otimes n}) \subset \mathsf{tr}_{\{1, \ldots, m\}}(\mathfrak{X}^{\bullet}(M))_{m+1} = 0$. In the particular case of pre\=/symplectic geometry ($m=1$) we have that $\mathsf{a}$ is strict: it is simply $\mathsf{a}_1$.

As a conclusion, strong \Linf actions from a Lie algebra on a symplectic manifold are the same as \Linf actions and as smooth actions.

Given an \Linf action (strong or not), we will call it {\bf multi\=/symplectic} if the image of $S^{\bullet}L$ consists of multi\=/symplectic multi\=/vector fields. We will say that the action is {\bf pre\=/hamiltonian} if the image of $S^{\bullet}L$ consists of hamiltonian multi\=/vector fields.


\subsection{Transfer in multi\=/symplectic geometry}

In this section we repeat what was done in subsections \ref{52} and \ref{53}, replacing  $\mathfrak{g}$ by $\mathfrak{X}_{\mathrm{sym}}^{\bullet}(M)$ and $D$ by $\mathsf{tr}_m(\Omega^{\bullet}(M))$ on a pre\=/$m$\=/symplectic manifold $(M, \omega)$. The results in that section were about \Linfms $\mathsf{a} \colon L \longrightarrow \mathsf{tr}_{\{1, \ldots, m\}}(\mathfrak{X}_{\mathrm{sym}}^{\bullet}(M), \nu_2)$. That means that, in this case, we will get information about multi\=/symplectic (non\=/strong) actions. We fix such an action from now on.

It is important to observe that we need the multi\=/vector fields to be {\it multi\=/symplectic}, otherwise $\widetilde{\omega}$ is not a cochain map.

The equivalent results of \ref{central} and \ref{central1} are the following:

\begin{cl}\label{centralms}
	Let $L \circlearrowleft (M,\omega)$ be a multi\=/symplectic action on the pre\=/m\=/symplectic manifold $(M, \omega)$. Then there is a one\=/to\=/one correspondence between:
	\begin{enumerate}
		\item Synchronized \Linf lifts of $\mathsf{a}$ to $\widetilde{\Xi}(M)$ and
		\item Cochain null\=/homotopies $(\mathsf{h}, \mathsf{f}:= \widetilde{\omega} \circ \hat{{\sf a}} = \widetilde{\omega} (\mathsf{a}_1, \ldots, \mathsf{a}_1)) \colon C(L) \longrightarrow \mathsf{tr}_m(\Omega^{\bullet}(M))$ such that $(\mathsf{h} \circ e)(L)=0$.
	\end{enumerate}
\end{cl} 

\begin{cl}\label{central1ms}
	Let $L \circlearrowleft (M,\omega)$ be a multi\=/symplectic action on the pre\=/m\=/symplectic manifold $(M, \omega)$ where $L$ is concentrated in degree $1$. Then there is a one\=/to\=/one correspondence between:
	\begin{enumerate}
		\item Synchronized \Linf lifts of $\mathsf{a}$ to $\widetilde{\Xi}^0(M)$ and
		\item Cochain null\=/homotopies $(\mathsf{h}, \mathsf{f}:= \widetilde{\omega} \circ \hat{{\sf a}} = \widetilde{\omega} (\mathsf{a}_1, \ldots, \mathsf{a}_1)) \colon C(L) \longrightarrow \mathsf{tr}_m(\Omega^{\bullet}(M))$.
	\end{enumerate}
\end{cl}

\begin{rk}\label{finalrk} 
	\par In subsection \ref{54} we saw that the obstructions for the two transfer theorems for morphisms (the direct and the converse one) to be the inverse of each other were $\widetilde{\omega} \circ \hat{\mathsf{a}} = \omega(\mathsf{a}_1, \ldots, \mathsf{a}_1)$ and $(\mathsf{h} \circ e)(L) = 0$.
	
	If the action is strict, the first of those two conditions is satisfied and the other one is conquered if $l_1 =0$. Both things happen in case $(M,\omega)$ is a pre\=/symplectic manifold and if $L$ is a non\=/graded Lie algebra.
\end{rk} 

We include the converse results and the proof of the general one for the sake of completeness:

\begin{pp}[Converse transfer result in multi\=/symplectic geometry]\label{cvms}\index{Transfer!converse result!in multi\=/symplectic geometry}
Let $\tilde{\mathsf{h}} \colon L \longrightarrow \widetilde{\Xi}(M)$ be a synchronized $L_{\infty}$\=/morphism lifting the action $\mathsf{a} \colon L \circlearrowleft M$. Let $\mathsf{h} \colon L \longrightarrow \Xi(M)$ be the associated \Linfm given by Lemma \ref{lmlift}. We view
$$\hat{\mathrm{h}}_j^{i} := \mathrm{h}_{j}^{m+1-i} \colon (S^j L)_{m+1-i} \subset C(L)^{i} \longrightarrow \Xi_j^{m+1-i-j}(M) \subset \Omega^{i-1}(M).$$
Then $\mathsf{a}$ is strict and $(\hat{\mathsf{h}},\mathsf{f}= \widetilde{\omega}(\mathsf{a}_1, \ldots, \mathsf{a}_1))$ is a homotopy from $(C(L),e)$ to $\mathsf{tr}_m(\Omega^{\bullet}(M))$.
\end{pp}

\begin{cl}\index{Transfer!converse result!in multi\=/symplectic geometry!for Lie$[1]\=/algebras$}\label{clcvms}
Let $\tilde{\mathsf{h}} \colon L \longrightarrow \widetilde{\Xi}^0(M)$ be a synchronized $L_{\infty}$\=/morphism lifting the action $\mathsf{a} \colon L \circlearrowleft M$ where $L$ is concentrated in degree $1$. Let $\mathsf{h} \colon L \longrightarrow \Xi^0(M)$ be the associated \Linfm given by Lemma \ref{lmlift}. We can view
$$\hat{\mathrm{h}}_{i} := \mathrm{h}^{m+1-i} \colon (S^\bullet L)_{m+1-i} = S^{m+1-i}L = C(L)^{i} \longrightarrow \Xi_{m+1-i-j}^0(M) \subset \Omega^{i-1}(M).$$
Then $\mathsf{a}$ is strict and $(\hat{\mathsf{h}},\mathsf{f}= \widetilde{\omega}(\mathsf{a}_1, \ldots, \mathsf{a}_1))$ is a homotopy from $(C(L),e)$ to $\mathsf{tr}_m(\Omega^{\bullet}(M))$.
\end{cl}

Those were \ref{maincv} and \ref{clcv} before. We include the proof of Proposition \ref{cvms}, since the other three results follow immediately from this one.

\dem
We want to apply Proposition \ref{maincv}. We take $\mathfrak{g} = \mathfrak{X}_{\mathrm{sym}}^{\bullet}(M)$ and $D=\mathrm{tr}_m(\Omega^{\bullet}(M))$. The map $\widetilde{\omega} \colon (C(\mathfrak{g}), d_{\textrm{CE}} ) \longrightarrow (D,d)$ is a cochain map by Proposition \ref{cochain} and $D^f$ is just $\mathrm{tr}_m(\Omega_{\mathrm{ham}}^{\bullet}(M))$ by definition of hamiltonian forms.

The map $p$ is no more than the choice of a hamiltonian multi\=/vector field. For any two multi\=/symplectic multi\=/vector fields $x$ and $y$ such that $\iota_x \omega = \iota_y \omega$ there are two possibilities: if $|x| \neq |y|$ then $\iota_x \omega = \iota_y \omega$ implies that $\iota_x \omega = \iota_y \omega=0$; if not $\widetilde{\omega}(x, z) = (-1)^{|x||z|} \iota_{z}\iota_x \omega = (-1)^{|y||z|} \iota_{z}\iota_y \omega= \widetilde{\omega}(y,z)$ for every $z$ product of multi\=/symplectic multi\=/vector fields. In any case $\widetilde{\omega}(x, -) = \widetilde{\omega}(y, -)$.

Since the action is multi\=/symplectic and it is concentrated in $\mathsf{a}_1$ we get the map $\mathsf{a} \colon L \longrightarrow \mathfrak{g}$ and we can apply Proposition \ref{maincv} to get the result.
\qed\\

\subsection{Momentum maps}

We are given an $m$\=/symplectic manifold $(M,\omega)$ and ${\sf a} \colon L \longrightarrow (\mathfrak{X}^{\bullet}(M), \nu_2)$ an \Linf action from an \Linf $(L,\{b_i\})$ concentrated in degrees $\{1, \ldots, r\}$ for some $r$ (possibly non finite). We have seen in Section \ref{61} that an \Linf lift of the action to $\widetilde{\Xi}(M)$ and a cochain null\=/homotopy $(\mathsf{h},\mathsf{f}=\widetilde{\omega}\circ\hat{\textsf{a}})$ to the complex of forms restricted to degrees $\{1, \ldots, m\}$ are good candidates to generalize the construction of a co\=/momentum map in symplectic geometry.

We begin with the definition of what a Homotopy Momentum Map is:

\begin{df}[Homotopy Momentum Map]\label{hmm}\index{Homotopy momentum!map}
Let $(L, \{l_i\})$ be an \Linf acting by ${\sf a}$ on the pre\=/$m$\=/symplectic manifold $(M,\omega)$. Let $C(L)^n := (S^{\bullet} L)_{m+1-n}$. If the action is multi\=/symplectic we get a cochain map $\mathsf{f} := \widetilde{\omega} \circ \hat{\mathsf{a}} \colon C(L) \longrightarrow \mathsf{tr}_m(\Omega^{\bullet}(M))$. A map $\mathsf{h} \in \, \underline{\textrm{Hom}}_{{\sf grVec}}(C(L), \Omega^{\bullet}(M))_{-1}$ is defined to be a Homotopy Momentum Map for the action if it defines a cochain null\=/homotopy
$$(\mathsf{h},\mathsf{f}) \colon C(L) \longrightarrow \mathsf{tr}_m(\Omega^{\bullet}(M)).$$
\begin{center}\leavevmode
\xymatrix{
0\ar[r]&(S^{\bullet}L)_{m+1} \ar^{\hat{l}}[r] \ar^{\widetilde{\omega}\circ\hat{\mathsf{a}}}[d]& (S^{\bullet}L)_{m}\ar^{\hat{l}}[r] \ar^{\widetilde{\omega}\circ\hat{\mathsf{a}}}[d] \ar^{\mathsf{h}}[dl]& \ldots \ar^{\hat{l}}[r] \ar^{\mathsf{h}}[dl] & (S^{\bullet}L)_{2} \ar^{\hat{l}}[r] \ar^{\widetilde{\omega}\circ\hat{\mathsf{a}}}[d] \ar^{\mathsf{h}}[dl]& (S^{\bullet}L)_{1} \ar[r] \ar^{\widetilde{\omega}\circ\hat{\mathsf{a}}}[d] \ar^{\mathsf{h}}[dl]&0\\ 
0 \ar[r]&\mathcal{C}^{\infty}(M)\ar^{d}[r] & \Omega^{1}(M) \ar^{d}[r] & \ldots \ar^{d}[r]& \Omega^{m-1}(M) \ar^{d}[r]& \Omega^{m}(M) \ar[r]&0 }
\end{center}
\end{df}

Observe that Homotopy Momentum Maps do not ensure the actions to be pre\=/hamiltonian. Eventhough they represent a trully generalization of the concept of co\=/momentum map in symplectic geometry, they go against the primary objective of momentum maps: having hamiltonian actions and giving the associated hamiltonian form to the image of $\mathsf{a}$.

That is the reason why we would ask, at least, the image of $L$ to be hamiltonian (if not the whole image of $S^{\bullet}(L)$). We give the following definition.

\begin{df}[Homotopy Hamiltonian Momentum Map]\label{hhmm}\index{Homotopy!hamiltonian!momentum!map}
	Let $\mathsf{h}$ be a Homotopy Momentum Map given by $(\mathsf{h},\mathsf{f}) \colon C(L) \longrightarrow \mathsf{tr}_m(\Omega^{\bullet}(M))$. $\mathsf{h}$ is called a Homotopy Hamiltonian Momentum Map if $(\mathsf{h} \circ l_1) (L)=0$.
\end{df}

If the action is concentrated in $\mathsf{a}_1$, the existence of a Homotopy Hamiltonian Momentum Map already ensures that the action was pre\=/hamiltonian.

\begin{pp}\label{stu}
	Let $\mathsf{a} \colon L \circlearrowleft (M,\omega)$ be an \Linf action. Let $\mathsf{h}$ be a Homotopy Hamiltonian Momentum Map for the action. If $\mathsf{a}$ is reduced to $\mathsf{a}_1$ then the action is pre\=/hamiltonian.
\end{pp}

\dem
Since $\mathsf{a}_{n\neq 1} = 0$ we have that $\iota_{\mathsf{a}_{n\neq 1}(x)}\omega = d(0)$ for every $x \in \, L^{\otimes n}$. Hence the action is, in that case, pre\=/hamiltonian if and only if $\iota_{\mathsf{a}_1 x}\omega$ is exact for every $x \in \, L$. 

But by assumption 
$$\iota_{\hat{\mathsf{a}}_1(x)} \omega = d (\mathsf{h}x) + (e \circ \mathsf{h}) x = d (\mathsf{h}) x + (\mathsf{h} \circ l_1)(x) = d (\mathsf{h}) x.$$ The action is then pre\=/hamiltonian.
\qed\\

Another point of view is given by the following definition:

\begin{df}[\mbox{\boldmath ${L_{\infty}}$}\=/Momentum Map]\label{lmm}\index{Linfinity mo @ ${L_{\infty}}$\=/momentum map}
Let $(L, \{l_i\})$ be an \Linf acting by ${\sf a}$ on the pre\=/$m$\=/symplectic manifold $(M,\omega)$. An \mbox{$L_{\infty}$}\=/Momentum Map, $\tilde{\mathsf{h}}$, is an \Linf lift of the action to $\widetilde{\Xi}(M)$:
\begin{center}
\begin{tikzpicture}[description/.style={fill=white,inner sep=2pt}]
\matrix (m) [matrix of math nodes, row sep=1.5em,
column sep=3.5em, text height=1.5ex, text depth=0.25ex]
{ & \widetilde{\Xi}(M) \\
L & \mathrm{tr}_{m}(\mathfrak{X}^{\bullet}(M)) \\};
\path[->,font=\scriptsize]
(m-2-1) edge node[auto] {${\sf a}$} (m-2-2)
(m-1-2) edge node[auto] {$\pi_1$} (m-2-2)
(m-2-1) edge node[auto] {$\tilde{\mathsf{h}}$} (m-1-2);
\end{tikzpicture}
\end{center}
\end{df}

{\bf Observe that the existence of an ${L_{\infty}}$\=/Momentum Map implies that the action is pre\=/hamiltonian (hence multi\=/symplectic) and strict.}

\begin{df}[Hamiltonian action]\index{Linfinity algebra @ $L_{\infty}$\=/algebra!action!hamiltonian}
Let $\mathsf{a} \colon L \circlearrowleft (M,\omega)$ be an \Linf action. We say that $\mathsf{a}$ is hamiltonian if there exists an ${L_{\infty}}$\=/Momentum for the action.
\end{df}

By Lemma \ref{lmlift} we know that an ${L_{\infty}}$\=/Momentum Map is the same as an \Linfm $h \colon L \longrightarrow \Xi(M)$ which satisfies the generating form condition. We know that the case where $r=1$ can be treated differently since the graded vector spaces $\Xi(M)$ and $\widetilde{\Xi}(M)$ are unnecessarily big.

\begin{df}[\mbox{\boldmath ${L_{\infty}}$}\=/Momentum Map]\label{lmm1}\index{Linfinity mo @ $L_{\infty}$\=/momentum map!for Lie$[1]$\=/algebras}
Let $(L, b_2)$ be a Lie$[1]$\=/algebra concentrated in degree $\{1\} $which is acting by ${\sf a}$ on the $m$\=/symplectic manifold $(M,\omega)$. An \mbox{$L_{\infty}$}\=/Momentum Map, $\tilde{\mathsf{h}}$ is a lift of the action to $\widetilde{\Xi}^0(M)$:
\begin{center}
\begin{tikzpicture}[description/.style={fill=white,inner sep=2pt}]
\matrix (m) [matrix of math nodes, row sep=1.5em,
column sep=3.5em, text height=1.5ex, text depth=0.25ex]
{ & \widetilde{\Xi}^0(M) \\
L & (\mathfrak{X}^{\bullet}(M),\nu_2) \\};
\path[->,font=\scriptsize]
(m-2-1) edge node[auto] {${\sf a}$} (m-2-2)
(m-1-2) edge node[auto] {$\pi_1$} (m-2-2)
(m-2-1) edge node[auto] {$\tilde{\mathsf{h}}$} (m-1-2);
\end{tikzpicture}
\end{center}
\end{df}

This Definition (\ref{lmm1}) is the \cite[Definition 5.1]{FRZ}, where they call it Homotopy Momentum Map. This apparently inconsistency of the terminology does not happen as we will see in the following paragraphs (observe that in \cite{FRZ} an action is an \Linfm towards $\mathfrak{X}^1(M)$ so that $\mathsf{a}$ is always reduced to $\mathsf{a}_1$). 

The main result of this section is that, in the hypothesis of Proposition \ref{stu}, not only the action is pre\-/hamiltonian but also hamiltonian. In particular we will see that for a non\=/graded Lie algebra acting on a pre\=/symplectic manifold the notions of ${L_{\infty}}$\=/Momentum Map, Homotopy Momentum Map and Homotopy Hamiltonian Momentum map coincide.

\begin{tm}\label{stu2}
Let $\mathsf{a} \colon L \circlearrowleft (M,\omega)$ be a strict \Linf action such that $(\mathsf{h} \circ l_1) (L)=0$. Then the following statements are equivalent:
\begin{enumerate}
	\item $\mathsf{h}$ is a Homotopy Momentum Map.
	\item $\mathsf{h}$ is a Homotopy Hamiltonian Momentum Map.
	\item $\phi \circ \mathsf{h}$ is an ${L_{\infty}}$\=/Momentum Map where $\phi \colon \widetilde{\Xi}(M) \longrightarrow \Xi(M)$.
\end{enumerate}
\end{tm}

\dem
The equivalence between $1$ and $2$ is given by Proposition \ref{stu} while the equivalence between $2$ and $3$ follows from Corollary \ref{centralms}.
\qed

As a Corollay we give obstructions for this last Theorem to hold independently of $\mathsf{a}$ and of $\mathsf{h}$:

\begin{cl}\label{eeoo}
Let $\mathsf{a} \colon L \circlearrowleft (M,\omega)$ be an \Linf action on a {\bf pre\=/symplectic} manifold such that $l_1 = 0$. Then the following statements are equivalent:
\begin{enumerate}
	\item $\mathsf{h}$ is a Homotopy Momentum Map.
	\item $\mathsf{h}$ is a Homotopy Hamiltonian Momentum Map.
	\item $\mathsf{h}$ is an ${L_{\infty}}$\=/Momentum Map.
\end{enumerate}
\end{cl} 

\dem
It follows from $(M, \omega)$ being pre\=/symplectic that $\mathsf{a}$ is reduced to $\mathsf{a}_1$. Since $l_1=0$ it is clear that $(\mathsf{h} \circ l_1) (L)=0$. Now we are in the hypothesis of Theorem \ref{stu2} which we can apply in order to get the desired result.
\qed\\

Obstructions for the existence of Homotopy Momentum Maps for Lie$[1]$\=/algebra actions have being given in \cite{FRZ}. A future line of investigation is the further study of the obstructions of the existence of Homotopy Momentum Maps in the general case (Definition \ref{hmm}) and of \mbox{$L_{\infty}$}\=/Momentum Maps (Definitions \ref{lmm} and \ref{lmm1}).

\subsection{Further remarks}

We have characterized all the synchronized \Linf lifts of multi\=/symplectic \Linf actions. But we have not included strong actions, nor characterized all cochain null\=/homotopies. We have the following list of exceptions:

	\begin{enumerate}
		
		\item If we start with a strong action $L \circlearrowleft_{\mathrm{st}} (M,\omega)$ we will not be able to talk about \Linfms lifting $\pi$ since the projection $\pi_{\mathrm{st}} \colon \widetilde{\Xi} (M)_{\mathrm{st}} \longrightarrow (\mathfrak{X}^{\bullet}(M), \{\nu_n\})$ is not an {\mbox{$L_{\infty}$}\=/morphism}). We could still talk about lifts outside of the category of \Linfs .
		
		\begin{df}[Strong \mbox{\boldmath ${L_{\infty}}$}\=/Momentum Map]\label{slmm}\index{Linfinity mo @ ${L_{\infty}}$\=/momentum map!strong}
			Let $(L, \{l_i\})$ be an \Linf stongly acting by ${\sf a}$ on the pre\=/$m$\=/symplectic manifold $(M,\omega)$. An \mbox{$L_{\infty}$}\=/Momentum Map, $\tilde{\mathsf{h}}$, is a lift (not an \Linf lift)of the action to $\widetilde{\Xi}(M)$:
			\begin{center}
				\begin{tikzpicture}[description/.style={fill=white,inner sep=2pt}]
				\matrix (m) [matrix of math nodes, row sep=1.5em,
				column sep=3.5em, text height=1.5ex, text depth=0.25ex]
				{ & \widetilde{\Xi}(M) \\
					L & \mathrm{tr}_{m}(\mathfrak{X}^{\bullet}(M)) \\};
				\path[->,font=\scriptsize]
				(m-2-1) edge node[auto] {${\sf a}$} (m-2-2)
				(m-1-2) edge node[auto] {$\pi_{\mathsf{st}}$} (m-2-2)
				(m-2-1) edge node[auto] {$\tilde{\mathsf{h}}$} (m-1-2);
				\end{tikzpicture}
			\end{center}
		\end{df}
		
		\item We can also talk about cochain maps $\mathsf{f}:= \widetilde{\omega} \circ \hat{\mathsf{a}} \colon C(L) \longrightarrow \mathsf{tr}_m(\Omega^{\bullet}(M))$ such that $\widetilde{\omega} \circ \hat{\mathsf{a}} \neq \widetilde{\omega}({\sf a}_1, \ldots, {\sf a}_1)$. Cochain null\=/homotopies $(\mathsf{h},\widetilde{\omega} \circ \hat{\mathsf{a}}) \colon C(L) \longrightarrow \mathsf{tr}_m(\Omega^{\bullet}(M))$ can be considered, but they do not give an \Linf lift and we cannot conclude that $\mathsf{a}$ is strict.
		
		\item There is even another modification on the construction. We can take homotopies $\mathsf{h}$ such that $\mathsf{h} \circ e \neq 0$. Those maps will not satisfy the generating form condition since $\iota_{\hat{\mathsf{a}} x}\omega = d \, \mathsf{h}(x) + \mathsf{h}(e (x))$ for all $x \in \, S^{\bullet}L$ (in particular $\iota_{\mathsf{a}_1 x}\omega = d \, \mathsf{h}(x) + \mathsf{h}(l_1 (x))$ for all $x \, \in L$). In this case the action is not required to be pre\=/hamiltonian.

	\end{enumerate}

But none of these arise in the case of non\=/graded Lie algebra acting on a pre\=/symplectic manifold.

\begin{pp}\label{ultimo}
Let $G$ be a Lie group with connectd Lie algebra $\mathfrak{g}$ acting by symplectomorphisms on the pre\=/symplectic manifold $(M, \omega)$. Then there is a one to one correspondence between:
\begin{enumerate}
	\item Co\=/momentum maps for the action.
	\item \mbox{$L_{\infty}$}\=/Momentum maps for the action.
	\item Strong\=/\mbox{$L_{\infty}$}\=/Momentum maps for the action.
	\item Homotopy Momentum Maps for the action.
	\item Homotopy Hamiltonian Momentum Maps for the action.
\end{enumerate}
\end{pp}

\dem
Since a Lie$[1]$\=/algebra has trivial $1$\=/bracket we are in the hypothesis of Corollary \ref{eeoo}. This shows the equivalence between $2$, $4$ and $5$. Since $\mathfrak{g}[-1]$ is concentrated in degree $\{1\}$ we have that strong \Linf actions, \Linf actions and symplectic actions are the same. In the first section of this chapter we saw that co\=/momentum maps could be understood as \Linf lifts of the action which shows the equivalence between $1$, $2$ and $3$.
\qed\\


\printindex

\end{document}